%% file: TransportationScienceTemplate_SSRN.tex
\documentclass[trsc,nonblindrev]{informs3_ssrn} % current default for manuscript submission

\OneAndAHalfSpacedXI % current default line spacing
\usepackage{subfigure}
\usepackage{tabularx}
\usepackage{graphicx}
\usepackage{algorithm}
\usepackage{multirow}
\usepackage{booktabs}
\usepackage{algpseudocode}
\usepackage{float}
\usepackage{enumitem}
\usepackage{lipsum} % for filler tex
\usepackage{array}
\usepackage{tikz}
\usepackage{caption}
\usepackage{tabularx}
\usepackage{mathrsfs}
\usepackage[utf8]{inputenc}
\newcolumntype{P}[1]{>{\centering\arraybackslash}p{#1}}
\newcolumntype{L}[1]{>{\raggedright\arraybackslash}p{#1}}

\usepackage{adjustbox}
\usepackage{natbib}
 \bibpunct[, ]{(}{)}{,}{a}{}{,}%
 \def\bibfont{\small}%
\TheoremsNumberedThrough     % Preferred (Theorem 1, Lemma 1, Theorem 2)
\EquationsNumberedThrough    % Default: (1), (2), ...
\begin{document}
%%%%%%%%%%%%%%%%

%\setcounter{page}{1}
%\DOI{10.1287/xxxx.0000.0000}%


%\RUNTITLE{}

\TITLE{Rural School Bus Routing and Scheduling}
%\TITLE{}

\ARTICLEAUTHORS{%
\AUTHOR{Prabhat Hegde, Vikrant Vaze}
\AFF{Thayer School of Engineering at Dartmouth College,\\ \EMAIL{Prabhat.Hegde.TH@dartmouth.edu},\EMAIL{Vikrant.S.Vaze@dartmouth.edu} \URL{}}
} % end of the block
\ABSTRACT{Long school bus rides adversely affect student performance and well-being. Rural school bus rides are particularly long, incentivizing parents to drive their children to school rather than to opt for the school bus. This in turn exacerbates the traffic congestion around schools, further compounding the problem of long bus rides, creating a vicious cycle. It also results in underutilized school buses and higher bus operating costs per rider. To address these challenges, this paper focuses on the design of rural school bus routes and schedules, a particularly challenging problem due to its unique operational complexities, including mixed loading and irregular road networks. We formalize a rural school bus routing and scheduling model that tackles these complexities while minimizing the {total commute time} of students. We develop an original road network-aware cluster-then-route heuristic that leverages our problem formulation to produce high-quality solutions and demonstrate its superior performance through extensive ablation studies. For real-world case studies, our approach outperforms status quo solutions by reducing {the commute times} of students by {22-25\%}. Our solutions also make the school bus more attractive, helping address both the underutilization of school buses and the prevalence of car commutes. Our routing and scheduling approach can improve school bus use by {14-15\%} and reduce car trips that induce congestion near schools by {10-13\%}. Many rural school districts share the operational characteristics modeled in this study, including long bus rides, high operational expenditures, mixed loading, and a high proportion of car-based school commutes, suggesting the broad applicability of our approach. Ultimately, by reducing student travel times, increasing school bus utilization, and alleviating congestion near schools, our approach enables rural school district planners to address transportation-related barriers to student performance and well-being.}
% \ABSTRACT{%
% Excessively long school bus rides can adversely affect student performance and well-being. These challenges are exacerbated by the common characteristics of rural school districts, including the prevalence of mixed loading, the irregular structure of road networks, and the drop-off-induced congestion of road traffic. This paper presents a new modeling framework and an original solution approach to produce high-quality rural school bus routes and schedules. We formulate a mixed-load school bus routing problem to minimize student travel time. Our solution approach consists of a road-network-aware cluster-then-route heuristic motivated by the lower-density rural road network. For real-world case study datasets, our approach outperforms previous benchmarks by reducing student travel time by 37-39\%, and demonstrates the potential that shorter bus trips have in increasing bus ridership and reducing road congestion.
% }%

%\KEYWORDS{deterministic inventory theory; infinite linear programming duality; 
%  existence of optimal policies; semi-Markov decision process; cyclic schedule}

% \HISTORY{}

\maketitle
%%%%%%%%%%%%%%%%%%%%%%%%%%%%%%%%%%%%%%%%%%%%%%%%%%%%%%%%%%%%%%%%%%%%%%

% Samples of sectioning (and labeling) in TRSC
% NOTE: (1) \section and \subsection do NOT end with a period
%       (2) \subsubsection and lower need end punctuation
%       (3) capitalization is as shown (title style).
%
%\section{Introduction.}\label{intro} %%1.
%\subsection{Duality and the Classical EOQ Problem.}\label{class-EOQ} %% 1.1.
%\subsection{Outline.}\label{outline1} %% 1.2.
%\subsubsection{Cyclic Schedules for the General Deterministic SMDP.}
%  \label{cyclic-schedules} %% 1.2.1
%\section{Problem Description.}\label{problemdescription} %% 2.

% Text of your paper here
\vspace{-6mm}

\section{Introduction} \label{sec1-Introduction}
\input{sources/1-Introduction}

\section{Relevant Literature} \label{sec2-RelevantLiterature}
\input{sources/2-LiteratureReview}
\section{Collaboration with School Districts} \label{sec3-Collaboration}
\input{sources/3-Collaboration}

\section{Model} 
\label{sec4-Model}
\input{sources/4-ModelDescription}
\section{Solution Approach} \label{sec5-SolutionApproach}
\input{sources/5-Algorithm}

\section{Computational Results} \label{sec6-ComputationalResults}
\input{sources/6-CaseStudies}

%\newpage
\section{Impacts of Mode Shift and Congestion Level} \label{sec7-Congestion}
\input{sources/7-Congestion}

\section{Conclusion} \label{sec8-Conclusion}
\input{sources/8-Conclusion}

% Acknowledgments here
\section*{Acknowledgment}
We acknowledge the support and collaboration of our partner school districts, whose practical problem statements and datasets served as the foundation for this research.
% We extend our gratitude to them for their valuable contributions and insights that have informed and enriched this project.
% Enter the text of acknowledgments here
%The authors would like to thank Jamie Teague and Karen Wright from the SAU-70 school district for their collaboration in helping us survey the parents of the school district, and for the data on the existing routes}% Leave this (end of acknowledgment)

% Appendix here
% Options are (1) APPENDIX (with or without general title) or 
%             (2) APPENDICES (if it has more than one unrelated sections)
% Outcomment the appropriate case if necessary
%
% \begin{APPENDIX}{<Title of the Appendix>}
% \end{APPENDIX}
%
%   or 
%
% \begin{APPENDICES}
% \section{<Title of Section A>}
% \section{<Title of Section B>}
% etc
% \end{APPENDICES}

% References here (outcomment the appropriate case) 

% CASE 1: BiBTeX used to constantly update the references 
%   (while the paper is being written).
\bibliographystyle{informs2014trsc} % outcomment this and next line in Case 1
\bibliography{SchoolBusRouting.bib} % if more than one, comma separated
\begin{APPENDICES}
\include{sources/9-SI}
\end{APPENDICES}
% CASE 2: BiBTeX used to generate mypaper.bbl (to be further fine tuned)
%\input{mypaper.bbl} % outcomment this line in Case 2

\end{document}

%% file: sources/1-Introduction.tex
Safe and reliable transportation to and from school is fundamental for educational access. The iconic yellow school bus serves as the backbone of this critical service in the United States, ensuring that distance and geography do not become barriers to learning opportunities. Almost half a million yellow school buses transport more than 25 million students ($\sim$55\% of the U.S. K-12 population) on a daily basis \citep{Burgoyne:2017:00}. The U.S. school bus fleet, which is more than twice the combined size of all other forms of U.S. mass transit (including bus, rail, and airline transportation), plays a central role in both the education and transportation sectors. However, despite this massive scale and importance, the operational decision-making and financing of school buses in the U.S. is done by individual school districts, making it challenging to develop unified approaches for school bus planning.

Almost half of all U.S. school districts, representing approximately 30\% of all schools and 20\% of all schoolchildren, are classified as rural \citep{Showalter:2019:00}. A rural service environment is different from an urban one due to lower population densities, longer bus routes, fewer schools per district, fewer riders per bus stop, more stops per route, more irregular road network structures, more homogeneous fleets, fewer or no one-way streets, fewer or no public transportation options to commute to school, structurally different road congestion patterns, and lower availability of alternative roads \citep{Tangiah:2008:00}. Urban school bus routing methods often focus on optimizing routes by minimizing the number of buses used within a dense, well-connected grid but do not adequately address the unique challenges posed by rural service environments. Rural districts require routing and scheduling strategies that account for longer distances, varied road conditions, and the need to balance operational efficiency with the well-being of students who face extended travel times. This paper aims to address these issues and makes the following contributions.
\begin{enumerate}[leftmargin=*]
    \item We develop a new mixed-integer linear optimization (MILO) model that, for a given number of school buses, minimizes the {total commute time} of students in a mixed-loading setting. A mixed-loading setting is one in which students from different schools can ride the same bus on the same trip. Our model is informed by analyses of data from an extensive survey of parents in our partner school district and broader evidence showing that many rural school districts share similar characteristics. {Our integrated model of the Rural School Bus Routing and Scheduling Problem simultaneously assigns students to stops, students to buses, and therefore stops to buses, while also determining the sequence of stops in each bus route. Bus schedules are inferred from the times when stops are served on a route, and student schedules are determined by their assigned stop, their assigned bus, and the corresponding bus schedule.}
    % Our integrated model formulation effectively addresses the three previously formalized subproblems of the school bus routing and scheduling problem: bus stop assignment, bus route generation, and school bus scheduling. %Importantly, our formulation includes key components of these subproblems, namely bus stop selection as a part of the assignment process and school bell time adjustment within scheduling. 
    \item We design an original cluster-then-route heuristic that solves the MILO model in a reasonable runtime. In the clustering step, { we develop a \textit{constrained $k$-means based clustering} approach} that assigns buses to service regions that contain candidate bus stops and students. In the routing step, for each service region, we propose an original \textit{Student Assignment, Stop Selection, Sequencing and Scheduling} (\textit{SAS4}) approach to determine the student-to-stop assignment, bus routes, bus schedules, {and student schedules}. 
    \item In addition to improving bus routes, we address the problems of underutilized buses and localized congestion caused by car commutes in the vicinity of schools by proposing an overall modeling framework that includes two interconnected procedures. The aforementioned cluster-then-route approach solves our MILO model to generate bus routes and schedules for a given set of school bus-using students and reduces the {total commute time} of those students. Shorter bus commutes can incentivize additional students to use the school buses, thus relieving congestion in the vicinity of schools by reducing car commutes, but also require re-optimization of school bus routes. Our iterative modeling and computational framework solves the resulting fixed-point problem that captures these interdependencies and generates bus routes, schedules, travel times, and travel mode choices that are in equilibrium with each other. 
    \item We perform extensive computational experiments for two school districts in rural New England to evaluate the tractability and performance of our modeling and computational framework. Our cluster-then-route approach solves the MILO model in reasonable runtimes. From a practical standpoint, our solutions lead to {22-25\%} savings in {total commute times} compared to the status quo in the two school districts. Moreover, we demonstrate the additional benefits due to the travel mode shift, driven by better routing and scheduling that reduce student commute times. Our approach can transport {14-15}\% more students while maintaining the same average { commute time} reduction per student. Furthermore, we show that reducing car commutes to school results in a {10-13}\% reduction in congestion in the schools' vicinity.
% (0.6 minutes decrease in District I, 0.1 minutes increase in District II)
\end{enumerate}
% \vspace{-6mm}
% \subsection{Distinctive features of rural school districts}\label{subsec:DistinctFeatures}

Rural school districts have varied socioeconomic and demographic compositions, which influence their operational planning \citep{Greenough:2015:00}. However, in terms of school bus travel, research has identified the following common operational challenges and policy concerns related to student well-being. 
\begin{itemize}[leftmargin=*]
    \item \textbf{Longer school bus commutes}: Rural school bus commutes are longer than their urban counterparts and often exceed 30 minutes. Students with longer school bus travel times tend to skip breakfast and report lower grades, poorer levels of fitness, fewer social activities, and worse study habits \citep{Lutz:2024:00, Gottfried:2021:00}. Longer commutes disproportionately burden students in lower income households and aggravate absenteeism \citep{Gottfried:2021:00}. In addition, rural school districts are likely to have more unpaved and hilly roads, are less likely to employ full-time bus superintendents, and experience more frequent inclement weather disruptions \citep{Howley:2001:00,Zars:1998:00,Spence:2000:00}. Exposure to these characteristics, in addition to longer commutes, compounds the various risks for students riding rural school buses \citep{Howley:2001:00,Zars:1998:00,Spence:2000:00}.

    \item \textbf{Higher operating expenses}: Bus operating costs per student are about 40\% higher in rural school districts than in urban ones. The marginal cost of accommodating an additional student in a rural school district is also about 2.5 times higher than in an urban school district. An additional mile traveled by bus in a rural school district costs about 10 times more than that in an urban one \citep{Ellegood:2024:00}. Despite these high costs, buses in rural school districts are often underutilized \citep{Stanley:2021:00}. This empirical evidence shows that the most cost-effective strategy for rural school districts is to plan better routes without increasing the number of school buses \citep{Ellegood:2024:00}.

    \item \textbf{Mixed loading}: Rural school buses are more likely to have mixed loading \citep{Ellegood:2015:00}. Rural school districts in the U.S. have undergone considerable consolidation in the past century, creating larger districts with fewer schools and more students per school \citep{Bard:2006:00}. Mixed loading is more beneficial for geographically larger school districts, where more bus stops are shared by students from different schools and where school locations are close together \citep{Ellegood:2015:00}, making it widely applicable in rural districts. Mixing elementary-school children with older ones can have adverse effects, including exposure to age-inappropriate content, bullying, and insufficient supervision of younger children in a combined age-group bus journey.

    \item \textbf{Higher proportion of car commutes}: Half of U.S. school-age children commute to school by car, and it is likely even higher in rural school districts due to greater car ownership \citep{Lidbe:2020:00}. When schools are close together, have similar start times, and a sparse road network around them, car drop-offs cause congestion near schools, creating major concerns, including \textit{(i)} disproportionate increases in the commute time of school buses that visit multiple schools; reducing the sleep time of children; \textit{(ii)} safety concerns while loading and unloading cars during rush hour; and \textit{(iii)} exposure of students to emissions. 
\end{itemize}

Thus, rural school bus routing and scheduling involves mixed loading with a fixed number of school buses and an emphasis on minimizing the {total commute times} of students, improving school bus utilization, and reducing congestion near schools. We propose a modeling framework to address these priorities.
Section \ref{sec2-RelevantLiterature} summarizes the literature relevant to our rural school bus routing and scheduling problem. { Section \ref{sec3-Collaboration} details the collaboration with our partner school districts.} Section \ref{sec4-Model} presents the notation used in this study and the mathematical model that minimizes the {total commute time of students}. Section \ref{sec5-SolutionApproach} describes a general solution approach for the routing and scheduling of mixed-loading school buses in rural school districts. Section \ref{sec6-ComputationalResults} reports the results of our solution approach for two real-world school districts. Section \ref{sec7-Congestion} proposes a modeling extension to capture the impact of altered routes on students' travel mode choices and localized road congestion in the vicinity of school locations. Section \ref{sec8-Conclusion} concludes the paper.

%% file: sources/2-LiteratureReview.tex
School bus routing and scheduling research has gained traction in the last few decades. Some recent studies focus on developing models and solution methods that are applicable to real-world instances \citep{Park:2010:00, Ellegood:2020:00}. These studies design bus routes and schedules that transport students from designated stops to schools and back, complying with regulatory requirements and operational constraints that are unique to each school district and its planning priorities. The problem can be decomposed into three subproblems: i) stop assignment (assignment of students to stops), ii) bus route generation (selecting stops and connecting them into routes), and iii) bus scheduling (assigning routes to buses and scheduling bus arrival/departure times at stops and schools). The design of effective school bus routes and schedules involves addressing one or more of these subproblems. Our operational context is that of a typical rural school district with mixed loading, where students experience long bus rides and buses are underutilized. These districts have irregular road networks and localized congestion due to car commutes. { We review the literature in four areas. We begin with multi-school mixed-loading school bus routing models, followed by solution approaches that are closely aligned with our instance of the rural school bus routing and scheduling problem; these represent the primary domains of our direct contribution. We then examine the literature on congestion in the context of school bus transportation and the relationship between student travel time and academic performance; this provides a broader context and suggests additional implications of our work.}
% The operational context of each school district is unique and might involve one or a combination of the formalized subproblems of the SBRP: Bus stop selection, Bus route generation, Bus route scheduling, School bell time adjustment, and strategic transportation policy issues \citep{Park:2010:00, Ellegood:2020:00}. Our operational context is that of a rural school district that has multiple schools that are served by a fleet of homogeneous buses that have mixed loading. 

\subsection{Multi-School Mixed-Loading School Bus Routing and Scheduling}
Although early work on mixed loading \citep{Bodin:1979:00, Chen:1990:00} was inspired by its frequent occurrence in rural areas (Figure \ref{Img:US5_NH_VTSchools}, as well as \citet{Ellegood:2015:00}), relatively little attention has been paid since then to school bus routing with mixed loading ($\sim$15\% of papers in \cite{Ellegood:2020:00}). Most studies focus on minimizing the number of buses \citep{Braca:1997:00,Park:2012:00,Ruiz:2015:00} or a composite cost function of the number of buses and the total distance traveled by the buses \citep{Caceres:2017:00, Lima:2016:00, Miranda:2018:00}, or a multi-objective framework with the number of buses as one objective \citep{Kang:2015:00, Mokhtari:2018:00}. However, in rural school districts, which typically have fewer schools and require a limited bus fleet, reducing the fleet by a single bus could considerably increase student {commute time} while challenging the maximum bus travel time and capacity constraints (Appendix \ref{Impact_of_fleet_size}). Some mixed-loading models minimize the total bus travel distance \citep{Campbell:2015:01, Ellegood:2015:00, Yao:2016:00}. Although reducing bus mileage with a fixed fleet helps lower operating costs in rural school districts \citep{Ellegood:2024:00}, it may not align with optimizing {bus commute} for students (as confirmed by the results in Section \ref{Subsec:Value_of_optimization}), nor does it address the problem of underutilization of buses. Furthermore, although a time loss function defined as the difference between bus and car travel times has been considered an objective \citep{Spada:2005:00}, it does not capture the actual burden on students, particularly those whose car travel times are already long.

We contribute a practically relevant objective applicable to rural school districts \textit{and} a solution approach to optimize it, jointly determining student-stop-bus assignments, route sequences, and implied schedules.

{\subsection{Solution Approaches to the Rural School Bus Routing and Scheduling Problem}
Solution approaches to the school bus routing and scheduling problem (SBRSP) are typically informed by the characteristics of the model, which, in turn, are shaped by the specific problem statements that vary from district to district. Our problem statement, inspired by our practical context (Section \ref{sec3-Collaboration}), applies to many rural school districts. It focuses on minimizing the total commute time of students, which requires explicit tracking of individual students. Furthermore, pre-assigning students to stops introduces suboptimality (see Section \ref{sec:BenefiitsAlgoDesign}), which motivates joint optimization of the student-to-stop assignment within the SBRSP. Even when students are already assigned to stops, obtaining exact solutions to multi-school mixed-loading SBRSP formulations is computationally intractable for real-world instances \citep{Caceres:2019:00,Spada:2005:00}. Thus, our integrated student-centric model, directly justified by practical problem needs, leads to a substantially more challenging optimization problem.

To address the inherent complexity of SBRSP, researchers commonly employ two-phase methods instead of insertion heuristics or savings-based heuristics \citep{Ellegood:2020:00}. The constrained vehicle routing literature suggests that, among the two-phase methods, cluster-first route-second methods outperform route-first cluster-second methods \citep{Cordeau:2007:00}. The routing phase addresses three components: stop selection, student-to-stop assignment, and connecting stops into routes. Sequential approaches handle them in different orders: Location-Allocation-Routing (LAR) \citep{Riera:2012:00, Riera:2013:00} first selects bus stops, followed by student assignment, and then routing; Allocation-Location-Routing (ALR) \citep{Chapleau:1985:00, Bowerman:1995:00} assigns students to potential stops to inform stop selection, followed by route generation; and Location-Routing-Allocation (LRA) \citep{Schittekat:2013:00} selects bus stops first, then plans routes, and finally assigns students. These sequential heuristics can compound suboptimality by fixing decisions early without considering their downstream effects on subsequent components.
% The spatial and network structure underlying SBRSP instances influences solution quality, yet existing approaches typically overlook road network topology, instead relying on simplified representations that consider only point-to-point Euclidean distances or travel times. A notable exception is \cite{Zeng:2019:00}, which explicitly leverages network structure but is tailored specifically for urban networks with regular grid-like connection patterns. Rural network topologies are fundamentally different, characterized by irregular, sparse connections that violate grid assumptions and exhibit tree-like substructures with highly variable connectivity where location-based clustering performs poorly. 
The complexity of multi-school mixed-loading SBRSP has led researchers to adopt evolutionary or trajectory-based metaheuristics, such as genetic algorithms \citep{Oluwadare:2018:00, Tangiah:2008:00, Kang:2015:00}, variable neighborhood search \citep{Bogl:2015:00, Lima:2016:00}, greedy random search \citep{Silva:2015:00, Lima:2016:00}, ant colony optimization \citep{Yao:2016:00, Mokhtari:2018:00}, iterated local search \citep{Lima:2016:00, Miranda:2018:00}, tabu search, and simulated annealing \citep{Spada:2005:00}. These approaches typically incorporate established heuristic components and operate independently of the underlying mathematical structure of the problem formulation.

Schools are incorporated into routes through three approaches: \textit{Insertion} approach adds schools anywhere along routes to minimize degradation of the objective \citep{Park:2012:00,Ruiz:2015:00, Lima:2016:00}; \textit{Sequential} approach visits all stops first and then schools in an optimized sequence \citep{Campbell:2015:01, Ellegood:2015:00}; \textit{Closed} approach begins and ends routes at the same school while visiting other schools just before returning to the starting school \citep{Silva:2015:00, Yao:2016:00}. These approaches require separate consideration of school nodes in solution design and do not optimize student pickups between school visits.

We contribute to the literature on solution methodology for SBRSP a new cluster-then-route matheuristic that addresses critical gaps in existing approaches for integrated student-centric models in rural, mixed-loading contexts. Our approach differentiates itself from the existing literature in four key ways. (1) It uses road network-aware clustering designed specifically for sparse, irregular rural topologies rather than relying on Euclidean distance-based partitioning. (2) It jointly optimizes stop selection, student assignment, and routing using an integrated mathematical formulation rather than sequential heuristics that compound suboptimality. (3) It directly leverages problem structure in a matheuristic rather than employing metaheuristics that operate independently of the underlying model. (4) Finally, it enables inter-school en route student pickups with optimized school sequencing decisions embedded within the solution approach itself, eliminating the need for separate consideration of school nodes during solution design.}

\subsection{Congestion Considerations in School Bus Transportation}
Traffic flows and congestion have been well studied \citep{May:1990:00,Daganzo:1997:00,Hoogendoorn:2001:00}. Researchers have also studied congestion near places of interest and special events in order to mitigate localized traffic build-up \citep{Beaudry:2010:00, Kwoczek:2014:00}. Research on public transit scheduling has shown that modeling the equilibrium between passenger flows and user demand improves overall system performance \citep{Pi:2019:00,Wei:2021:00}. There is widespread recognition that congestion negatively influences students and increases various costs of school bus transportation \citep{LaVigne:2007:00, Lu:2017:00}. However, the proposed solutions to this problem have been limited to school choice \citep{Wilson:2017:00} or school location decisions \citep{Lu:2017:00, Vitale:2019:00}. Although the school bus routing literature does not directly address traffic flow and congestion, there is research on stochastic link travel times \citep{Yan:2015:00, Sun:2018:00, Babei:2019:00} and uncertainty in journey times and ridership \citep{Levin:2016:00, Caceres:2017:00}. None of these models capture how congestion-induced increases in bus travel times affect the decision to opt for the school bus. This feature is especially important under mixed loading with identical school start times because school buses travel longer distances on congested links than cars. {We address localized, time-specific real-world congestion challenges near schools by introducing} a new, unexplored dimension to the SBRSP: the impact of congestion-dependent travel times on students' mode choices and bus ridership patterns.

\subsection{Student Travel Time, Sleep, and Academic Performance}
The American Academy of Pediatrics recommends a school start time of 8:30 AM or later for middle and high schoolers to achieve an optimal level of sleep \citep{AAP:2014:00}, aligning with the circadian rhythm of adolescents \citep{Crowley:2007:00}. Even a modest postponement of school start time can have a positive influence on students' sleep, mood, performance, and behavior \citep{Hissel:2018:00, Owens:2010:00}. While staggering school start times can help districts move toward these recommendations, there are social costs and non-monetary challenges associated with changing school start times \citep{Kirby:2011:00, Owens:2014:00} that make adhering to the guideline impractical for many school districts.  From a routing perspective, in geographically large school districts, picking up distant students on separate trips of the same bus incurs high costs of deadheading operations.
% Staggering school start times is a good way to adhere to this recommendation. Staggering also enables a reduction in bus transportation costs by reusing buses to individually serve schools with sufficiently different start times. Recent SBRSP research \citep{Bertsimas:2019:00, Bertsimas:2023:00, Delarue:2023:00} has co-optimized routes and school start times to aid this effort. However, only 18\% U.S. public schools start at or after 8:30 AM \citep{Wheaton:2015:00}. 

Shortening commute times is a possible improvement to the status quo in such districts. In addition to better alignment with the adolescent circadian rhythm, the reduction of sleep loss can result in improved academic performance \citep{Fredriksen:2004:00,Curcio:2006:00} and reduce morning stress. Shorter commutes can improve well-being by allowing pre-school socialization, a few additional minutes to ensure that students eat breakfast \citep{Taras:2005:00, Rampresaud:2005:00}, and flexibility for unexpected emergencies or delays \citep{Short:2015:00,Wahistrom:2022:00}. We propose a routing and scheduling solution that aims to capitalize on these benefits by minimizing {the total commute time for students.}

%% file: sources/3-Collaboration.tex
{This work originated from outreach by a member of the school board of a rural New England school district seeking to improve their transportation operations. In an email, they summarized their situation as follows: ``I'm on the school board and we have, as you may know, a widely distributed district and a traffic pattern that goes from free flowing in the outskirts to standstill in town at the times we need transportation. We also have constraints, for example, no kid should be on a bus for more than a certain number of minutes, and anyone within a certain distance of their home needs bus access." Following preliminary discussions, we forged a partnership with the school district and worked closely with their business administrators, who enumerated the operational constraints and challenges to help us understand the full scope of the problem.

Our partnership with the district involved a series of meetings that motivated a detailed examination of their existing operations. Understanding operational specifics, such as the geographical spread of the district, traffic patterns, and the possibility (or lack thereof) of changing school start times, enabled us to identify the key factors affecting routing decisions. We learned that the district conducts an annual transportation survey, and based on comments received in previous surveys, the district administrators knew that student travel time was the main concern among parents. As part of our collaboration, they sought our input on the survey design before administering it again, ensuring that the questions captured information that would inform our routing model. Determining which features to incorporate into our model was an iterative process: we met regularly with the school district administrators to refine our understanding of what was essential and subsequently presented our solution approach and results for their feedback. Furthermore, our partner district brought to our attention the fact that several neighboring districts shared similar features, highlighting the potential for a broader impact of our work. Subsequently, we succeeded in acquiring data from a second partner school district, also in rural New England, to demonstrate the generalizability of our modeling and computational approach and our practical insights.

The remainder of this section reflects the insights gained through this collaboration. Section \ref{Subsec:ContextRegulationsConstraints} draws on our discussions with district administrators to provide an overview of the operational context, the regulations that govern school transportation, and the resulting constraints on routing decisions. Section \ref{Subsec:SurveyandCounts} presents the findings of the parent survey, which we helped design, along with vehicle counts in front of schools that informed our choice of objective function and motivated efforts to address congestion. Section \ref{Subsec:BroaderApplicabilityofModeling} discusses the modeling considerations that arose from our partner district's setting but are not addressed in the existing literature, and their broader applicability in school districts with similar characteristics.}

\subsection{{Operational Context, Regulations, and Constraints}}\label{Subsec:ContextRegulationsConstraints}

{The district operates three schools (elementary, middle, and high) and manages a fleet with a fixed number of school buses. All schools in the district start at 8:00 AM and, for local logistical reasons, changing these start times is not an option the district would like to consider. To make the best use of resources given this constraint, the school district uses mixed loading, where buses serve students attending different schools on the same trip. The district knows which students require transportation and which schools they attend. For safety reasons, students cannot be transferred between buses. The district faces the challenge of designing bus routes that determine: (1) the assignment of students to stops and buses, (2) the sequence in which each bus visits stops and schools, and (3) the timing of all pick-ups and drop-offs. Existing routes were planned manually, leading to potentially significant suboptimalities that our optimization approach seeks to address.

The morning commute is particularly critical for the school district. Several students are driven to school by parents, creating congestion near schools. In a mixed-loading setting where all schools start at the same time, this congestion disproportionately impacts school bus travel. In addition, morning routes operate under strict time constraints, as schools require students to arrive before instruction begins. Therefore, the district prioritizes the problem of morning commutes. State regulations establish guidelines for maximum bus ride times. The district also faces practical limitations on the available driver shifts, requiring buses to complete their routes within a specified time window. Some buses are underutilized due to long commute times.

Candidate bus stops are located at fixed, known locations that must be regularly monitored to ensure compliance with 22 evaluation criteria to adhere to state regulations. It is critically important that stops are at locations where there is sufficient time to embark and disembark. Identifying new candidates for bus stop locations would require the district to invest resources to ensure that new stops meet these criteria; hence, it is beyond the scope of this work. The school district recognizes that having many stops slows buses down and desires bus-stop consolidation to reduce the total number of stops. This consolidation was already underway in the district---two routes with an excessive number of stops were recently manually shortened.

From a legal standpoint, school districts are not required to provide transportation to students living within walking distance of school, defined as 2 miles in most states. However, the district currently operates two buses with predominantly ``in-town" service that cater to these ``walkers," and there have been complaints from some parents about students living close to school yet experiencing very long bus rides. Although state law requires the district to provide a stop within 1.2 miles of the residence of any student living outside the 2 mile radius, this threshold is too permissive for our partner district. They prefer to restrict walking to approximately 5--6 minutes per student, corresponding to a walking distance limit of 0.3 miles.

These operational realities, regulations, and constraints form the foundation for our routing model.}

\subsection{Insights from Parent Survey and Vehicle Counts} \label{Subsec:SurveyandCounts}

{The school district surveyed parents to assess what drives decisions about school bus use. In a sample of 790 parents, 41\% regularly used the school bus, 23\% occasionally used it, and the remainder did not use it. Both parents of school bus commuters (children who regularly take the bus) and private commuters (those who sometimes or never use the bus) overwhelmingly agree that the duration of the journey is the most important factor affecting whether the child takes the bus (Figure \ref{Img:SchoolSurvey}). For parents of private commuters, other criteria were considerably less important, suggesting that shorter bus journeys encourage bus use.
% \vspace{-4mm}
\begin{figure}[h]\label{Img:SchoolSurveyandData}
    \centering
    \subfigure[{The importance of five factors in school commute decisions: duration, punctuality, student behavior, safety, and comfort, to the school bus commuters and private commuters.}]{
        \includegraphics[width=0.47\textwidth]{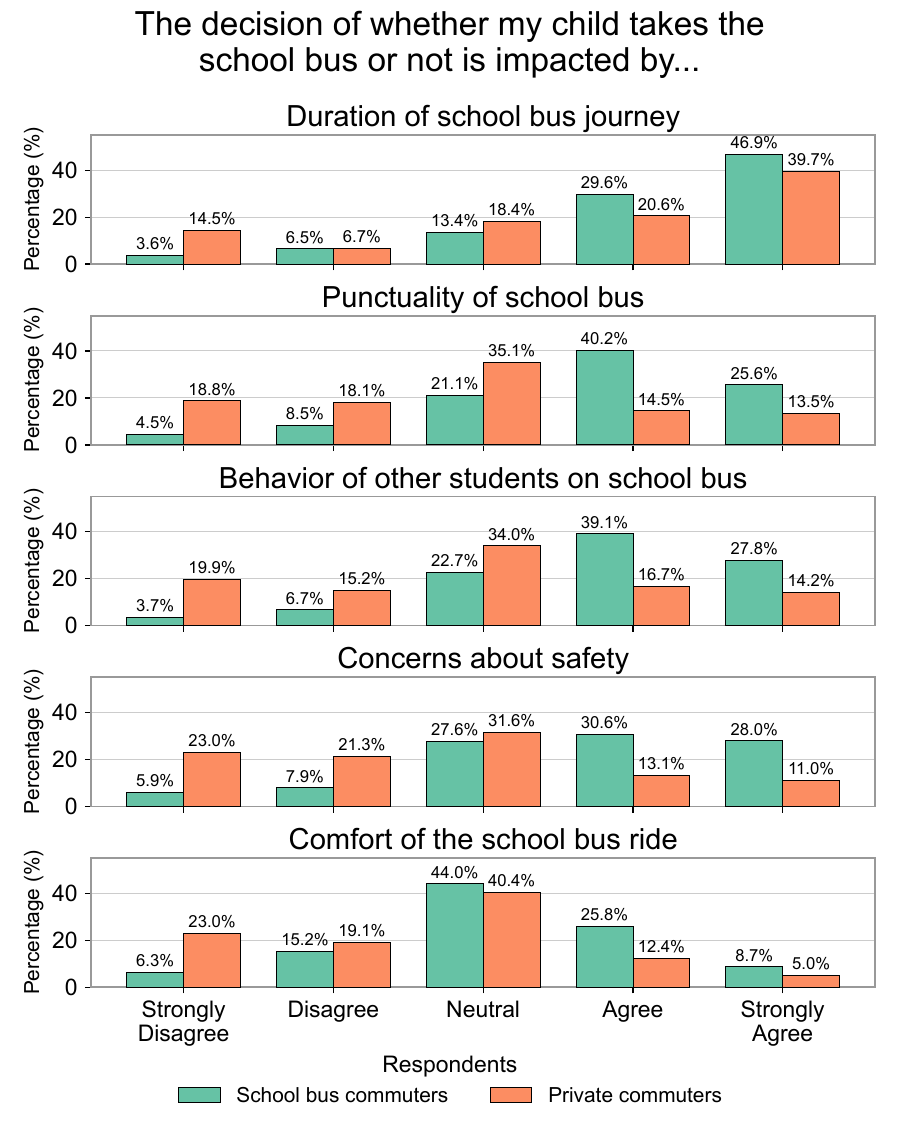}
        \label{Img:SchoolSurvey}
    }
    \hfill
    \subfigure[Number of vehicles on the streets that led to the Elementary School, Middle School, and High School on days when schools were open (solid lines) versus closed (dashed lines).]{
        \includegraphics[width=0.47\textwidth]{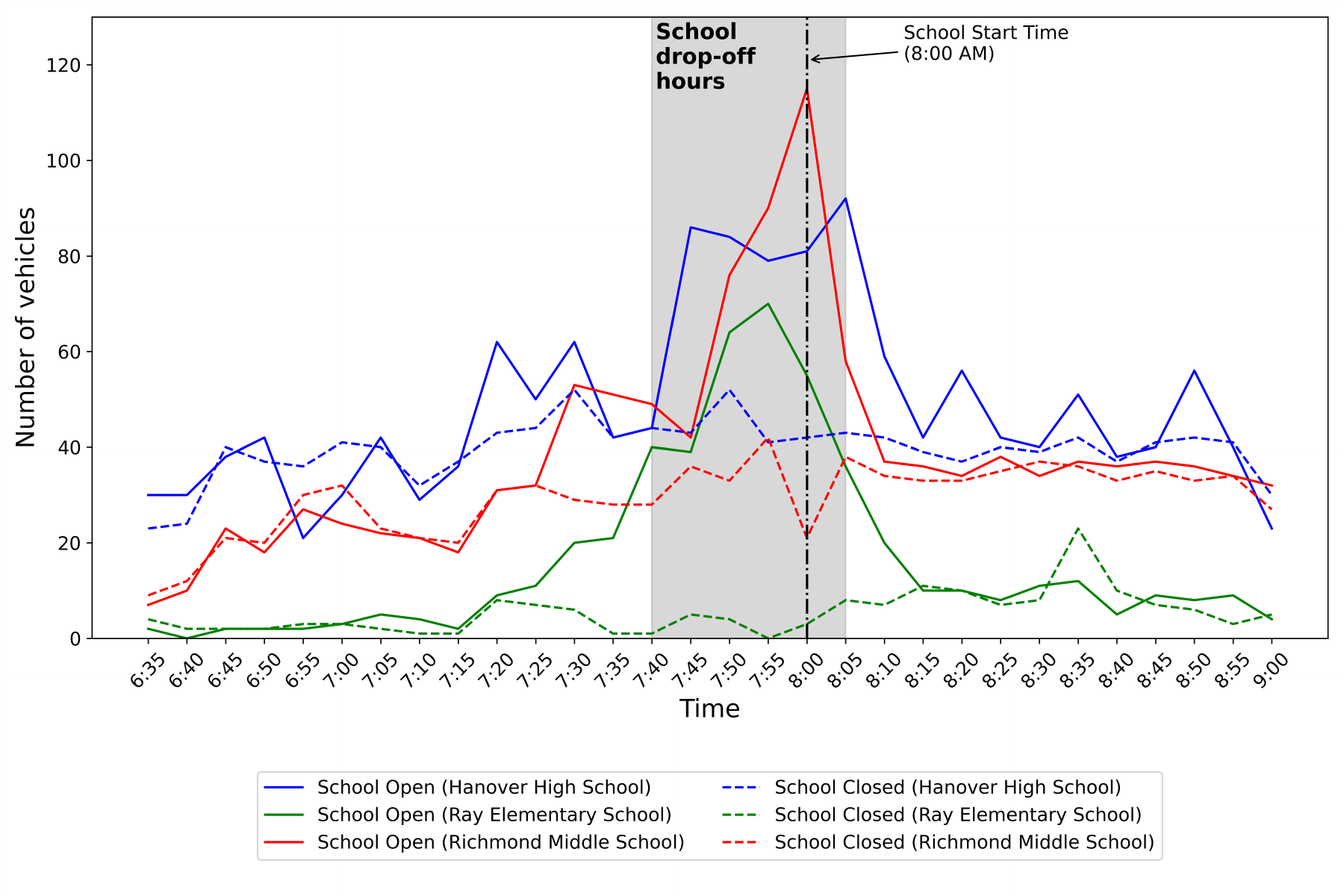}
        \label{Img:Congestion}
    }
    \caption{Insights from real-world survey of parents in a school district and vehicle counts data near schools}
    \label{fig:main}
\end{figure}
% \vspace{-7mm}
To assess congestion near schools, a concern raised by many parents, we manually collected vehicle counts (Figure \ref{Img:Congestion}) on nearby roads on two consecutive Mondays: one with school in session versus another when it was closed. In the 25-min window around school start time, we observed a large peak in vehicles near all three schools on the day school was open, while this peak was completely absent when school was closed. Neither day was a federal or state holiday, and all businesses remained open; the only difference was the operation of schools. These data suggest that reducing car commutes can alleviate congestion near rural schools and decrease travel times, an insight confirmed by computational results in Section \ref{sec7-Congestion}.}
%\vspace{-3mm}
\subsection{Broader Applicability and Contribution of Our Modeling Efforts}
\label{Subsec:BroaderApplicabilityofModeling}
{Mixed loading is particularly beneficial in rural service environments with a small number of schools, {especially when those schools are located in close proximity to one another} \citep{Ellegood:2015:00}. In addition, a smaller school count supports a smaller fixed school bus fleet. Figure \ref{Img:US5_NH_VTSchools} shows that most school districts in New Hampshire and Vermont are rural, usually with fewer than five schools, and that the schools are typically in close proximity to each other, which can contribute to localized congestion. These characteristics make most of these school districts well positioned to benefit from our SBRSP model, our solution approach, and our iterative framework to reduce localized congestion.}
\par\vspace{-3mm}
\noindent
\begin{minipage}[t]{0.48\textwidth}
  \vspace{0pt}
  \normalsize
  \normalfont
  Some school districts have considered staggering school start times in combination with designing bus routes and schedules \citep{Bertsimas:2019:00}, motivated by both operational considerations and health benefits. However, in practice, only 18\% of U.S. public schools start at or after 8:30 AM \citep{Wheaton:2015:00}. Many school districts (including our partners) face insurmountable logistical challenges that prevent the implementation of later or staggered start times \citep{Kirby:2011:00,Owens:2014:00}. Our model is informed by the situation in our partner school districts but is widely applicable to many other rural school districts
  (Figure \ref{Img:US5_NH_VTSchools}).
\end{minipage}%
\hfill
\begin{minipage}[t]{0.48\textwidth}
  \vspace{-8.7mm}
  \centering
  \includegraphics[width=\textwidth]{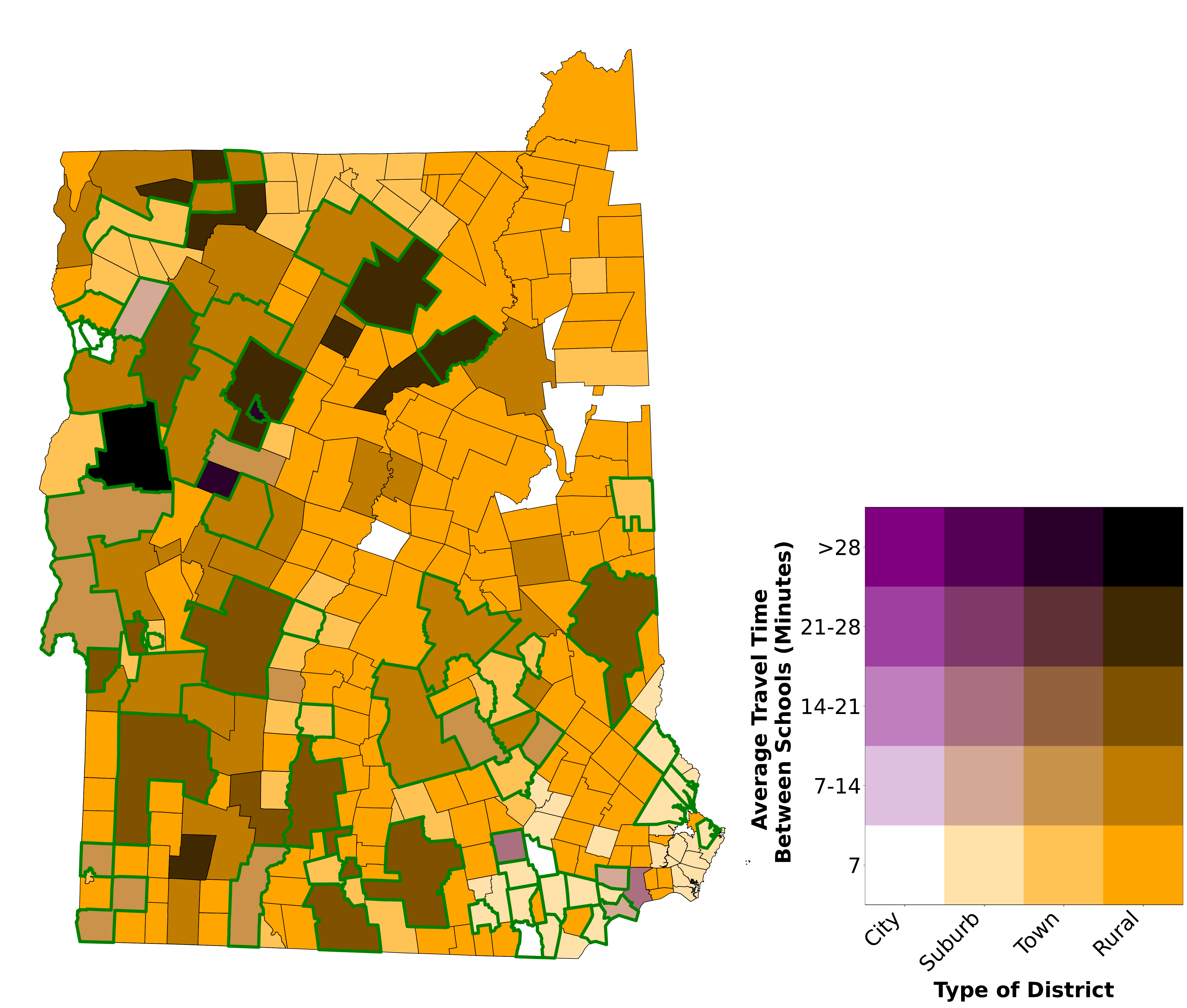}
  \captionof{figure}{New Hampshire and Vermont school districts,
  color-coded as city, suburb, town and rural as per the National
  Center for Educational Statistics classification
  \citep{Showalter:2019:00}. Shades show the average inter-school
  travel time. Districts outlined in green have more than five schools.}
  \label{Img:US5_NH_VTSchools}
\end{minipage}
\vspace{2mm}

{\hspace{3mm} Our original SBRSP formulation addresses several gaps in the existing scientific literature. In the existing school bus routing and scheduling models, the assignment of students to stops is often overlooked because students are pre-assigned to stops, and in some cases, stops are assigned exclusively to students attending the same school \citep{Spada:2005:00, Park:2012:00}. However, each additional stop could add a considerable amount of travel time for students already on the bus. We provide flexibility to the model for assigning students to stops and buses, which enables better solutions. Previous research also identifies assumptions for multi-school mixed-loading models: \textit{(i)} students from different schools can ride the same bus on the same trip; \textit{(ii)} all students are picked up by a bus before any drop-off from that bus; \textit{(iii)} each bus stop is served by only one bus; \textit{(iv)} each school can be served by multiple buses, although a bus can only visit a school once; \textit{(v)} a bus is not required to visit all schools \citep{Park:2012:00, Lima:2016:00}. Our conversations with rural school bus planners showed that assumptions \textit{(i)}, \textit{(iv)}, and \textit{(v)} are reasonable from an operational perspective, but \textit{(ii)} and \textit{(iii)} could lead to suboptimal solutions, particularly for stops geographically located between schools. Therefore, we relax assumptions \textit{(ii)} and \textit{(iii)}, which is of practical interest to many school districts where schools are in close proximity to each other and students can be picked up en route from one school to another.}

%% file: sources/4-ModelDescription.tex
Our mixed-integer linear optimization model for rural SBRSP with a fixed number of buses and mixed loading %Section \ref{Formulation} defines the notation and formulates the mathematical model.
is based on a network with nodes at the candidate bus stops ($\mathcal{N}$) and schools ($\mathcal{M}$), as well as the origin ($O$) and destination ($D$) depots. Directed arcs connect the origin to each stop ($\mathcal{A}_1$), each stop to every other stop ($\mathcal{A}_2$), each stop to each school ($\mathcal{A}_3$), each school to every other school ($\mathcal{A}_4$), each school to each stop ($\mathcal{A}_5$), and each school to the destination ($\mathcal{A}_6$). We denote by $\Delta_{ij}$ the bus travel time along arc $(i,j)$. The network is served by a fixed fleet of buses ($\mathcal{K}$). Appendix \ref{Impact_of_fleet_size} evaluates the impact of modifying the size of this fleet. Each bus $k$ has a capacity $C_k$, and must satisfy the route duration limit $T^B$, while each student’s bus ride time is bounded by $T^S$. $\mathcal{S}$ denotes the set of students taking the bus. $\mathcal{S}_i$ and $\mathcal{S}_m$ are the subsets of students who live within a short walking distance of stop $i$ and attend school $m$, respectively. $F_{is}$ denotes the time for student $s$ to walk to stop $i$. $G_{ms}$ is a binary parameter equal to 1 if student $s$ attends school $m$; 0 otherwise. $A_1, B_1$ and $A_2, B_2$ are coefficients of affine functions used to model pick-up and drop-off times \citep{Braca:1997:00}. $L_m$ is the number of students attending school $m$, and $\Lambda$ is the school bell time. 

The binary decision variables $e_{isk}$ assign students to stops and buses, $x_{ijk}$ determine whether a bus travels an arc, and $r_{ik}$ indicate whether a bus visits a node. {The integer decision variables $v_{mk}$ and $w_{ijk}$ represent the number of students dropped off at a given school by a given bus, and the number of students traveling on each bus on each arc, respectively. The continuous decision variables $t_{ik}$ denote the time at which each bus serves each node}. The continuous variables {$p_{s}$ and $d_{s}$ track the pick-up and drop-off times of each student}, while $\tau_{isk}$ and $\kappa_{msk}$ are their stop- and school-specific versions. Table \ref{Table:Notation} summarizes the notation. 

\textbf{Objective Function}
{\begin{equation}\label{Sec3:ObjectiveFunction}
    \text{Minimize} \sum_{s \in \mathcal{S}} \big(d_s - p_s\big) + \sum_{s \in \mathcal{S}} (\Lambda - d_s) + \sum_{i \in \mathcal{N}} \sum_{s \in \mathcal{S}_i} \sum_{k \in \mathcal{K}} F_{is} \cdot e_{isk}
\end{equation}
The objective function \eqref{Sec3:ObjectiveFunction} minimizes the total commute time of students, which is the sum of three components. The first term represents the total travel time of the students on the bus, measured as the difference between their drop-off and pick-up times. The second term captures the total idle time of the students waiting at school between their drop-off time and the start of school at time $\Lambda$, while the third term accounts for the total walking time of the students from their residences to their assigned bus stops. This student-centered aggregate objective follows the classical utilitarian principle \citep{Bertsimas:2011:00}. Since $\Lambda$ is a constant, the objective function \eqref{Sec3:ObjectiveFunction} can be simplified to the following equivalent form.}
\begin{subequations}
\renewcommand{\theequation}{\theparentequation.\arabic{equation}}
{\begin{equation}\label{Sec3:ObjectiveFunctionReduced}
    \text{Maximize} \sum_{s \in \mathcal{S}} p_s - \sum_{i \in \mathcal{N}} \sum_{s \in \mathcal{S}_i} \sum_{k \in \mathcal{K}} F_{is} \cdot e_{isk}
\end{equation}}

\begin{table}
\caption{Notation}
\label{Table:Notation}
\scriptsize
\begin{tabular}{>{\raggedright}p{2.3cm} >{\raggedright}p{1cm} p{12.2cm}}
\hline
\textbf{Component} & \textbf{Type} & \textbf{Description} \\
\hline
%\{d\}$ & Set & Destination \\
$\mathcal{N}$ & Set & Nodes representing bus stops \\
$\mathcal{M}$ & Set & Nodes representing schools\\
$\mathcal{S}$ & Set & Students\\
$\mathcal{K}$ & Set & School buses\\
$\mathcal{A}_1$ & Set & Zero-cost directed arcs connecting origin $O$ to bus stops. $\mathcal{A}_1 = \{(O,i): \forall i \in \mathcal{N}\}$\\
$\mathcal{A}_2$ & Set & Directed arcs connecting bus stops to other bus stops. $\mathcal{A}_2 = \{(i,j), \forall i,j \in \mathcal{N}, i \neq j \}$\\
$\mathcal{A}_3$ & Set & Directed arcs connecting bus stops to schools. $\mathcal{A}_3 = \{(i,j), 
    \forall i \in \mathcal{N},  \forall j \in \mathcal{M}\}$\\
$\mathcal{A}_4$ & Set & Directed arcs connecting schools to other schools. 
$\mathcal{A}_4 = \{(i,j), 
    \forall i,j \in \mathcal{M}, i\neq j\}$\\
$\mathcal{A}_5$ & Set & Directed arcs connecting schools to bus stops. $\mathcal{A}_5 = \{(i,j), \forall i \in \mathcal{M},  \forall j \in \mathcal{N}\}$\\
$\mathcal{A}_6$ & Set & Zero-cost directed arcs connecting schools to destination $D$. $\mathcal{A}_6 = \{(i,D), \forall i \in \mathcal{M} \}$ \\
$\mathcal{A}$ & Set & Arcs. $\mathcal{A} = \mathcal{A}_1 \cup \mathcal{A}_2 \cup \mathcal{A}_3 \cup \mathcal{A}_4 \cup \mathcal{A}_5 \cup \mathcal{A}_6$ \\
$\mathcal{N}'$ & Set & Nodes. $\mathcal{N}' = \{O\} \cup \mathcal{N} \cup \mathcal{M} \cup \{D\}$\\
$\mathcal{S}_i\subseteq \mathcal{S}$ & Set & Students that live within a maximum allowed walking distance to stop $i\in\mathcal{N}$\\
$\mathcal{S}_m\subseteq \mathcal{S}$ & Set & Students that attend school $m\in\mathcal{M}$\\
\hline
    $O, D \in \mathcal{N}'$ & Par. & Origin and destination depot nodes \\
    $A_1, A_2$ & Par. & Intercept of the linear expression for the time required to board and deboard students respectively\\
    $B_1, B_2$ & Par. & Slope of the linear expression for the time required to board and deboard students respectively\\
    $C_k$ & Par. & Student carrying capacity of bus $k \in \mathcal{K}$\\
    %$D$ & Par. & Maximum walking distance from home to bus stop\\
    % $T$ & Par. & Maximum travel time allowed for each bus route\\
    % $c_{ij}^S \in \mathbb{R}$ & Par. &Per-student cost of student time when traveling along arc $(i,j) \in A$ on the bus\\
    %$D_{is} \in \mathbb{R}^+$ & Par. & Distance from the home of student $s \in S$ to stop $i \in \mathcal{N}$\\
    $G_{ms} \in\{0,1\}$ & Par. & 1 if student $s \in S$ attends school $m \in \mathcal{M}$, 0 otherwise\\
    $Q_i \in \{-1,0,1\}$ & Par. & -1 if $i=D$, 1 if $i=O$, 0 otherwise, $\forall i \in \mathcal{N}'$\\
    $\Delta_{ij} \in \mathbb{R}^+$ & Par. & Bus travel time for arc $(i,j) \in \mathcal{A}$\\
    $L_m \in \mathbb{Z}^+$ & Par. & Number of students attending school $m \in \mathcal{M}$\\
    {$T^B$} & {Par.} & {Maximum time allowed for each bus route}\\
    {$T^S$} & {Par.} & {Maximum permissible student time on the bus}\\
    {$F_{is}$} & {Par.} & {Walking time for student $s \in \mathcal{S}_i$ to stop $i \in \mathcal{N}$}\\
    {$\Lambda$} & {Par.} & {Exogenous bell time of the schools}\\
\hline
    $e_{isk} \in\{0,1\} $ & Var. & 1 if student $s \in \mathcal{S}_i$ is picked up by bus $k \in \mathcal{K}$ at bus stop $i \in \mathcal{N}$, 0 otherwise\\
    $\tau_{isk} \in \mathbb{R}^+$  & Var. & Time at which student $s \in \mathcal{S}_i$ is picked up by bus $k \in \mathcal{K}$ at bus stop $i \in \mathcal{N}$ if it picks up that student at that bus stop, 0 otherwise\\
    $r_{ik} \in\{0,1\} $ & Var. & 1 if bus $k \in \mathcal{K}$ visits node $i \in \mathcal{N}'$, 0 otherwise\\
    $t_{ik} \in \mathbb{R}^+$ & Var. & Time when bus $k \in \mathcal{K}$ reaches node $i \in \mathcal{N}'$\\
    $x_{ijk} \in\{0,1\}$ & Var. & 1 if bus $k \in \mathcal{K}$ travels on arc $(i,j) \in \mathcal{A} $, 0 otherwise\\
    ${w_{ijk} \in \mathbb{Z}^+}$ & Var. & Number of students on bus $k\in \mathcal{K}$ when traveling on arc $(i,j)\in \mathcal{A}$ if it travels on arc $(i,j)$, 0 otherwise\\
    $\kappa_{msk} \in \mathbb{R}^+$  & Var. & Time at which student $s \in \mathcal{S}_m$ is dropped off by bus $k \in \mathcal{K}$ at school $m \in \mathcal{M}$ if it drops off that student at that school, 0 otherwise\\    
    ${v_{mk} \in \mathbb{Z}^+}$ & Var. & Number of students dropped off by bus $k \in \mathcal{K}$ at school $m \in \mathcal{M}$ \\
    % $u_{ik} \in \mathbb{R}^+$ & Var. & Number of students on bus $k \in \mathcal{K}$ right after the bus leaves after stopping at node $i \in \mathcal{N}'$\\
    {$p_{s}\in \mathbb{R}^{+}$} & {Var.} & {Time when student $s \in \mathcal{S}$ is picked up by a bus}\\
    {$d_{s}\in \mathbb{R}^{+}$} & {Var.} & {Time when student $s \in \mathcal{S}$ is dropped off}\\ 
\hline
\end{tabular}
\end{table}

\textbf{Origin and destination constraints}
\begin{align}
% &\sum_{i \in \mathcal{N}} x_{oik} = r_{ok} \quad \forall k \in \mathcal{K} \label{ODConstraint:EachBusCanHaveOnlyOneStartPoint}\\
% &\sum_{m \in \mathcal{M}} x_{mdk} = r_{dk} \quad \forall k \in \mathcal{K} \label{ODConstraint:EachBusEndsAtASchool}\\
%&t_{ik} \leq {T^B}\cdot (1-x_{Oik}) \quad \forall i \in \mathcal{N}, \forall k \in \mathcal{K}\label{ODConstraint:StartTime}\\
%&t_{Dk} \leq {T^B} \quad \forall k \in \mathcal{K} \label{ODConstraint:TravelTimeLimit}\\
&r_{Ok} = 1,  r_{Dk} =1, t_{Dk} = {T^B} \quad \forall k \in \mathcal{K}\label{ODConstraint:DepotStartTime}
\end{align}
% Constraints \eqref{ODConstraint:EachBusCanHaveOnlyOneStartPoint}
% ensure that each bus can have only one starting stop. Constraints \eqref{ODConstraint:EachBusEndsAtASchool} ensure that the bus can go to the destination from only one school.
%Constraints \eqref{ODConstraint:StartTime} ensure that time is counted from the first stop visited by a bus. 
%Constraints \eqref{ODConstraint:TravelTimeLimit} enforce an upper bound on the trip's end time. 
{Constraints \eqref{ODConstraint:DepotStartTime} require that all buses visit the origin and destination nodes and reach the destination node by time $T^B$. This effectively means that, in our model, we measure time starting from $\Lambda - T^B$ as time 0.}
%this sets to 0 We compute the actual pick-up and drop-off times for all students, as well as the arrival times of each bus at every stop and school, by setting the arrival time of each bus at its last visited school equal to the bell time $\Lambda$.

\textbf{Student pick-up and drop-off constraints}
\begin{align}
&\sum_{i \in \mathcal{N}}\sum_{s \in \mathcal{S}_i \cap \mathcal{S}_m} e_{isk} = v_{mk} \quad \forall m \in \mathcal{M}, k \in \mathcal{K}\label{PickUpConstraintNew:All_students_dropped_off_by_bus_are_picled_up}\\
% &\sum_{i \in \mathcal{N}} \sum_{s \in \mathcal{S}_i} e_{isk} = \sum_{m \in \mathcal{M}} v_{mk} \quad \forall k \in \mathcal{K} \label{PickUpConstraintNew:Students_Picked_Up_Have_To_Be_Dropped_Off}\\
&\sum_{k \in \mathcal{K}} v_{mk} = L_m \quad \forall m \in \mathcal{M}\label{PickUpConstraint:Students_Droppedoff_AtSchool_NumberOfStudentsThatGoToSchool}
\end{align}

Constraints \eqref{PickUpConstraintNew:All_students_dropped_off_by_bus_are_picled_up} ensure that each bus drops off a number of students at each school equal to the total number of students who attend that school and are picked up by that bus. 
%Constraints \eqref{PickUpConstraintNew:Students_Picked_Up_Have_To_Be_Dropped_Off} ensure that all students picked up by a bus are dropped off at the schools. 
Constraints \eqref{PickUpConstraint:Students_Droppedoff_AtSchool_NumberOfStudentsThatGoToSchool} ensure that the correct number of students are dropped off at each school.

\textbf{Bus ride constraints}
\begin{align}
&\sum_{j\in \mathcal{N} \cup \mathcal{M} \cup \{D\}:(i,j)\in \mathcal{A}} x_{ijk} - \sum_{j\in \mathcal{N} \cup \mathcal{M} \cup \{O\}:(j,i)\in \mathcal{A}} x_{jik} = Q_i \quad \forall i \in \mathcal{N}^\prime, k \in \mathcal{K} \label{FBConstraint:FlowBalanceConstraint}\\
&r_{ik} + r_{jk} \geq 2 \cdot x_{ijk} \quad  \forall (i,j) \in \mathcal{A}, k \in \mathcal{K}\label{BR:Constraint-Consecutive_stops}\\
&\sum_{j \in \mathcal{N}^\prime : (i,j) \in \mathcal{A}} x_{ijk} = r_{ik} \quad \forall i \in \mathcal{N} \cup \mathcal{M}, k \in \mathcal{K}\label{BR:Constarint-Bus_has_to_come_from_somewhere_to_pick_up_or_drop-off}\\
&r_{ik} \leq \sum_{s \in \mathcal{S}_{i}} e_{isk} \quad \forall i \in \mathcal{N}, k \in \mathcal{K}\label{BR:Constarint-Bus_stops_students_have_to_be_picked}\\
&e_{isk} \leq r_{ik} \quad \forall i \in \mathcal{N}, \forall s \in \mathcal{S}_i, \forall k \in \mathcal{K}\label{BR:Constarint-Bus_stops_students_have_to_be_picked_1}\\
% &r_{mk} = \sum_{(j,m) \in \mathcal{A}_3 \cup \mathcal{A}_4 } x_{jmk} \quad \forall m \in \mathcal{M}, k \in \mathcal{K}\label{BR:Constarint-Bus_has_to_come_from_somewhere_to_drop-off}\\
&\sum_{l\in\mathcal{N}':(l,i) \in \mathcal{A}_1 \cup \mathcal{A}_2 \cup \mathcal{A}_5} w_{lik}+ \sum_{s \in \mathcal{S}_i} e_{isk} - w_{ijk} \leq C_k\cdot (1-x_{ijk}) \quad \forall (i,j) \in \mathcal{A}_2 \cup \mathcal{A}_3 , k \in \mathcal{K}\label{BR:Constraint-Bus_visits_stop_students_before_after_positive}\\
&-\left(\sum_{l\in\mathcal{N}':(l,i) \in \mathcal{A}_1 \cup \mathcal{A}_2 \cup \mathcal{A}_5} w_{lik}+ \sum_{s \in \mathcal{S}_i} e_{isk} - w_{ijk}\right) \leq C_k\cdot (1-x_{ijk}) \quad \forall (i,j) \in \mathcal{A}_2 \cup \mathcal{A}_3, k \in \mathcal{K}\label{BR:Constraint-Bus_visits_stop_students_before_after_negative}\\
&\sum_{l\in\mathcal{N}':(l,i) \in \mathcal{A}_3 \cup \mathcal{A}_4} w_{lik} - v_{ik} - w_{ijk} \leq C_k\cdot (1-x_{ijk}) \quad \forall (i,j) \in \mathcal{A}_4 \cup \mathcal{A}_5 \cup \mathcal{A}_6, k \in \mathcal{K}\label{BR:Constraint-Bus_visits_school_students_before_after_positive}\\
&-\left(\sum_{l\in\mathcal{N}':(l,i) \in \mathcal{A}_3 \cup \mathcal{A}_4} w_{lik} - v_{ik} - w_{ijk}\right) \leq C_k\cdot (1-x_{ijk}) \quad \forall (i,j) \in \mathcal{A}_4 \cup \mathcal{A}_5 \cup \mathcal{A}_6, k \in \mathcal{K}\label{BR:Constraint-Bus_visits_school_students_before_after_negative}\\
&G_{ms} + \sum_{i\in\mathcal{N}: s\in\mathcal{S}_i} e_{isk} \leq r_{mk}+1 \quad \forall m \in \mathcal{M}, s \in \mathcal{S}, k \in \mathcal{K}\label{BR:Constraint-Bus_visits_school_if_student_is_on_board}\\
&\sum_{i \in \mathcal{N}: s \in \mathcal{S}_i} \sum_{k \in \mathcal{K}}e_{isk} = 1 \quad \forall s \in \mathcal{S}\label{BR:Student_unique_bus_unique_stop}\\
&t_{ik} + A_1 + B_1 \cdot \sum_{s \in \mathcal{S}_i}e_{isk} + \Delta_{ij} - t_{jk} \leq {T^B} \cdot (1-x_{ijk})\quad  \forall (i,j) \in \mathcal{A}_2 \cup \mathcal{A}_3, k \in \mathcal{K}\label{BR:Constraint-pick_up_time}\\
%&{t_{ik} + \Delta_{ij} + A_1 + B_1 * \sum_{s \in \mathcal{S}_j}e_{jsk} - t_{jk} \leq T \cdot (1-x_{ijk})\quad  \forall (i,j) \in \mathcal{A}_2 \cup \mathcal{A}_5, k \in \mathcal{K} \nonumber} \\
&t_{ik} + A_2 + B_2 \cdot v_{ik} + \Delta_{ij} - t_{jk} \leq {T^B} \cdot (1-x_{ijk})\quad  \forall (i,j) \in \mathcal{A}_4 \cup \mathcal{A}_5  \cup \mathcal{A}_6, k \in \mathcal{K}\label{BR:Constraint-drop_off_time}\\
%&{t_{ik} + \Delta_{ij} + A_1 + B_1 * v_{jk} - t_{jk} \leq T \cdot (1-x_{ijk})\quad  \forall (i,j) \in \mathcal{A}_3 \cup \mathcal{A}_4, k \in \mathcal{K} \nonumber} \\
& {d_s - p_s \leq T^S \quad \forall s \in \mathcal{S}} \label{TTConstraint:MaximumBusRideTime}\\
& w_{ijk} \leq C_k \cdot x_{ijk} \quad \forall (i,j) \in A, k \in \mathcal{K}\label{CapacityConstraint:NumberofStudentsInBusConsecutiveStops1}
\end{align}

Constraints \eqref{FBConstraint:FlowBalanceConstraint} balance the flow at each bus stop and school and ensure that each bus route starts at the origin and ends at the destination. Constraints \eqref{BR:Constraint-Consecutive_stops} ensure that a bus traveling on an arc visits both its end nodes. Constraints \eqref{BR:Constarint-Bus_has_to_come_from_somewhere_to_pick_up_or_drop-off} ensure that a bus visits a stop or a school if and only if it travels on exactly one outgoing arc from that stop or school. Constraints \eqref{BR:Constarint-Bus_stops_students_have_to_be_picked} and \eqref{BR:Constarint-Bus_stops_students_have_to_be_picked_1} require that at least one student is picked up by a bus at a stop if and only if that bus visits that stop. Constraints \eqref{BR:Constraint-Bus_visits_stop_students_before_after_positive} and \eqref{BR:Constraint-Bus_visits_stop_students_before_after_negative} (respectively, \eqref{BR:Constraint-Bus_visits_school_students_before_after_positive} and \eqref{BR:Constraint-Bus_visits_school_students_before_after_negative}) ensure that the number of students on the bus after visiting a stop (resp., school) is equal to the number of students on the bus before the stop (resp., school) was visited by that bus plus (resp., minus) the number of students picked up (resp., dropped off) at that stop (resp., school).
%Constraints \eqref{BR:Constraint-Bus_visits_school_students_before_after_positive} and \eqref{BR:Constraint-Bus_visits_school_students_before_after_negative} ensure that after visiting a school, the number of students on the bus is equal to the number of students on the bus before the school was visited minus the number of students dropped off at that school.
Constraints \eqref{BR:Constraint-Bus_visits_school_if_student_is_on_board} ensure that if a bus picks up a student who attends a particular school, then the bus must visit that school. Constraints \eqref{BR:Student_unique_bus_unique_stop} require that each student is picked up by exactly one bus at exactly one stop. Constraints \eqref{BR:Constraint-pick_up_time} and \eqref{BR:Constraint-drop_off_time} calculate the time at which a bus visits each stop and each school while accounting for the bus travel time and the time taken to pick up and drop off students. { Constraints \eqref{TTConstraint:MaximumBusRideTime} limit the bus ride time for each student to the maximum allowed value, ensuring a baseline level of horizontal equity \citep{Litman:2017:00} across all students.} Constraints \eqref{CapacityConstraint:NumberofStudentsInBusConsecutiveStops1} ensure that the number of students on a bus does not exceed the bus capacity and that students do not travel on an arc that is not traveled by that bus.

% \textbf{School bus capacity constraints}
% \begin{align}
% &\sum_{s \in \mathcal{S}_i}e_{isk} \leq C_k \quad \forall i \in \mathcal{N} ,\forall k \in \mathcal{K}\label{CapacityConstraint:StudentsPicked/Capacity_at_stop}\\
% &v_{mk} \leq C_k\cdot r_{mk} \quad \forall m \in \mathcal{M} , \forall k \in \mathcal{K} \label{CapacityConstraint:StudentsDropped/Capacity_at_stop}\\
% & w_{ijk} \leq C_k \cdot x_{ijk} \quad \forall (i,j) \in A, k \in \mathcal{K}\label{CapacityConstraint:NumberofStudentsInBusConsecutiveStops1}
% \end{align}

% Constraints \eqref{CapacityConstraint:StudentsPicked/Capacity_at_stop} - \eqref{CapacityConstraint:NumberofStudentsInBusConsecutiveStops1} enforce the bus capacity constraints on the number of students on the bus, the number of students picked up at a bus stop, the number of students dropped off at a school, and the number of students on the bus while traveling between two nodes, respectively.

\textbf{Student schedule constraints}
\begin{align}
&\tau_{isk} + {T^B} \cdot (1 - e_{isk}) \geq {{p_s}} \quad \forall i \in \mathcal{N}, s \in \mathcal{S}_i, k \in \mathcal{K}\label{StudentTimeConstraintNew:PickUpTimeGreater}\\
&\kappa_{msk} \leq {{d_s}} \quad \forall m \in \mathcal{M}, s \in \mathcal{S}_m, k \in \mathcal{K}\label{StudentTimeConstraint:Drop-offTimeLesser}\\
% &t_{mk} * (G_{ms} + q_{sk}-1) \leq d_{sk} \quad \forall m \in \mathcal{M}, \forall s \in \mathcal{S},\forall k \in \mathcal{K}\label{Constraint:Drop-offTimeLesser}\\
% &T \cdot \sum_{i \in \mathcal{N}: s \in \mathcal{S}_i}e_{isk} \geq {p_{sk}} \quad \forall s \in \mathcal{S}, k \in \mathcal{K}\label{StudentTimeConstraintNew:PickUpTimeUpperLimit}\\
% &T \cdot \sum_{i \in \mathcal{N}: s \in \mathcal{S}_i} e_{isk} \geq {d_{sk}} \quad \forall s \in \mathcal{S}, k \in \mathcal{K}\label{StudentTimeConstraintNew:Drop-offTimeUpperLimit}\\
&\tau_{isk} \leq {T^B} \cdot e_{isk} \quad \forall i \in \mathcal{N}, s \in \mathcal{S}_i, k \in \mathcal{K}\label{StudentTimeConstraintNew:tau0=>t0}\\
&\tau_{isk} \leq t_{ik}\quad \forall i \in \mathcal{N}, s \in \mathcal{S}_i, k \in \mathcal{K}\label{StudentTimeConstraint:TauCannotBeGreaterthanT}\\
&\tau_{isk} \geq t_{ik} - {T^B} \cdot (1 - e_{isk}) \quad \forall i \in \mathcal{N}, s \in \mathcal{S}_i, k \in \mathcal{K}\label{StudentTimeConstraintNew:Forcet=tau->p=1}\\
&\kappa_{msk} \leq {T^B} \cdot \sum_{i \in \mathcal{N}: s \in \mathcal{S}_i} e_{isk} \quad \forall m \in \mathcal{M}, s \in \mathcal{S}_m, k \in \mathcal{K}\label{StudentTimeConstraintNew:kappa0=>t0}\\
&\kappa_{msk} \leq t_{mk}\quad \forall m \in \mathcal{M}, s \in \mathcal{S}_m, k \in \mathcal{K}\label{StudentTimeConstraint:KappaCannotBeGreaterthanT}\\
&\kappa_{msk} \geq t_{mk} - {T^B} \cdot (1 - \sum_{i \in \mathcal{N}: s \in \mathcal{S}_i} e_{isk}) \quad \forall m \in \mathcal{M}, s \in \mathcal{S}_m, k \in \mathcal{K}\label{StudentTimeConstraintNew:Forcet=Kappa->q=1}\\
% &t_{mk} * (G_{ms} -1) + \kappa_{msk} \leq d_{sk} \quad \forall m \in \mathcal{M}, \forall s \in \mathcal{S}_m,\forall k \in \mathcal{K}\label{StudentTimeConstraint:Drop-offTimeLesser}\\
&{T^B}\cdot(G_{ms} + e_{isk} - 2) \leq (t_{mk} - t_{ik}) \quad \forall m \in \mathcal{M}, i \in \mathcal{N}, k \in \mathcal{K}, s \in \mathcal{S}_i\label{StudentTimeConstraintNew:pickupafterdropoff}
\end{align}

Constraints \eqref{StudentTimeConstraintNew:PickUpTimeGreater}-\eqref{StudentTimeConstraint:Drop-offTimeLesser} {link} the pick-up (respectively, drop-off) time of a student to the time when the corresponding bus visits the corresponding stop (resp. school). Constraints \eqref{StudentTimeConstraintNew:tau0=>t0}-\eqref{StudentTimeConstraintNew:Forcet=tau->p=1} (resp., \eqref{StudentTimeConstraintNew:kappa0=>t0}-\eqref{StudentTimeConstraintNew:Forcet=Kappa->q=1}) set the pick-up (resp., drop-off) time variables for each student to the corresponding bus stopping times, if the bus picks up (resp., drops off) that student at that stop (resp., school) and 0 otherwise. {Because constraints \eqref{StudentTimeConstraint:Drop-offTimeLesser} need not be binding at optimality, the actual drop-off time of student $s$, $d_s$, is calculated in post-processing as $\sum_{m \in \mathcal{M}} \sum_{k \in \mathcal{K}}$ $\kappa_{msk}$ using the $\kappa_{msk}$ values obtained at the optimal solution}. Constraints \eqref{StudentTimeConstraintNew:pickupafterdropoff} prevent a bus from picking up a student after that bus has already visited the school of that student. 
Finally, constraints \eqref{eq:DoD:1}-\eqref{eq:DoD:5} define the allowable ranges of the values of the decision variables.
\vspace{-10mm}

\begin{align}
&e_{isk} \in \{0,1\}, \tau_{isk} \in \mathbb{R}^+ \quad \forall s \in \mathcal{S}_i, i \in \mathcal{N}, k \in \mathcal{K} \label{eq:DoD:1}\\
& r_{ik} \in \{0,1\}, t_{ik} \in \mathbb{R}^+ \quad \forall i \in \mathcal{N}', k \in \mathcal{K}\label{eq:DoD:4}\\
&x_{ijk} \in \{0,1\}, {w_{ijk} \in \mathbb{Z}^+} \quad \forall (i, j) \in \mathcal{A}, k \in \mathcal{K} \label{eq:DoD:2}\\
&\kappa_{msk}, {p_s}, {d_s} \in \mathbb{R}^+, {v_{mk} \in \mathbb{Z}^+} \quad \forall m \in \mathcal{M}, s \in \mathcal{S}_m, k \in \mathcal{K}\label{eq:DoD:5}
\end{align} 
\end{subequations}

%% file: sources/5-Algorithm.tex
\begin{figure}[h!]
    \centering
    \includegraphics[width=0.85\textwidth]{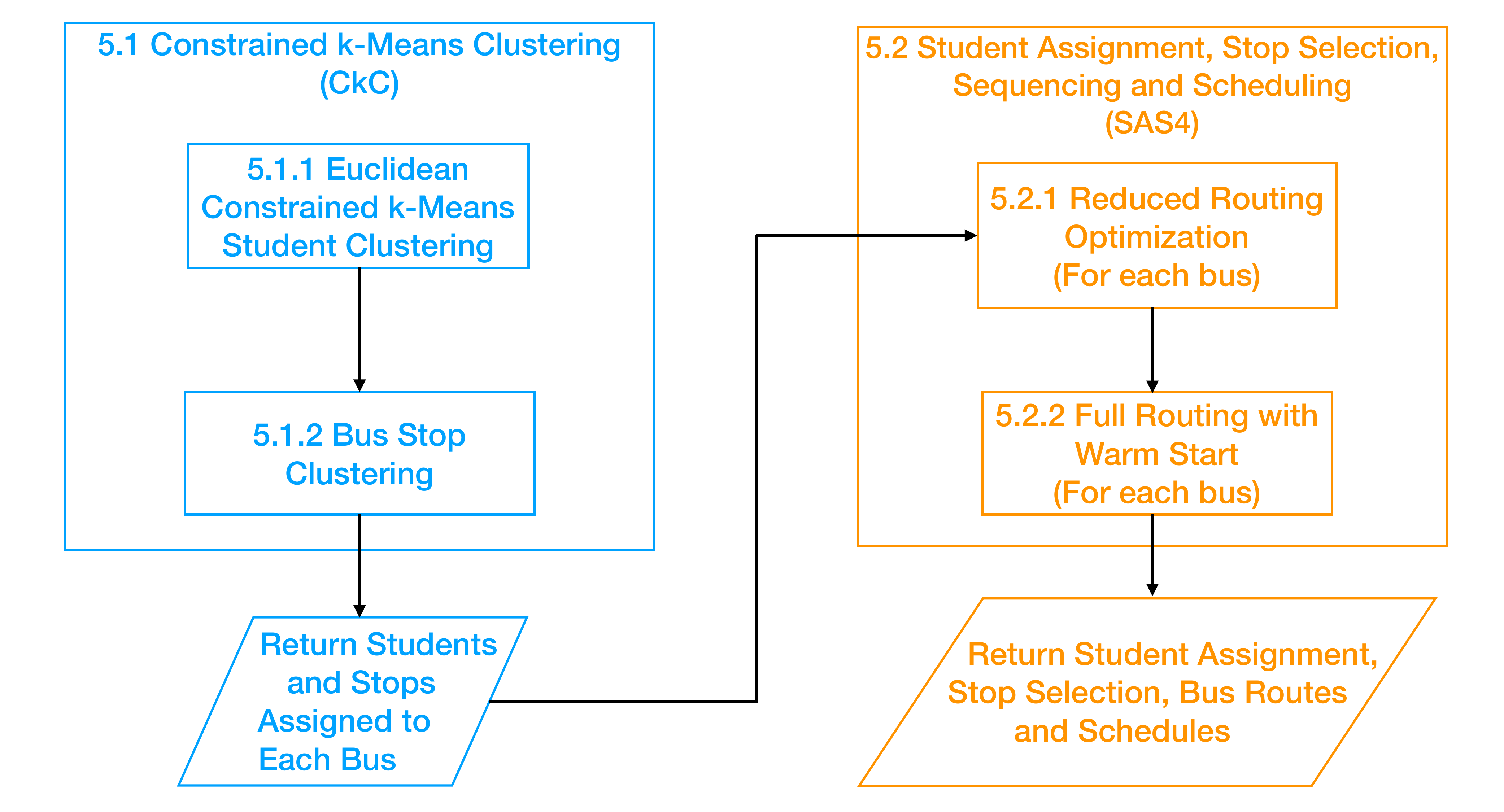}
    \caption{Our overall solution approach: A new two-phase cluster-then-route heuristic.}
    \label{fig:SBRPSolutionApproach}
\end{figure}

In this section, we present our solution approach (shown in Figure \ref{fig:SBRPSolutionApproach}) to the formulation \eqref{Sec3:ObjectiveFunctionReduced}-\eqref{eq:DoD:5} presented in Section \ref{sec4-Model}. Section \ref{subsection:constrained-$k$-means-main} details the clustering step, which introduces our {\textit{Constrained $k$-Means Clustering} (\textit{CkC})} approach, which assigns buses to service regions that contain stops and students. Section \ref{sec:S4A} details the routing step, where we explain our \textit{Student Assignment, Stop Selection, Sequencing and Scheduling (SAS4)} approach for designing routes of individual buses in their assigned service regions. %Our \textit{cluster-then-route} heuristic includes two sequential phases: a) {a \textit{Constrained $k$-means Clustering} (\textit{CkC})} phase (on left) creates and assigns service regions to individual buses (Section \ref{subsection:constrained-$k$-means-main}), and b) a \textit{Student Assignment, Stop Selection, Sequencing and Scheduling} (\textit{SAS4}) phase (on right) designs routes for each bus within its service region (Section \ref{sec:S4A}).

% \subsection{Overview}\label{subsec:Overview}
% \vspace{-5mm}

{ \subsection{Constrained \( k \)-Means Clustering} \label{subsection:constrained-$k$-means-main}
In this section, we present a Euclidean constrained $k$-means clustering approach to assign students and stops to buses. {This clustering approach relies on student locations and does not explicitly incorporate school locations. It is appropriate for school districts with the proximity characteristics discussed in Section \ref{sec3-Collaboration}, including our partner districts and many similar rural school districts, where schools are generally located in close proximity, making them amenable to mixed loading}. This approach is computationally tractable and avoids unbalanced clusters by using upper limits on cluster sizes set to bus capacities. Appendix \ref{appendix:RoadNetworkAwareClustering} describes an extended version of this clustering approach that incorporates road network information to further improve solution quality. Section \ref{subsection:constrained-$k$-means} describes the formulation for the Euclidean constrained $k$-means clustering of student locations. Section \ref{subsubsec:BusStopClustering} describes our approach to assign stops to student clusters.}

\subsubsection{Euclidean Constrained \( k \)-Means {Student} Clustering}
\label{subsection:constrained-$k$-means}
Consider a set $\mathcal{S}$ of students and a set $\mathcal{K}$ of buses. Let $C_k$ be the student carrying capacity of bus $k\in\mathcal{K}$. Let $(x_s, y_s)$ be the geographic location of student $s \in \mathcal{S}$.
%Let $\delta(x_i, c_k)$ represent the straight-line distance between student $i$ and the centroid of the service region of bus $k \in \mathcal{K}$.
Let $z_{sk}$ be a binary variable equal to 1 if and only if student $s$ is assigned to bus $k$. Let $(\bar{x}_k,\bar{y}_k)$ be the centroid of the locations of all the students assigned to bus $k$. The Euclidean constrained $k$-means clustering formulation is as follows \citep{Bradley:2000:00}.
\begin{subequations}
\renewcommand{\theequation}{\theparentequation.\arabic{equation}}
\begin{align}
    \text{Minimize} &\sum_{s \in \mathcal{S}} \sum_{k \in \mathcal{K}} z_{sk} \cdot ((x_s - \bar{x}_k)^2+(y_s - \bar{y}_k)^2) \label{eq:objective} \\
    \text{s.t.} \quad
    &\sum_{k \in \mathcal{K}} z_{sk} = 1 \quad \forall s \in \mathcal{S} \label{eq:assignment} \\
    &\sum_{s \in \mathcal{S}} z_{sk} \leq C_k \quad \forall k \in \mathcal{K} \label{eq:capacity} \\
    &\bar{x}_k = \frac{\sum_{s \in \mathcal{S}} z_{sk} x_s}{\sum_{s \in \mathcal{S}} z_{sk}} \quad \forall k \in \mathcal{K} \label{eq:centroidx}\\
    &\bar{y}_k = \frac{\sum_{s \in \mathcal{S}} z_{sk} y_s}{\sum_{s \in \mathcal{S}} z_{sk}} \quad \forall k \in \mathcal{K} \label{eq:centroidy}\\
    &z_{sk} \in \{0,1\} \quad \forall s \in \mathcal{S}, \forall k \in \mathcal{K}\label{Dod:Eucleadean$k$-means1}\\
    &\bar{x}_k, \bar{y}_k \in \mathbb{R}^{+} \quad \forall k \in \mathcal{K} \label{Dod:Eucleadean$k$-means2}
\end{align}
\end{subequations}
The objective function \eqref{eq:objective} minimizes the sum of squares of the distances between the location of each student and the centroid of the locations of all students assigned to the same bus as that student. The assignment constraints \eqref{eq:assignment} assign each student to exactly one bus. The constraints \eqref{eq:capacity} require that each bus $k$ be assigned to at most $C_k$ students to respect the bus capacity. The constraints \eqref{eq:centroidx}-\eqref{eq:centroidy} calculate the centroid of the locations of all the students assigned to a bus. Finally, the constraints \eqref{Dod:Eucleadean$k$-means1}-\eqref{Dod:Eucleadean$k$-means2} define the domains of the decision variables. {The initial centroid coordinates $\{\bar{x}^0_k, \bar{y}^0_k\}_{k \in \mathcal{K}}$ are selected using $k$-means++ initialization, which sequentially selects a subset of size $|\mathcal{K}|$ of student locations as initial centroids, where each randomized selection is based on probabilities favoring locations far from those already selected (see  \citet{Arthur:2006:00} for details).} This formulation is {then} solved iteratively as follows.
\begin{enumerate}
    \item \textbf{Assignment step}: At the start of each iteration $n$, fix the centroid coordinates \(\{\bar{x}^n_k,\bar{y}^n_k\}_{k \in \mathcal{K}}\) to their values from the previous iteration. Solve for the student-to-bus assignments \(z_{sk}^n\) that minimize the objective \eqref{eq:objective} subject to constraints \eqref{eq:assignment}, \eqref{eq:capacity}, and \eqref{Dod:Eucleadean$k$-means1}, keeping the centroids fixed.
    \item \textbf{Centroid update step}: Update the centroids using the assignments \(z_{sk}^n\) from the assignment step.
    \begin{subequations}
    \renewcommand{\theequation}{\theparentequation.\arabic{equation}}
    \begin{align}
    \bar{x}^{n+1}_k &= \frac{\sum_{s \in \mathcal{S}} z_{sk}^n x_s}{\sum_{s \in \mathcal{S}} z_{sk}^n} \quad \forall k \in \mathcal{K} \\
    \bar{y}^{n+1}_k &= \frac{\sum_{s \in \mathcal{S}} z_{sk}^n y_s}{\sum_{s \in \mathcal{S}} z_{sk}^n} \quad \forall k \in \mathcal{K}
    \end{align}
    \end{subequations}
    \item \textbf{Convergence check}: Iterate between the assignment and the centroid update steps until the centroid change satisfies the following convergence criterion (where \(\varepsilon\) is a small predefined tolerance value).
    \begin{equation}
    (\bar{x}^{n+1}_k - \bar{x}^n_k)^2 + (\bar{y}^{n+1}_k - \bar{y}^n_k)^2
     \leq \varepsilon \quad \forall k \in \mathcal{K}
    \end{equation}
\end{enumerate}
{ At each iteration, the assignment step (Step 1) optimizes the objective for fixed centroids, and the centroid update step (Step 2) optimizes the objective for fixed assignments. Consequently, the objective value is nonincreasing across iterations, and the procedure is guaranteed to terminate in a finite number of iterations at a solution that is locally optimal with respect to these two steps \citep{Bradley:2000:00}.}

Let $q_s=(x_s,y_s)$ represent the geographic location of each student $s \in \mathcal{S}_k$, where $\mathcal{S}_k$ is the set of students assigned to bus $k$. Let $\alpha_s$ be nonnegative coefficients summing to 1. Then, for each bus $k$, we define a service region $H_k$ as the convex hull of the locations of students assigned to that bus, as follows.
\begin{align} \label{eq:CH}
H_k = \left\{ q \in \mathbb{R}^2 : q = \sum_{s \in \mathcal{S}_k} \alpha_s \cdot q_s, \quad \sum_{s \in \mathcal{S}_k} \alpha_s = 1, \quad \alpha_s \geq 0 \:\:\:\: \forall s \in \mathcal{S}_k \right\} \quad \forall k \in \mathcal{K}.
\end{align}

\subsubsection{Bus Stop Clustering}\label{subsubsec:BusStopClustering}
According to our problem statement, each student must have access to a bus stop within the maximum allowed walking distance. Students who reside near the edges of the convex region $H_k$ may not satisfy this requirement; for some of these students, all bus stops within the maximum allowed walking distance may be outside $H_k$. To meet this requirement, the bus stop clustering step expands the service region $H_k$. Specifically, we shift all edges of the region $H_k$ outward from its center by a fixed distance equal to the maximum allowed walking distance to a bus stop, creating a new service region $H_k^{\prime}$.

Note that the overlapping parts of these expanded service regions may now contain some bus stops, overcoming a key restriction of previous studies that required each stop to be served by a single bus ({assumption \textit{(iii)} of common multi-school mixed-loading models in Section \ref{Subsec:BroaderApplicabilityofModeling}}). In addition, this approach is particularly beneficial for stops located between schools ({relaxing assumption \textit{(ii)} of common multi-school mixed-loading models mentioned in Section \ref{Subsec:BroaderApplicabilityofModeling}}, which requires all pick-ups to be done before any drop-offs by a bus begin), where students headed to schools in different directions can be picked up by different buses at the same stop -- a routing efficiency enabled by our overall modeling and solution approach.

\subsection{Student Assignment, Stop Selection, Sequencing and Scheduling (\textit{SAS4})} \label{sec:S4A}
Our SAS4 approach addresses the inherent computational challenge in our routing problem. As shown in Section \ref{sec:BenefiitsAlgoDesign}, the formulation \eqref{Sec3:ObjectiveFunctionReduced}-\eqref{eq:DoD:5} is intractable; it cannot directly produce feasible solutions, even for a single bus problem. This complexity stems from the need for explicit modeling of student-to-stop assignments rather than using pre-assignments. { We overcome this challenge by developing an original two-step routing approach guided by our formulation. We first reduce the complexity of our formulation by implementing three reductions and solve the resulting reduced routing formulation to optimality for a single bus, as described in Section \ref{Subsec:Reduced_Routing_Optimization}. The three reductions involve elimination of arcs $\tilde{\mathcal{A}}_5$ (which denote arcs $\mathcal{A}_5$ for a single bus formulation), modification of the objective function to consider only the travel time of students between stops and schools---excluding any boarding/deboarding times at stops or schools, walking time to stops, or idle time at schools---and pre-assignment of students to stops such that the number of stops is minimized. Then, in the second step, as described in Section \ref{Subsec:WarmStart}, we solve a single bus version of our full routing formulation by warm-starting it with the optimal solution of the reduced routing formulation.}

\subsubsection{Reduced Routing Optimization} \label{Subsec:Reduced_Routing_Optimization} We implement three reductions to the single bus version of our model (shown in Appendix \ref{SingleBusMILO}) to formulate a reduced problem and solve it optimally. The single bus version maintains a one-to-one equivalence with the formulation \eqref{Sec3:ObjectiveFunctionReduced}-\eqref{eq:DoD:5}, ensuring that our reduced routing formulation captures the same practical considerations. The three reductions are as follows.
\begin{enumerate}
    \item \underline{Reduction 1 - Remove network arcs}: We restrict the solution space to require the completion of all pick-ups before any drop-offs by that bus. This is accomplished by removing the arcs (denoted by $\tilde{\mathcal{A}_5}$ in Appendix \ref{SingleBusMILO} and $\mathcal{A}_5$ in formulation \eqref{Sec3:ObjectiveFunctionReduced}-\eqref{eq:DoD:5}) that connect schools to stops and by eliminating the now redundant constraints \eqref{StudentTimeConstraintNew:pickupafterdropoff_singleBus} of the single bus formulation in Appendix \ref{SingleBusMILO}.
    \item {\underline{Reduction 2 - Ignore non-travel times}: We modify the objective function \eqref{Sec3:ObjectiveFunctionReduced} of our original formulation to include only the travel time of the students between stops and schools, excluding any walking time, time spent waiting for boarding/deboarding, and idle time. The new objective is as follows.}
    \begin{align} \label{eq:ModifiedObjFun}  
        \sum_{(i,j) \in \tilde{\mathcal{A}} \setminus \tilde{\mathcal{A}_5}} \Delta_{ij} \cdot w_{ij}
    \end{align}
    \item \underline{Reduction 3 - Fix student assignments via stop minimization}: {As part of this reduction, we first obtain student assignments ($e_{is}$) and stop selections ($r_i$) by solving a stop-minimizing student pre-assignment problem.} Let $\tilde{\mathcal{N}}$ be the set of stops within the service region of a bus, $\tilde{\mathcal{S}}$ be the set of students assigned to that bus, and $\tilde{\mathcal{S}}_{i} \subseteq \tilde{\mathcal{S}}$ be the subset of students within a maximum allowed walking distance to stop $i \in \tilde{\mathcal{N}}$. Let $e_{is}$ be binary decision variables that indicate whether student $s \in \mathcal{S}_i$ is assigned to stop $i \in \tilde{\mathcal{N}}$. $r_{i}$ is a binary variable that indicates whether stop $i \in \tilde{\mathcal{N}}$ is selected. The stop minimizing student pre-assignment is formulated as follows.
    \begin{subequations}
    \renewcommand{\theequation}{\theparentequation.\arabic{equation}}    
    \begin{align}
        \text{Minimize} &\sum_{i \in \tilde{\mathcal{N}}} r_i \label{stopminimizinfstudentassignment:OF}\\
        \text{s.t.} \quad &\sum_{i \in \tilde{\mathcal{N}}: s \in \tilde{\mathcal{S}_i}} e_{is} = 1 \quad \forall s \in \tilde{\mathcal{S}} \\
        & e_{is} \leq r_i \quad \forall i \in \tilde{\mathcal{N}}, s \in \tilde{\mathcal{S}_i}\\
        &e_{is}, r_i \in \{0,1\} \quad \forall i \in \tilde{\mathcal{N}}, s \in \tilde{\mathcal{S}}_i\label{stopminimizinfstudentassignment:end}
    \end{align}
    \end{subequations}
    {We then fix student assignments ($e_{is}\: \forall i \in \tilde{\mathcal{N}}, s \in \tilde{\mathcal{S}}_i$) and stop selections ($r_i\: \forall i \in \tilde{\mathcal{N}}$) in the single bus reduced routing optimization formulation to those in the optimal solution to formulation \eqref{stopminimizinfstudentassignment:OF}-\eqref{stopminimizinfstudentassignment:end}.}
\end{enumerate}

\subsubsection{Full Routing with Warm Start}\label{Subsec:WarmStart}
A feasible solution to the reduced routing problem obtained by applying the three reductions described in Section \ref{Subsec:Reduced_Routing_Optimization} is clearly feasible for the single bus version of our original problem in Appendix \ref{SingleBusMILO}. We use the optimal solution to the reduced routing problem as an initial solution for each cluster to warm-start and solve a single bus MILO formulation in Appendix \ref{SingleBusMILO}.

%% file: sources/6-CaseStudies.tex
{We apply our solution approach to our partner school districts in rural New England. Appendix \ref{AppA} provides the case study details. This section assesses the performance of our solution approach and shows the value it provides to our partner school districts. The rest of this section is organized as follows.} Section \ref{sec:BenefiitsAlgoDesign} presents a systematic ablation study by selectively removing algorithmic features, allowing us to identify which algorithmic components are critical for obtaining high-quality solutions. Section \ref{Subsec:Value_of_optimization} shows the improvements that our solutions generate compared to the existing routes and schedules currently used by these school districts. Section \ref{subsec:StructureofRoutes} analyzes the structures and commute time distributions of existing and proposed solutions to identify key differences. Building on insights from Sections \ref{Subsec:Value_of_optimization} and \ref{subsec:StructureofRoutes}, Section \ref{Subsec:DriversofImprovement} systematically explores the main operational factors that drive the superior performance of our solutions. All experiments are performed on a server with 256 GB RAM and a Linux Operating System, implemented using Python 3.12.4 along with Gurobi Solver 11.0.3.

\subsection{Benefits of Algorithm Design} \label{sec:BenefiitsAlgoDesign}
Table \ref{AblationTable} shows the results of an ablation study comprising {six} algorithmic configurations. In the ``Base" configuration, all elements of our solution method shown in Figure \ref{fig:SBRPSolutionApproach} are active. In the other configurations, we disable key features of our methodology to evaluate their impacts on the quality of the solution and runtime. The tested algorithmic configurations (Config) are as follows. Furthermore, Appendix \ref{appendix:RoadNetworkAwareClustering} shows the additional benefits of incorporating road network-aware clustering as an extension of our solution approach. 
\begin{enumerate}
    % \item \textbf{No Euclidean constrained k-means \& No clustering problem size reduction (No \ref{subsection:constrained-$k$-means} and No \ref{subsubsec:ClusteringReduction}).} Directly solving the \textit{road network-aware constrained k-means} formulation \eqref{Eq:HybridkMeans_Objective}-\eqref{DoD:RNAwareKMeans2} in the clustering phase without constraining the initial assignments of the students through Equations \eqref{students_in_same_cluster_certain} and \eqref{student_in_a_given_cluster_certain}. 
    % \item \textbf{No road network-awareness (No \ref{subsubsec:ClusteringReduction} and No \ref{Subsec:RoadNetworkAwareConstrainedKMeans}).} Performing the clustering phase exclusively with the Euclidean constrained k-means formulation in Section \ref{subsection:constrained-$k$-means}. 
    % \item \textbf{No clustering problem size reduction (No \ref{subsubsec:ClusteringReduction}}). Using the complete solution approach, except for the clustering problem size reduction step and instead warm-starting \textit{road network-aware constrained k-means clustering} using the clusters from the Euclidean constrained k-means clustering results in Step \ref{subsection:constrained-$k$-means}.
    \item \textbf{No removal of arcs $\mathcal{A}_5$ (No Reduction 1) in \ref{Subsec:Reduced_Routing_Optimization}.} Using the complete solution approach, except for the elimination of arcs connecting schools to bus stops.
    \item \textbf{No modification to the routing objective (No Reduction 2) in \ref{Subsec:Reduced_Routing_Optimization}.} Using the complete solution approach, except for the alteration of the objective function shown in Equation \eqref{eq:ModifiedObjFun}.
    \item \textbf{No student-to-stop pre-assignment (No Reduction 3) in \ref{Subsec:Reduced_Routing_Optimization}.} Using the complete solution approach, with the exception of the stop-minimizing student pre-assignment.
    \item \textbf{No reduced routing, but with pre-assignment (No \ref{Subsec:Reduced_Routing_Optimization} but with formulation \eqref{stopminimizinfstudentassignment:OF}-\eqref{stopminimizinfstudentassignment:end}}). Using the complete solution approach, except Step \ref{Subsec:Reduced_Routing_Optimization}, and directly solving a single bus problem in Appendix \ref{SingleBusMILO} while fixing the pre-assignments based on the optimal solution to formulation \eqref{stopminimizinfstudentassignment:OF}-\eqref{stopminimizinfstudentassignment:end}.
    \item \textbf{No reduced routing (No \ref{Subsec:Reduced_Routing_Optimization}).} Using the complete solution approach, except Step  \ref{Subsec:Reduced_Routing_Optimization} and directly solving a single bus problem in Appendix \ref{SingleBusMILO} without any warm-start. 
\end{enumerate}
\begin{table}[h]
\centering
\caption{Ablation study results. ``No sol.'' indicates configurations where no solution was found for any of the clusters. ``No sol.$^*$'' indicates configurations where no solution was found for at least one of the clusters.}
\label{AblationTable}
\begin{tabular}{lrrrrrr}
\toprule
& \multirow{2}{*}{\begin{tabular}[c]{@{}r@{}}Feature(s)\\disabled\end{tabular}} 
& \multicolumn{2}{c}{District I} 
& \multicolumn{2}{c}{District II} \\
\cmidrule(lr){3-4} \cmidrule(lr){5-6}
Ablation & & Gap \% & Time (h) & Gap \% & Time (h) \\
\midrule
Base & - 
& 0 & 25.5 
& 0 & 19.8 \\

1. No removal of arcs $\mathcal{A}_5$ 
& \textit{\underline{rdn.1}} in \ref{Subsec:Reduced_Routing_Optimization} 
& 0.58 & 26.9 
& 0.62 & 20.9\\

2. No modification to the routing objective
& \textit{\underline{rdn.2}} in \ref{Subsec:Reduced_Routing_Optimization} 
& No sol.$^*$ & 29.1 
& No sol.$^*$ & 22.4 \\

3. No student-to-stop pre-assignment 
& \textit{\underline{rdn.3}} in \ref{Subsec:Reduced_Routing_Optimization} 
& No sol.$^*$ & 36.7 
& No sol.$^*$ & 28.6 \\

4. No reduced routing but with pre-assignment 
& \ref{Subsec:Reduced_Routing_Optimization} with Eqs \\ 
& \eqref{stopminimizinfstudentassignment:OF}-\eqref{stopminimizinfstudentassignment:end} 
& 27.38 & 25.3 
& 19.55 & 19.7 \\

5. No reduced routing 
& \ref{Subsec:Reduced_Routing_Optimization} 
& No sol. & 25.2 
& No sol. & 19.5 \\
\bottomrule
\end{tabular}
\end{table}
Gap \textit{\%} equals \text{$\frac{O-O^{*}}{O^{*}}$}, where $O$ and $O^{*}$ are the objective values obtained by a Config and the ``Base'' Config, respectively. The gaps quantify how much worse the solution quality of a Config is compared to that of the ``Base'' Config. Table \ref{AblationTable} reports results with a 5,000 s runtime limit for the reduced routing step (\ref{Subsec:Reduced_Routing_Optimization}) and a 10,000 s limit for the full routing step (\ref{Subsec:WarmStart}). The ``Base'' Config has runtimes of 19.8-25.5 hours.

{The removal of arcs $\mathcal{A}_5$ (Config 1) has a small (0.58\%-0.62\%) adverse impact on the quality of the solution. Moreover, the longer time taken to obtain solutions in Step \ref{Subsec:Reduced_Routing_Optimization} leads to 5.5\%-5.6\% higher total runtime for this configuration. Config 4 substantially degrades the solution quality, resulting in a gap of 19.55\%-27.38\%. In Configs 2 and 3, at least one cluster yields no solution in the routing phase of the algorithm, while trying to run the model without any reductions (Config 5) yields no feasible solution for any of the clusters. Thus, each algorithmic component provides considerable value, together enabling us to tractably use our MILO formulation to obtain high-quality solutions. Moreover, as shown in Appendix \ref{appendix:RoadNetworkAwareClustering}, incorporating the road network-aware clustering extension consistently yields a further improvement in terms of additional cost reductions of approximately 1.3\%--2.7\% compared to this ``Base'' configuration.}

% Moreover, as shown in Appendix \ref{appendix:RoadNetworkAwareClustering}, incorporating the road network-aware clustering extension can provide a further improvement of blah blah compared to the ``Base'' configuration.}

% The removal of arcs $\mathcal{A}_5$ (Config 4) has no impact on the quality of the solution obtained. However, the longer time taken to obtain solutions in Step \ref{Subsec:Reduced_Routing_Optimization} leads to 4.7\%-6.1\% higher total runtime for this configuration. Ignoring road network structure (Config 2) worsens solution quality by 0.99\%-2.54\%. Config 7 substantially degrades the solution quality, resulting in a gap of 19.27\%-27.04\%. Attempting to directly solve the road network-aware constrained $k$-means clustering (Config 1) proves to be infeasible, as it yields no feasible solution even after running the \textit{road network-aware constrained k-means} formulation (Equations \eqref{Eq:HybridkMeans_Objective}-\eqref{DoD:RNAwareKMeans2}) for 24 hr. In Config 3, 5, 6 and 8, at least one bus yields no solution in phase 2 of the algorithm. Thus, each algorithmic component provides considerable value, which justifies the associated complexity.

\subsection{Value of Optimization}\label{Subsec:Value_of_optimization}
Table \ref{Table:Metrics_baseline} lists aggregate performance metrics for our solutions and the status quo. Our solutions provide substantial { commute time savings in both districts. Average commute time is reduced by 21.87\% (from 33.55 to 26.21 min) in District I and 24.53\% (from 41.75 to 31.50 min) in District II. Furthermore, our solution reduces both the average ride times (from 23.87 to 16.82 min in District I and from 30.09 to 20.41 min in District II) and the average idle times prior to the start of the school (from 7.43 to 6.35 min in District I and from 9.41 to 7.88 min in District II).} Our solutions have a slight increase in walking time (0.8 min in District I and 0.96 min in District II) due to fewer bus stops (20.11\% in District I and 25.96\% in District II).
\vspace{-4mm}
{
\renewcommand{\arraystretch}{1.3}
\begin{table}[h]
\centering
\footnotesize
\caption{Aggregate Metrics Demonstrating the Value of Optimization.}
\label{Table:Metrics_baseline}
\begin{tabular}{>{\raggedright\arraybackslash}p{4.23cm}
                >{\raggedleft\arraybackslash}p{1.3cm}
                >{\raggedleft\arraybackslash}p{1.3cm}
                >{\raggedleft\arraybackslash}p{1.3cm}
                >{\raggedleft\arraybackslash}p{1.3cm}
                >{\raggedleft\arraybackslash}p{1.3cm}
                >{\raggedleft\arraybackslash}p{1.3cm}
                >{\raggedleft\arraybackslash}p{1.3cm}
                >{\raggedleft\arraybackslash}p{1.3cm}}
\toprule
& \multicolumn{4}{c}{District I} & \multicolumn{4}{c}{District II} \\
\cmidrule(lr){2-5} \cmidrule(lr){6-9}
Metric & Status Quo & Our Solution & Change & \% Change
       & Status Quo & Our Solution & Change & \% Change \\
\midrule

\shortstack[l]{Total commute time (min)}
 & 18,854.67 & 14,730.76 & -4,123.91 & -21.87
 & 18,120.29 & 13,674.62 & -4,445.67 & -24.53 \\

Avg commute time (min/student)
 & \small 33.55 & \small 26.21 & \small -7.34 &
 & \small 41.75 & \small 31.50 & \small -10.25 & \\

\addlinespace[0.5em]

\shortstack[l]{Total ride time (min)}
 & 13,413.89 & 9,450.33 & -3,963.56 & -29.55
 & 13,057.28 & 8,859.08 & -4,198.20 & -32.15 \\

Avg ride time (min/student)
 & \small 23.87 & \small 16.82 & \small -7.05 &
 & \small 30.09 & \small 20.41 & \small -9.68 & \\

\addlinespace[0.5em]

\shortstack[l]{Total idle time (min)}
 & 4,175.89 & 3,568.66 & -607.23 & -14.54
 & 4,082.37 & 3,419.22 & -663.15 & -16.24 \\

Avg idle time (min/student)
 & \small 7.43 & \small 6.35 & \small -1.08 &
 & \small 9.41 & \small 7.88 & \small -1.53 & \\

\addlinespace[0.5em]

\shortstack[l]{Total walking time (min)}
 & 1,264.91 & 1,711.77 & 446.86 & 35.33
 & 980.64 & 1,396.32 & 415.68 & 42.39 \\

Avg walking time (min/student)
 & \small 2.25 & \small 3.05 & \small 0.80 &
 & \small 2.26 & \small 3.22 & \small 0.96 & \\

\addlinespace[0.5em]

\shortstack[l]{No. of stops}
 & 179 & 143 & -36 & -20.11
 & 208 & 154 & -54 & -25.96 \\

\addlinespace[0.5em]

No. students per stop
 & 3.14 & 3.93 & 0.79 &
 & 2.09 & 2.82 & 0.73 & \\
\addlinespace[0.5em]

Total bus travel time (min)
& 403.4 & 367.92 & -35.48 & -8.8
& 422.50 & 367.99 & -54.51 & -12.9\\
Avg bus travel time (min)
& 44.82 & 40.88 & -3.94 & -8.8
& 60.36 & 52.57 & -7.79 & -12.9\\
\bottomrule
\end{tabular}
\end{table}
}

\begin{figure}
    \centering
    \includegraphics[width=\linewidth]{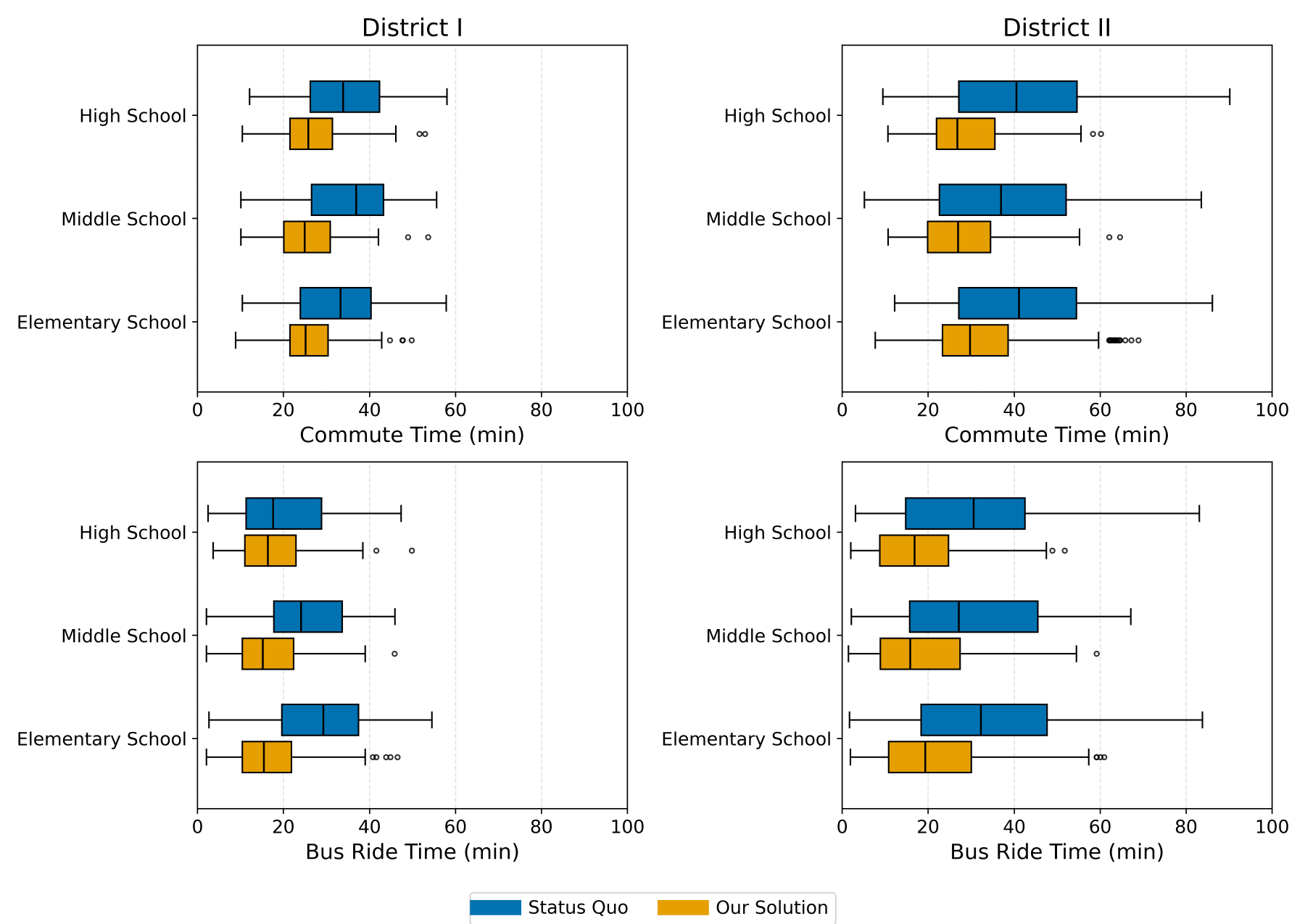}
    \caption{Total commute time and bus ride time comparison across schools and districts}
    \label{fig:Boxplot}
\end{figure}
\vspace{2mm}

\subsection{Structure of the Resulting Routes {and Commute Time Distributions}}\label{subsec:StructureofRoutes}

Section \ref{Subsec:Value_of_optimization} demonstrated that our solutions achieve significant improvements in {total commute times}. We now examine how these gains are distributed across students in Figure \ref{fig:Boxplot} in this section and in Table \ref{Table:CommuteTimePercentilesMinutes} in Appendix \ref{Section:Structure}. { Appendix \ref{Section:Structure} also presents the status quo and optimized routes for both districts in Figure~\ref{fig:routecomparison}, along with route-level operational metrics in Figure \ref{fig:two_plots}, including student bus ride time, student idle time, bus travel time, bus travel distance, student travel distance, and bus stopping time. In both districts, the optimized routes appear more consolidated and spatially coherent, with fewer fragmented segments and less overlap in service areas. The routes are more direct, with longer continuous segments and fewer detours, as evident in Figure~\ref{fig:routecomparison}. Figure~\ref{fig:two_plots} indicates that the total bus travel distances remain comparable, even as the bus route times decrease. It also shows a reduction in bus stopping times. Section \ref{Subsec:DriversofImprovement} leverages these observations to explain the key mechanisms driving improvements in our solution relative to the status quo.}

{Across the six schools in the two districts, Figure \ref{fig:Boxplot} shows substantial reductions in student commute times under our optimized routing solution, consistently ranging from 24\% to 34\%. Our solution reduces the median {commute time} for elementary school students by 24\% in District I (from {33.3 to 25.2 min}) and by 28\% in District II (from {41.1 to 29.8 min}). Middle school students experience median commute-time reductions of 32\% in District I (from {36.9 to 25.0 min}) and 27\% in District II (from {37.0 to 27.0 min}). High school students similarly benefit, with median decreases of 24\% in District I (from {33.9 to 25.8 min}) and 34\% in District II (from {40.5 to 26.8 min}). Importantly, the largest gains occur for students with the longest commutes. Elementary school commute times at the $95^\textit{th}$ percentile drop by almost 10 min in District I (from 50.4 to 40.8 min) compared to an 8 min decrease in the corresponding median values, and by 16 min in District II (from 78.8 to 62.9 min) compared to an 11 min decrease in the median. Middle school commute times at the $95^\textit{th}$ percentile drop by almost 13 min in District I (from 50.9 to 38.1 min) compared to a 12 min decrease in the median values, and by over 30 min in District II (from 76.3 to 45.9 min) compared to a 10 min decrease in the median; similarly, high school commute times at the $95^\textit{th}$ percentile drop by over 11 min in District I (from 51.8 to 40.6 min) compared to an 8 min decrease in the median, and by over 19 min in District II (from 71.7 to 52.4 min) compared to a less than 14 min decrease in the median.
% Middle school students see reductions of 25.0\% in District I and 39.8\% in District II, while high school students experience decreases of 21.7\% and 27.0\%, respectively.

Overall, the optimized solution delivers consistent commute-time improvements across districts and school levels, with particularly strong benefits for students facing the longest commute times.}

\subsection{{Key Dimensions and Mechanisms of Improvement Compared to the Status Quo}}\label{Subsec:DriversofImprovement}
{The status quo routes represent manually constructed plans developed by the school districts under the current operational practice (Section \ref{sec3-Collaboration}). In this section, we dive deep into three key complementary dimensions of improvement achieved by our optimized routes and identify two primary mechanisms that enable the improvements. Figure \ref{fig:DriversOfImprovement} summarizes these three complementary dimensions: students’ initial commute burden, the competitiveness of bus travel relative to direct commutes, and the existing pick-up positions.

The first major dimension is the substantial reduction in commute time for students with the longest commutes. Panel A of Figure \ref{fig:DriversOfImprovement} compares each student’s {final commute time under our optimized solution against their initial commute time under the status quo; points below the diagonal indicate reductions in commute time, while those above it indicate increases}. The results show a clear pattern: students with longer initial commutes experience substantially greater reductions. In particular, students with status quo commute times above the $75^\textit{th}$ percentile see average reductions of 22 min in District I and 37 min in District II, compared to only about 7 and 12 min, respectively, for the remaining students—indicating that the improvements are highly concentrated among those with the longest commutes. Overall, the majority of students (68\%-74\%) experience reductions in commute time, while only 26\%-32\% face modest increases.

A second dimension of improvement is the reduction in the excessive travel-time burden for students whose bus commutes are disproportionately long relative to a direct commute by car (measured as detour ratio, which is the ratio of a student’s total commute time using the bus to the time needed to drive to school). Whether a student chooses to ride the bus highly depends on the total commute time, as discussed in Section \ref{sec3-Collaboration} and modeled in Section \ref{sec7-Congestion}. Panel B of Figure \ref{fig:DriversOfImprovement} relates commute time improvements to the status quo detour ratio. The figure shows that students with the highest status quo detour ratios experience the largest reductions in commute time under our optimized solutions. In this way, the optimization primarily benefits the students for whom bus travel is the least competitive with a direct commute. Median improvements increase from 8–18\% (3–6 min on average) for detour ratios below 5, to 31–50\% (14–27 min) for detour ratios between 5 and 10, and to 38–64\% (19–35 min) for detour ratios above 10.
% on the order of {\color{blue}x\%} among the most extreme cases.

A third dimension of improvement is the especially large reduction in the cumulative delay experienced by students who board early and remain on the bus through many subsequent stops. As discussed in Section \ref{sec3-Collaboration}, even before we forged our partnership with the school districts, they had already identified reducing the number of stops as a practical goal to streamline operations and alleviate long ride times. Panel C of Figure \ref{fig:DriversOfImprovement} compares commute time improvements against a third dimension of improvements, namely students' status quo pick-up order. The results reveal that students who were previously picked up earliest, those in the first 25\% of stops along their assigned routes and therefore experiencing the greatest number of downstream stops, see the largest reductions in commute time under our solutions, with median improvements of 38–54\% (14–30 min). In contrast, students picked up near the end of routes, those in the final 25\%, may experience mild increases, typically less than 5 min, in commute time. Together, these three panels demonstrate that the improvements achieved by our optimization are not uniform but are driven by targeted reductions in the most burdensome components of the status quo system: long commutes, extreme detours relative to direct commutes, and cumulative delays borne disproportionately by early boarders.}

{More broadly, these results point to the potential for additional improvements through further reductions in stop frequency. Allowing longer maximum walking distances would increase flexibility in stop placement and could further reduce bus commute times. However, our primary results are not driven by such trade-offs. As shown in Table \ref{Table:Metrics_baseline}, while our approach modestly reduces the number of bus stops, the resulting increase in average walking time is less than one minute. Moreover, consistent with the preferences expressed by the partner school district (Section \ref{sec3-Collaboration}), the walking times remain well below state guidelines and the district's stated desires. To ensure that we are not systematically shifting the burden onto students, we constrain the average walking time for each bus to remain within a small tolerance of $\chi = 20$ seconds, as explored in Appendix \ref{Walking_time_driving_time}. In both districts, even with keeping average walking times nearly the same, the commute time consistently reduces by around 18\% (somewhat less than the original reductions of 22-25\% but still substantial) while the average walking time increase by a negligible 12-13 seconds across the two districts. Increasing the walking limit from 0.30 to 0.45 miles yields an additional reduction of approximately 3-4 minutes in average student commute time across both districts, while average walking time increases from about 3 min to roughly 4.5–4.7 min. The largest gains occur between 0.30 and 0.35 miles, where commute times decrease by about 1.5–2 min, whereas further increases from 0.40 to 0.45 miles yield only marginal improvements of less than one minute. Together, these findings show that our approach delivers substantial improvements even under conservative walking assumptions, while also identifying a clear pathway for further gains if these constraints are moderately relaxed. 
% Appendix \ref{fig:RideTimeWalkingTimeTradeOff} shows the results of this analysis.

We now describe the two key mechanisms of the improvements achieved by our solutions: reduced stopping time and higher effective driving speeds. We begin the discussion by observing an 8.8-12.9\% reduction in the travel time of the buses (as shown in Table \ref{Table:Metrics_baseline}) without a commensurate reduction in miles traveled. In District I, the average distance traveled per bus decreases by only 1.6\% (12.58 to 12.38 miles), while in District II it decreases by 4.5\% (17.60 to 16.81 miles), confirming that the substantial improvements do not arise from shorter distance routes. Instead, they are driven by more efficient operations. First, average stopping time per bus decreases by 20.8\% in District I (739 s to 585 s) and by 32\% in District II (923 s to 628 s), reflecting substantial reductions in start–stop delays due to stop consolidation. Second, average driving speeds increase by 7.5\% in District I (20.95 to 22.53 mph), indicating that buses are able to utilize faster roads, partly due to the reduction of ``in-town'' routes. In contrast, speeds increase only modestly in District II (by 2.0\%, from 23.48 to 23.95 mph), indicating that reduced start-stop delays is the dominant contributor of improvements. Together, these results show that our solutions achieve meaningful commute time reductions primarily through operational efficiencies rather than reductions in travel distance.

% To further disentangle these effects, Appendix \ref{Walking_time_driving_time} evaluates system performance under varying maximum walking distance thresholds. First, under tighter walking constraints, we confirm that the observed improvements are not simply driven by increased walking. Second, by relaxing these constraints, we quantify the additional gains achievable if longer walking distances were permitted. Together, these results show that our approach delivers meaningful improvements under conservative walking assumptions while also identifying the potential for further gains if these constraints are slightly relaxed.}
Finally, Appendix \ref{Impact_of_fleet_size} shows that although increasing fleet size could further reduce total commute time, the marginal gains may not justify the additional operational costs.}

\begin{figure}[H]
    \centering
\includegraphics[width=0.95\textwidth]{sources/Images/combined_scatter_plot_final_vs_initial.pdf}
    \caption{Final versus initial student commute times (Panel A) and commute-time improvement versus status quo detour ratio and status quo pick-up order (Panels B and C).}
    \label{fig:DriversOfImprovement}
\end{figure}

%% file: sources/7-Congestion.tex
So far, our model (Section \ref{sec4-Model}), solution approach (Section \ref{sec5-SolutionApproach}), and computational results (Section \ref{sec6-ComputationalResults}) have focused on minimizing the total commute times of students under the assumption that the students' decision to use a particular mode (bus or car) is not affected by improvements in bus commute times. However, in practice, shorter bus commutes are likely to make bus travel more attractive to students, especially in light of the empirical evidence in Section \ref{Subsec:SurveyandCounts} (Figure \ref{Img:SchoolSurvey}). As a result, some students may switch from car commutes to school buses, leading to secondary benefits from fewer vehicles on the road, reduced localized congestion near schools, and increased bus utilization. As more students adopt the bus, the roads in the vicinity of the schools become less congested, further reducing bus commute times. Car commutes also become faster. However, congestion near schools disproportionately affects school buses because they visit multiple schools. Thus, reduced congestion near schools provides a greater reduction in the disutility and travel times of school buses than car commutes and incentivizes a mode shift from car commutes to school buses. Quantifying the extent of these potential mode shifts is challenging due to the limited behavioral data. In practice, school districts typically revise bus routes on a seasonal basis, making an iterative approach that accounts for congestion and mode choice both feasible and beneficial. In order to capture these relationships and quantify the aforementioned secondary benefits, in this section, we present an iterative approach and a more comprehensive view of how optimized bus routes may influence the overall travel behavior of students.

\subsection{Algorithm}
We now describe our iterative algorithm to quantify the additional benefits of school bus routing and scheduling optimization. The algorithm captures three key relationships. First is the effect of changes in travel times on mode choice and demand for bus travel. {Section \ref{subsec:modechoice} describes the mode choice model and its calibration, while Section \ref{Subsec:ModeChoiceUpdate} explains how mode choices are updated within the iterative algorithm}. Second, it approximates the effects of changes in car commute patterns on road congestion and travel times (Section \ref{subsec:congestion}). Finally, it optimizes bus routes and schedules (Section \ref{sec5-SolutionApproach}) in response to changes in student demand and road travel times. Section \ref{subsec:overalliteralgo} describes the overall structure of this iterative algorithm. Note that students who regularly bike or walk to school are not considered in this analysis, as these transportation modes do not impact congestion or ridership and are typically used by students living closer to schools. 

\subsubsection{Mode Choice Model} \label{subsec:modechoice}
Consistent with the survey findings in Section \ref{Subsec:SurveyandCounts}, we divide the student population into three sets: $\mathcal{S}_A$ (students who ``always'' ride the bus), $\mathcal{S}_S$ (those who ``sometimes'' ride the bus), and $\mathcal{S}_N$ (those who ``never'' ride the bus). Our mode choice model, described by Equation \eqref{Eq:ModeChoice}, assumes that the students in the set $\mathcal{S}_S$ are the only ones who can consider changing their travel mode. The survey additionally suggests that some students in $\mathcal{S}_N$ could also switch to the bus if bus commute times decrease. In that sense, our results represent a lower bound on the benefits of our approach, and the actual benefits could be even larger. Let $T_{cs}$ and $T_{bs}$ be the anticipated commute times by car and bus, respectively, for student $s \in \mathcal{S}_S$. As per the multinomial logit model, the probability of $s$ choosing the bus is as follows. 
\begin{equation}\label{Eq:ModeChoice}
P_{bs} =
\frac{1}{1 + \exp\!\Bigl(-A\bigl(T_{cs} - T_{bs}\bigr)\Bigr)} \quad \forall s \in \mathcal{S}_S
\end{equation}
where $A \geq 0$ is the utility coefficient of {commute} time, indicating sensitivity to differences in {commute} time. The value of $A$ is not readily available, but it is key to approximating the impact of commute time on mode choice and, thus, to estimating the additional benefits of our approach. Hence, $A$ must be calibrated.

\textbf{Calibration of the utility coefficient of {commute} time ($A$)}: Let {$Y_{b}^{SQ}$} be the status quo number of bus riders, {and let $T_{cs}^{SQ}$ and $T_{bs}^{SQ}$ be the status quo car and bus commute
times, respectively, for student $s$. Let $P_{bs}^{SQ}$ denote the bus-choice probability obtained from Equation \eqref{Eq:ModeChoice} using these status quo commute times}. We calibrate $A$ so that the size of the always-bus group plus the expected number of bus riders in the $\mathcal{S}_S$ group matches the status quo number of bus riders {$Y_{b}^{SQ}$} under the status quo times, as follows.
\begin{align}\label{Eq:Calibration}
\lvert \mathcal{S}_A \rvert + \sum_{s \in \mathcal{S}_S} {P_{bs}^{SQ}} = {Y_{b}^{SQ}}
\end{align}

Substituting Equation \eqref{Eq:ModeChoice} into Equation \eqref{Eq:Calibration} and rearranging the terms, we obtain the following.
\begin{align}\label{Eq:TBD}
\sum_{s \in \mathcal{S}_S} \frac{1}{1 + \exp\!\Bigl(-A\bigl({T_{cs}^{SQ}} - {T_{bs}^{SQ}}\bigr)\Bigr)} = {Y_{b}^{SQ}} - \lvert \mathcal{S}_A \rvert
\end{align}

{This is a single equation in one variable, $A$, and it is solved using binary search to calibrate $A$ for each district. Using the calibrated value of $A$, we calculate $P_{bs}^{SQ}$ for each student in $\mathcal{S}_S$ and sort the students in descending order of $P_{bs}^{SQ}$. We define the cutoff probability $P^C=P_{bs'}^{SQ}$, where $s'$ is the student in the $(Y_b^{SQ}-\lvert\mathcal{S}_A\rvert)^\text{th}$ position in this order. Under the status quo commute times, students with $P_{bs}^{SQ}\geq P^C$ are classified as bus riders, and those below the cutoff as car commuters. This construction translates the observed aggregate ridership into a deterministic individual student-level mode assignment.}

\subsubsection{Road Congestion and Travel Time Update} \label{subsec:congestion}

Once we determine which students do not ride the bus, we calculate their shortest car paths from home to school. Similar to \citet{Wei:2021:00}, we draw on the speed-flow relationships from \citet{Schrank:2015:00} to estimate congestion and update road segment travel times. {Table~\ref{tab:speed_flow}, adapted from \citet{Schrank:2015:00}, presents the speed equations as a function of traffic volume per lane for freeways and arterial streets across different congestion levels.

{\begin{table}[htbp]
\centering
\caption{{Segment travel speed as function of traffic volume per lane (Schrank et al. 2015)}}\label{tab:speed_flow}
\begin{tabular}{llcc}
\toprule
Road Type & Congestion Level & Traffic Volume per Lane & Speed Equation \\
\midrule
Freeway & Uncongested & Under 15,000 & 60 \\
        & Medium & 15,001-17,500 & 70-(0.9$\times$ ADT/Lane) \\
        & Heavy & 17,501-20,000 & 78-(1.4$\times$ ADT/Lane) \\
        & Severe & 20,001-25,000 & 96-(2.3$\times$ ADT/ Lane) \\
        & Extreme & Over 25,000 & 76-(1.46$\times$ ADT/Lane) \\
\midrule
Arterial Street & Uncongested & Under 5,500 & 35 \\
                & Medium & 5,501-7,000 & 33.58-(0.74$\times$ ADT/Lane) \\
                & Heavy & 7,001-8,500 & 33.80-(0.77$\times$ ADT/Lane) \\
                & Severe & 8,501-10,000 & 31.65-(0.51$\times$ ADT/Lane) \\
                & Extreme & Over 10,000 & 32.57-(0.62$\times$ ADT/Lane) \\
\bottomrule
\multicolumn{4}{l}{{\small\textit{Note.} ADT/Lane stands for the average daily traffic per lane, in thousands.}} \\
\end{tabular}
\end{table}}

For roads with speed limits between 35 and 60 mph, we derive speed equations by linear interpolation of the coefficients in Table~\ref{tab:speed_flow}. We define an interpolation weight $w$ based on the speed limit $S$ as follows.
\begin{equation}
w = \frac{S - 35}{60 - 35} = \frac{S - 35}{25},
\end{equation}
where the denominator represents the difference between the freeway speed (60 mph) and the arterial street speed (35 mph) assumed by \citet{Schrank:2015:00}. Thus, $w = 0$ corresponds to arterial street conditions, and $w = 1$ corresponds to freeway conditions. For a speed equation given by $\text{Speed}=\alpha-(\beta \times \text{ADT/Lane})$ (where ADT stands for the average daily traffic), the interpolated coefficients are calculated as follows.
\begin{subequations}
\renewcommand{\theequation}{\theparentequation.\arabic{equation}}
\begin{align}
\alpha_{\text{interpolated}} &= \alpha_{\text{arterial}} + w \times (\alpha_{\text{freeway}} - \alpha_{\text{arterial}}), \\
\beta_{\text{interpolated}} &= \beta_{\text{arterial}} + w \times (\beta_{\text{freeway}} - \beta_{\text{arterial}}).
\end{align}
\end{subequations}

For all roads in the network, we first calculate the status quo ADT values based on the average speeds from OpenStreetMap \citep{Haklay:2008:00}. The roads in the immediate vicinity of schools bear the brunt of the school drop-off traffic, creating a peak traffic pattern near the morning drop-off times. Therefore, for these road segments near the schools, we further adjust the status quo ADT values to account for this peak traffic pattern using the traffic count data obtained in Section \ref{sec3-Collaboration}. In contrast, the remaining rural roads in the network have a relatively flat time-of-day profile and hence do not need such a peak traffic adjustment.

Within each travel time update in our overall iterative algorithm, we first update the path of each student commuting by car using the shortest path from the student's home to the student's school, leveraging Dijkstra's algorithm. We then update the road segment flows based on these updated car travel paths by summing all vehicles that use a segment. Next, we substitute the road segment flows into the speed equations with the aforementioned interpolated coefficients to calculate the speeds and travel times on each segment. Finally, these segment travel times are aggregated along each car's path to calculate the total car commute times for each student. In the next subsection, we explain our overall algorithm that iterates between the mode choice model, road congestion and travel time updates, and school bus routing and scheduling optimization.}

\subsubsection{Mode Choice Update}\label{Subsec:ModeChoiceUpdate}
{At each mode choice update, we recompute $P_{bs}$ for every student $s\in\mathcal{S}_S$ using Equation \eqref{Eq:ModeChoice} and the updated car and bus commute times, while holding $A$ and $P^C$ fixed at their calibrated values. We then assign students with $P_{bs}\geq P^C$ to the bus and those with $P_{bs}<P^C$ to the car, while students in $\mathcal{S}_A$ and $\mathcal{S}_N$ retain their original modes. Under this deterministic rule, changes in $P_{bs}$ affect ridership only when the bus-choice probability $P_{bs}$ of a student $s$ crosses the threshold $P^C$.}

\subsubsection{Overall Algorithm} \label{subsec:overalliteralgo}
\begin{figure}[h!]
    \centering
    \includegraphics[width=0.48\textwidth]{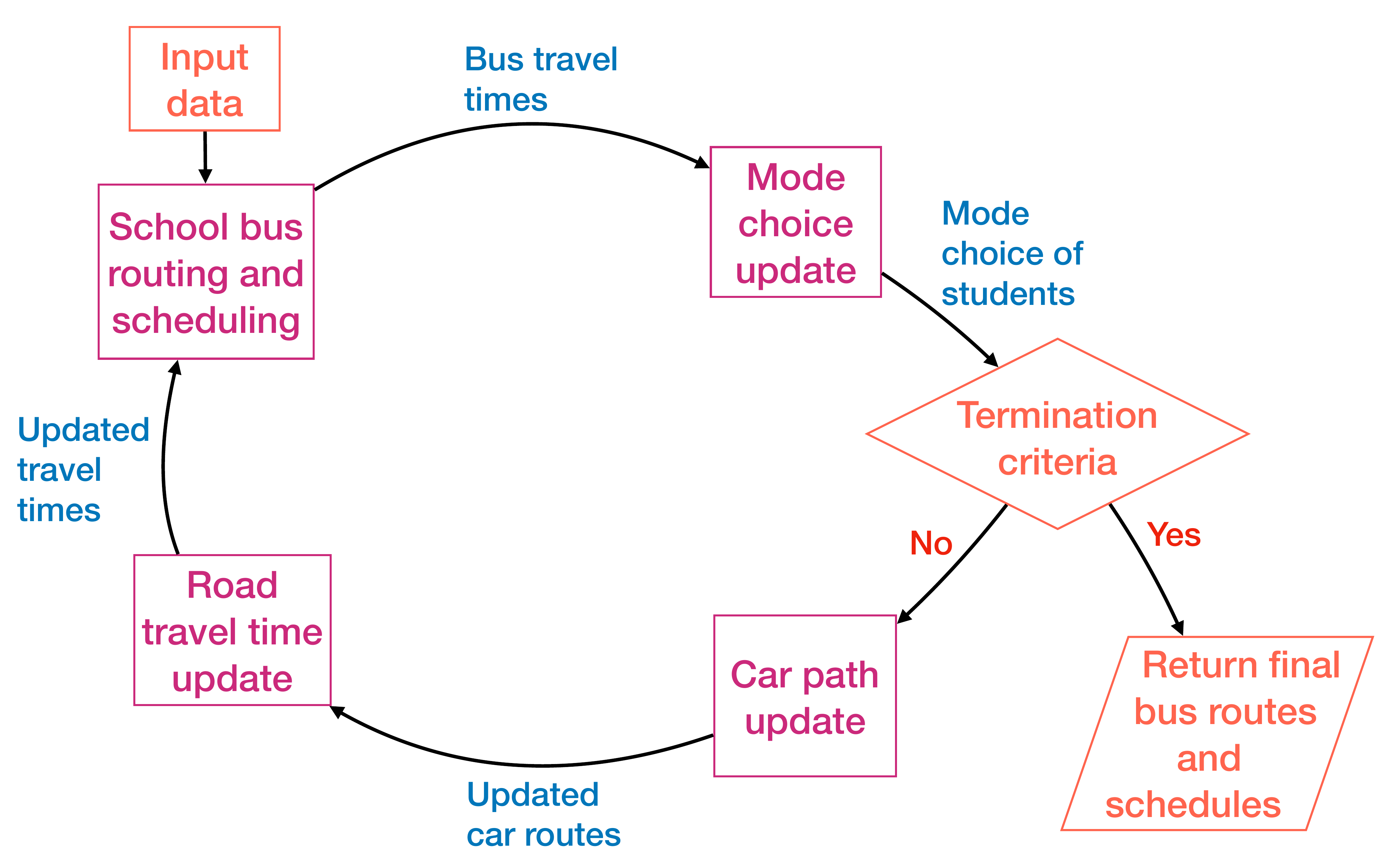}
    \caption{Iterative school bus routing and scheduling optimization with mode choice and travel time updating}
    \label{fig:SBRPFramework}
\end{figure}

{Figure \ref{fig:SBRPFramework} illustrates our overall iterative algorithm, which cycles through mode choice updates, car path updates, road travel time updates, and school bus routing and scheduling optimization. The algorithm begins with an initial set of inputs, including status quo travel times, which are calculated using the status quo ADT values and the speed-flow relationships for each road segment (as described in Section \ref{subsec:congestion}).
%student home locations, and school destinations,
First, the algorithm generates a baseline routing and scheduling plan by solving the SBRSP formulation  \eqref{Sec3:ObjectiveFunctionReduced}-\eqref{eq:DoD:5} presented in Section \ref{sec4-Model} using the solution approach described in Section \ref{sec5-SolutionApproach}. Next, we apply the mode choice update process described at the end of Section \ref{subsec:modechoice} to determine which students opt for the bus. Those who do not ride the bus commute by car and contribute to road congestion. We then calculate the paths of these car commutes based on the shortest travel time from home to the school that the student attends, update travel times to reflect this change in traffic flow as described in the last paragraph of Section \ref{subsec:congestion}, and re-run the bus route and schedule optimization, adjusting the bus routes and schedules. We repeat this cycle, updating mode choice, car travel paths, segment travel times, and bus routes and schedules until the number of students choosing the bus remains unchanged in successive iterations.}

\subsection{Results and Discussion}

\renewcommand{\arraystretch}{1.3}
\begin{table}[h]
\centering
\footnotesize % smaller font size for fitting width
\caption{Iterative Rural SBRSP with Mode Choice and Congestion Updating - Results}\label{Table:Metrics_congestion}
\begin{tabular}{>{\raggedright\arraybackslash}p{2.5cm}
                >{\centering\arraybackslash}p{1.1cm}
                >{\centering\arraybackslash}p{1.1cm}
                >{\centering\arraybackslash}p{1.1cm}
                >{\centering\arraybackslash}p{0.9cm}
                >{\centering\arraybackslash}p{0.9cm}
                >{\centering\arraybackslash}p{1.1cm}
                >{\centering\arraybackslash}p{1.1cm}
                >{\centering\arraybackslash}p{1.1cm}
                >{\centering\arraybackslash}p{0.9cm}
                >{\centering\arraybackslash}p{0.9cm}}
\toprule
& \multicolumn{5}{c}{District I} & \multicolumn{5}{c}{District II} \\
\cmidrule(lr){2-6} \cmidrule(lr){7-11}
Metric & Status Quo & Routing & Routing + Mode & $\Delta R$ (\%) & $\Delta RM$ (\%) & Status Quo & Routing & Routing + Mode & $\Delta R$ (\%) & $\Delta RM$ (\%) \\
\midrule
\shortstack[l]{Total commute\\ time (min)}
    & 18,854.67 & 14,730.76 &  17,138.38& -21.87 & -9.10 & 18,120.29 & 13,674.62 & 15,773.42 & -24.53 &  -12.95\\
\scriptsize \textit{Avg per student (min)}
    & \scriptsize 33.55 & \scriptsize 26.21 & \scriptsize 26.53 & \scriptsize -21.87 & \scriptsize -20.92 & \scriptsize 41.75 & \scriptsize 31.50 & \scriptsize 31.87  & \scriptsize -24.53 & \scriptsize -23.66 \\
\addlinespace[0.5em]
\shortstack[l]{Total ride\\time (min)}
    & 13,413.89 & 9,450.33 &  11,032.21 & -29.55 & -17.76 & 13,057.28 & 8,859.08 & 10,459.34 & -32.14 & -19.90\\
\scriptsize \textit{Avg per student (min)}
    & \scriptsize 23.87 & \scriptsize 16.82 & \scriptsize 17.07 & \scriptsize -29.55 & \scriptsize -28.49 & \scriptsize 30.09 & \scriptsize 20.41 & \scriptsize 21.13  & \scriptsize -32.14 & \scriptsize -29.78 \\
\addlinespace[0.5em]
\shortstack[l]{Total idle\\time (min)}
    &  4,175.89 &  3,568.66 & 4,098.88 & -14.54 & -1.84& 4,082.37 & 3,419.22  & 3,815.45 & -16.24 & -6.54\\
\scriptsize \textit{Avg per student (min)}
    & \scriptsize 7.43 & \scriptsize 6.35 & \scriptsize 6.35 & \scriptsize -14.54 & \scriptsize -14.67 & \scriptsize 9.41 & \scriptsize 7.88 & \scriptsize 7.71  & \scriptsize -16.24 & \scriptsize -18.07 \\
\addlinespace[0.5em]
\shortstack[l]{Total walking\\time (min)}
    & 1,264.91 & 1,711.77 &  2007.29 & 35.33 & 58.69 & 980.64 & 1,396.32 & 1,498.63 & 42.39 & 52.82 \\
\scriptsize \textit{Avg per student (min)}
    & \scriptsize 2.25 & \scriptsize 3.05 & \scriptsize 3.11 & \scriptsize 35.56 & \scriptsize 37.78 & \scriptsize 2.26 & \scriptsize 3.22 & \scriptsize 3.03  & \scriptsize 42.48 & \scriptsize 34.07 \\
\addlinespace[0.5em]
\shortstack[l]{No. of bus\\ commuters}
    & 562 & 562 & 646 & 0 & 14.95 & 434 & 434 & 495 & 0 & 14.06 \\
\addlinespace[0.5em]
\shortstack[l]{No. of car\\ commuters}
    & 639 & 639 &  555 & 0 & -13.15 & 597 & 597 & 536 & 0 & -10.22 \\
\addlinespace[0.5em]
\shortstack[l]{Utilization ratio}
    & 0.81 & 0.81 &  0.93 & 0 & 14.95 & 0.80 & 0.80 & 0.92 & 0 & 14.06 \\
\addlinespace[0.5em]
\shortstack[l]{No. of stops}
    & 179 & 143 & 146 & -20.11 & -18.44 & 208 & 154 & 160 & -25.96 & -23.08 \\
\addlinespace[0.5em]
Students per stop
    & 3.14 & 3.93 &  4.42 & 25.17 & 40.93 & 2.09 & 2.82 & 3.09 & 35.06 & 48.27 \\
\addlinespace[0.5em]
\shortstack[l]{Total {CCTS} (min)}
    & 6,668.8 & 6,668.8 & 6,071.7 & 0.00 & -8.95
    & 8,365.9 & 8,365.9 & 7,059.46 & 0.00 & -15.62 \\

\scriptsize \textit{Avg per student (min)}
    & \scriptsize 10.44 & \scriptsize 10.44 & \scriptsize 10.94 & \scriptsize 0.00 & \scriptsize 4.79
    & \scriptsize 14.01 & \scriptsize 14.01 & \scriptsize 13.17 & \scriptsize 0.00 & \scriptsize -6.01 \\
% \shortstack[l]{Total BTT (min)}
%     & 0 & 0 &  0 & 0 & 0 & 0 & 0 & 0 & 0 & 0 \\
% \scriptsize \textit{Avg per bus (min)}
%     & \scriptsize 0 & \scriptsize 0 & \scriptsize 0 & \scriptsize 0 & \scriptsize 0 & \scriptsize 0 & \scriptsize 0 & \scriptsize 0  & \scriptsize 0 & \scriptsize 0 \\
\bottomrule
\end{tabular}
\end{table}

\begin{figure}[h]\label{gs_dist1}
    \centering
    \subfigure{
        \includegraphics[width=0.47\textwidth]{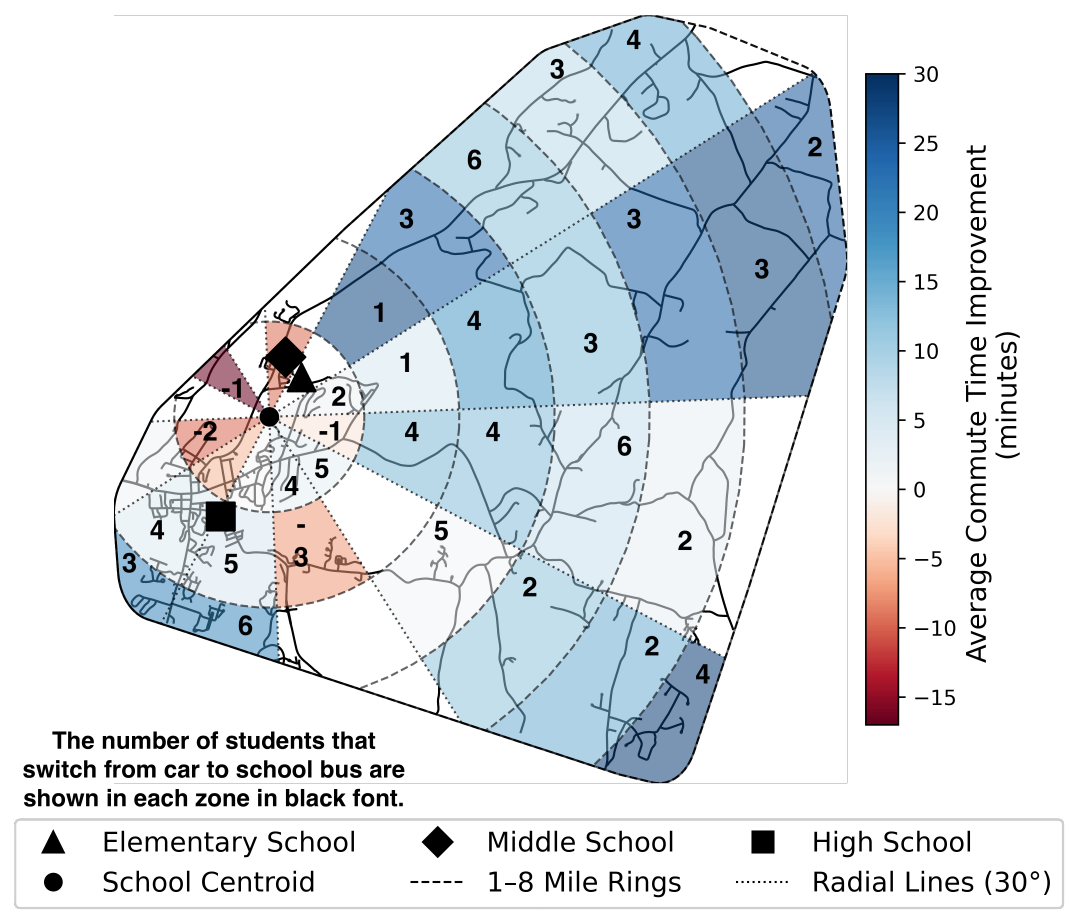}
        \label{Img:GeospaialSpread_Hanover}
    }
    \hfill
    \subfigure{
        \includegraphics[width=0.47\textwidth]{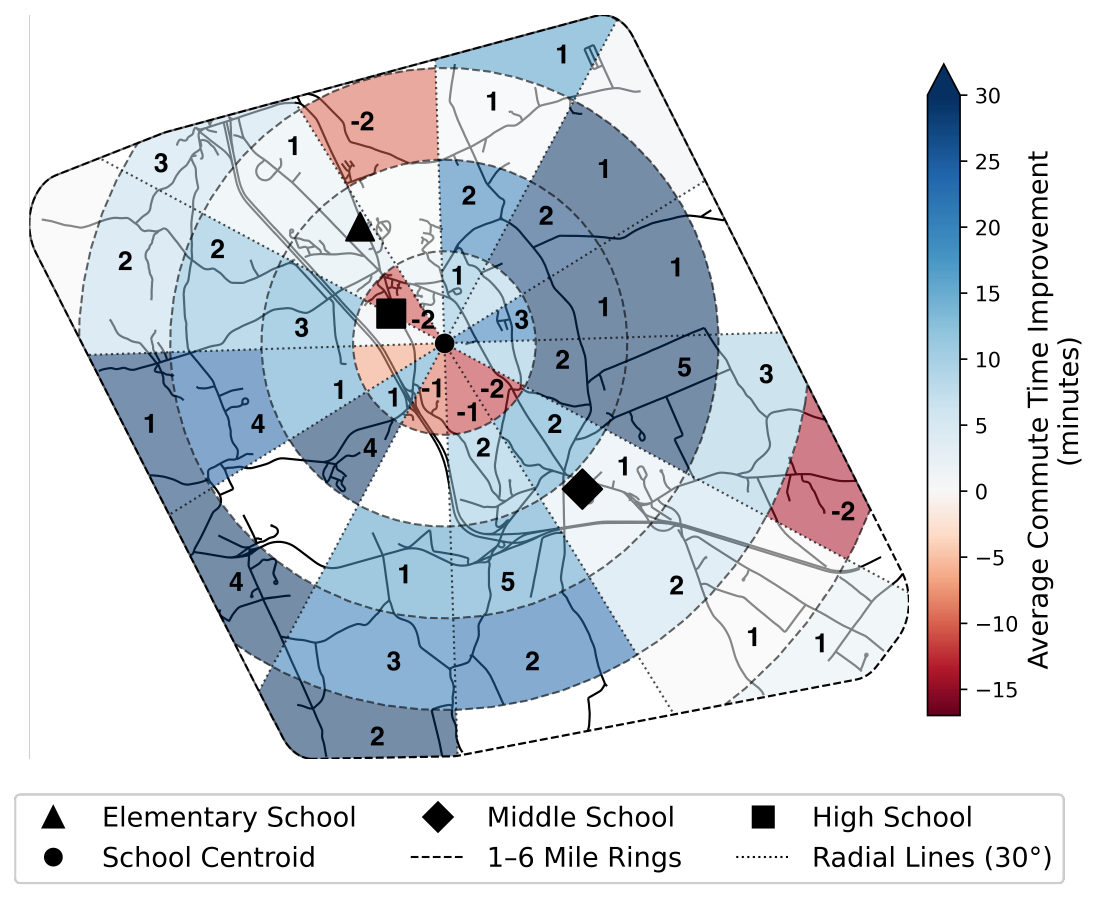}
        \label{Img:GeospatialSpread_Hopkinton}
    }
    \caption{Decreasing total commute times versus increasing bus ridership}
    \label{Img:GeospatialSpread}
\end{figure}
Table \ref{Table:Metrics_congestion} displays the results of the overall iterative optimization algorithm, including {total commute time, total bus ride time, total idle time, total walking time}, the number of bus and car commuters, utilization ratio, number of stops, number of students per stop, and total {car commute time of students (CCTS)}. \textit{Status Quo} corresponds to the existing school bus routes and schedules. \textit{Routing} shows the results of a single run of the cluster-then-route algorithm, focusing solely on route and schedule optimization without considering the mode shift. \textit{Routing + Mode} presents the results of our iterative framework. The \textit{$\Delta R$ (\%)} and \textit{$\Delta RM$ (\%)} columns quantify the percentage improvement achieved by the \textit{Routing} and \textit{Routing + Mode} solutions, respectively, compared to the \textit{Status Quo}.

{
Our solutions have substantially more students choosing the bus, raising the number of bus commuters from 562 to 646 in District I and from 434 to 495 in District II. This expanded bus ridership alleviates the underutilization of rural school buses and reduces operational expenses per student. In addition, between 61 and 84 cars are removed from the roads per day during peak hours, reducing localized congestion. Importantly, these gains are achieved with only marginal changes in average total commute time per student. In District I, the average bus ride time per student slightly increases from 16.82 to 17.07 min, and idle time remains nearly unchanged (at 6.35  min), while walking time also remains nearly unchanged (changing from 3.05 to 3.11 min). Similarly, in District II, average ride time, idle time, and walking time remain comparable. Overall, the incorporation of mode choice is found to substantially improve bus utilization and reduce car traffic while maintaining essentially the same level of commute quality.

Figure \ref{Img:GeospatialSpread} pinpoints where mode choice changes ridership and how those changes align with commute-time impacts. In both districts, the dominant pattern is synergy: zones shaded in blue (commute-time reductions) tend to be the same zones where we see increased ridership—students are more willing to adopt the bus when their total commute time improves. Conversely, the few zones associated with reduced ridership align with red shading, indicating that the additional access and operational burden in those locations translates into a small commute-time penalty that discourages bus use. The spatial structure of this ``win–lose” pattern differs by district. In District I, ridership declines are concentrated near the schools. In District II, the losing regions include a couple of distant pockets at the edge of the service area where commute times increase and ridership falls, which is plausible given District II’s overall longer baseline commutes (higher average commute times in Table 5), a more heterogeneous geography{, and slightly longer inter-school distances.} %For far-out students, serving them by bus may require detours and added walking/idle components that can dominate any routing efficiencies, making the bus less attractive in those specific corridors even as the district sees overall net ridership gains.

}

%% file: sources/8-Conclusion.tex
Rural school districts face unique challenges in designing bus routes and schedules, including long commutes, underutilized buses, localized congestion, irregular road networks, and mixed loading. This paper introduces a student-centric approach for designing rural school bus routes and schedules within multi-school districts with a fixed fleet.
%Mixed-load multi-school bus routing formulations are computationally intractable with commercial solvers and typically employ objectives that ignore student-level decisions.
We formulate an integrated mixed-integer linear optimization (MILO) model that incorporates student-level variables to minimize the {total commute time}. We develop an original two-phase cluster-then-route heuristic that leverages the road network structure and our MILO formulation to generate high-quality solutions for real-world instances. Furthermore, we combine our routing and scheduling approach with a mode choice model within an iterative framework to capture students' responses to the new bus routes and schedules. Improved routes and schedules attract more students to bus transportation, amplifying the systemwide benefits through a positive feedback loop. The results of two real-world rural school districts suggest that our approach reduces the {total commute time} by {$\sim 22$--$25\%$}. Additionally, accounting for student mode choice unlocks the potential to transport {$\sim 14$--$15\%$} more students.

In future research, our solution approach could be extended through time-dependent clustering to explore staggered school start times in rural settings as an additional lever to address congestion in the vicinity of schools. {For districts with greater inter-school travel distances, the clustering procedure could also account explicitly for students' destination schools, thereby improving the generalizability of our approach}. In addition, future work can collect and use disaggregate mode choice data to train more accurate mode choice models. Finally, even though our solution approach considerably outperforms the status quo, an optimality gap may still exist. This motivates further research into other heuristics or exact techniques to handle the complexity of multi-school mixed-loading rural school bus routing and scheduling.

%% file: sources/9-SI.tex
\section{Extended Solution Approach - Road Network-Aware Clustering}\label{appendix:RoadNetworkAwareClustering}

% {This appendix describes an extended approach to improve upon the Euclidean constrained $k$-means clustering solution from Section \ref{subsection:constrained-$k$-means} using road network information. Given the sparseness of road networks in rural school districts, students within a bus service region may appear geographically close, but the structure of the road network can result in substantial road distances between them. To address this, we identify a group of students whose reassignment to another cluster is likely to provide the largest benefit in terms of reducing the total commute time. Then, we solve a reduced form of the road network-aware constrained $k$-means clustering problem by focusing on reassigning this identified group of students to other clusters while fixing the cluster membership of the remaining students and respecting the upper limits on cluster sizes.}

{This appendix presents an extended approach to improve upon the Euclidean constrained $k$-means clustering solution introduced in Section \ref{subsection:constrained-$k$-means} by incorporating road network information. In rural school districts, the sparsity and structure of the road network can lead to situations where students who are geographically close are separated by substantial travel distances along the network.

Ideally, clustering would be performed directly using road network distances. However, as we demonstrate later in Section \ref{enhanced_approach_plus_ablation}, solving a fully road-network-aware constrained $k$-means clustering problem directly is computationally intractable for the instance sizes considered in our case studies. To address this challenge, we develop a extended approach that builds on the initial Euclidean $k$-means clustering solution and refines it using road network information.

Specifically, in Section \ref{Subsec:RoadNetworkAwareConstrainedKMeans}, we first present a road-network-aware constrained $k$-means clustering formulation. Then, in Section \ref{subsubsec:ClusteringReduction}, we present a reduced road-network-aware constrained $k$-means formulation, obtained by fixing the assignments of a subset of students to make it computationally tractable. Finally, Section \ref{enhanced_approach_plus_ablation} presents the overall approach and an extended ablation study to evaluate the impact of these additional components on solution quality.}

\subsection{Road Network-Aware Constrained $k$-means Clustering}\label{Subsec:RoadNetworkAwareConstrainedKMeans}
Consider a set $\mathcal{S}$ of students and a fleet $\mathcal{K}$ of buses. Let $C_k$ be the student carrying capacity of bus $k \in \mathcal{K}$. Let $\delta_{ss'}$ be the road network-based travel time between students $s$ and $s'$, where $s, s' \in \mathcal{S}$. Let $x_{ss'}$ be a binary variable that is 1 if and only if students $s \in \mathcal{S}$ and $s' \in \mathcal{S}$ are in the same cluster. Let $y_{sk}$ be a binary variable that is 1 if and only if student $s \in \mathcal{S}$ is assigned to bus $k \in \mathcal{K}$. The full road network-aware constrained $k$-means clustering formulation is as follows.
% \subsection{Reduced Road Network-Aware Constrained $k$-means Clustering}
% Consider a set $\mathcal{K}$ of buses and initial assignments of a set $\mathcal{S}$ of students provided by the Euclidean constrained $k$-means clustering step (Section \ref{subsection:constrained-$k$-means}). We define $\mathcal{S}_k$ as the set of students initially assigned to each bus $k\in\mathcal{K}$, creating a partition of all students: by definition, $\bigcup_{k\in \mathcal{K}} \mathcal{S}_k = \mathcal{S}$ and $\mathcal{S}_k \cap \mathcal{S}_{k'} = \emptyset \: \forall k, k' \in\mathcal{K}, k\neq k'$. Note that while the sets of students initially assigned to different buses do not overlap, their convex hulls (representing the bus service regions) may, and often do, have an overlap.

% Let $\overline{C}_k$ be the student carrying capacity of the bus $k \in \mathcal{K}$. Let $\delta_{ss'}$ be the road network-based travel time between students $s, s' \in \mathcal{S}$. Let $x_{ss'}$ be a binary variable which is 1 if and only if students $s, s' \in \mathcal{S}$ are in the same cluster. Let $y_{sk}$ be a binary variable which is 1 if and only if student $s \in \mathcal{S}$ is assigned to bus $k \in \mathcal{K}$. The full road network-aware constrained $k$-means formulation is as follows.
\begin{subequations}
\renewcommand{\theequation}{\theparentequation.\arabic{equation}}
\begin{align}
    \text{Minimize} &\sum_{s \in \mathcal{S}} \sum_{s' \in \mathcal{S}} \delta_{ss'}\cdot x_{ss'} \label{Eq:HybridkMeans_Objective}\\
\text{s.t.} \quad & x_{ss'} \geq y_{sk} + y_{s'k} - 1 \quad \forall s \in \mathcal{S}, s' \in \mathcal{S}, k \in \mathcal{K} \label{Eq:HybridkMeans_ClusterImprovement}\\
& \sum_{s \in \mathcal{S}} y_{sk} \leq C_k \quad \forall k \in \mathcal{K} \label{Eq:HybridkMeans_CapacityConstraint}\\
& \sum_{k \in \mathcal{K}} y_{sk} = 1 \quad \forall s \in \mathcal{S} \label{Constraint:UniqueClusterForPoint}\\
&x_{ss'} \in \{0,1\} \quad \forall s, s' \in \mathcal{S}, s \neq s'\label{DoD:RNAwareKMeans1}\\
&y_{sk} \in \{0,1\} \quad \forall s \in \mathcal{S}, k \in \mathcal{K}\label{DoD:RNAwareKMeans2}
\end{align}
\end{subequations}
The objective function \eqref{Eq:HybridkMeans_Objective} minimizes the sum of the distances between all pairs of students assigned to the same bus. Constraints \eqref{Eq:HybridkMeans_ClusterImprovement} ensure that two students assigned to the same bus are in the same cluster. Constraints \eqref{Eq:HybridkMeans_CapacityConstraint} enforce the bus capacity limits. Constraints \eqref{Constraint:UniqueClusterForPoint} ensure that each student is assigned to exactly one bus. Constraints \eqref{DoD:RNAwareKMeans1}-\eqref{DoD:RNAwareKMeans2} define the domains of the decision variables.

% The reduced form of this full formulation is obtained by fixing the assignment of students who satisfy two conditions: (i) they are not in the overlap of multiple bus service regions, and (ii) they are in the largest connected component of the road network within their cluster. We refer to students satisfying condition (i) as \textit{intra-service region students} and students satisfying condition (ii) as those connected to the \textit{main road network component} of their cluster. We now formally define the set $\mathcal{S}_{\mathcal{G}_k}$ of students satisfying both conditions for each bus $k$.
{
\subsection{Reduced Road Network-Aware Constrained $k$-means { and Bus Stop Clustering}}\label{subsubsec:ClusteringReduction}

Solving the full road network-aware constrained $k$-means formulation is computationally challenging for the problem sizes considered in our case studies. To address this, we fix the assignments of a subset of students whose cluster membership is unlikely to change under road network distances and allow reassignment only for the remaining students.

Let $\mathcal{S}_k$ denote the set of students initially assigned to bus $k \in \mathcal{K}$. For each bus $k$, let $\mathcal{S}_{\mathcal{G}_k} \subseteq \mathcal{S}_k$ denote the subset of students whose assignments are fixed in the reduced problem.

Accordingly, the reduced formulation is obtained by imposing the following constraints:
\begin{subequations}
\begin{align}
&x_{ss'} = 1 \quad \forall s, s' \in \mathcal{S}_{\mathcal{G}_k}, \, k \in \mathcal{K} \label{students_in_same_cluster_certain}\\
&y_{sk} = 1 \quad \forall s \in \mathcal{S}_{\mathcal{G}_k}, \, k \in \mathcal{K} \label{student_in_a_given_cluster_certain}
\end{align}
\end{subequations}

We now describe how the set $\mathcal{S}_{\mathcal{G}_k}$ is constructed. Let $q_s = (x_s, y_s)$ denote the geographic location of student $s \in \mathcal{S}$. For each bus $k$, let $H_k$ denote the service region (convex hull of assigned students) defined by Equation \eqref{eq:CH} in Section~\ref{subsection:constrained-$k$-means}.

\textbf{Inter-service region students.}
We define $\mathcal{S}_{H_k}$ as the set of students assigned to bus $k$ who are located within the intersection of the service region $H_k$ with the service region of at least one other bus, and $\mathcal{S}_{\tilde{H}_k}$ as the set of students assigned to bus $k$ who are not located in the service region of any other bus.
\begin{align}
    \mathcal{S}_{H_k} &= \bigcup_{j \neq k} \{ s \in \mathcal{S}_k: q_s \in H_k \cap H_j \}, \\
    \mathcal{S}_{\tilde{H}_k} &= \mathcal{S}_k \setminus \mathcal{S}_{H_k}.
\end{align}

\textbf{Intra-service region students.}
Let $\mathscr{G}_{H_k} = (\mathscr{N}_{H_k}, \mathscr{A}_{H_k})$ denote the road network restricted to the service region $H_k$, and let $\mathscr{G}_{H_k}^{\text{max}}$ denote its largest connected component. We define $\mathcal{S}_{\mathcal{G}_k}$ as the set of students within $\mathcal{S}_{\tilde{H}_k}$ whose residential locations are at the nodes or on the arcs of subgraph $\mathscr{G}_{H_k}$.
\begin{align}
    \mathcal{S}_{\mathcal{G}_k} = \{ s \in \mathcal{S}_{\tilde{H}_k} : q_s \in \mathscr{G}_{H_k}^{\text{max}} \}.
\end{align}

Students in $\mathcal{S}_{H_k}$ (inter-service region students) and those in $\mathcal{S}_{\tilde{H}_k} \setminus \mathcal{S}_{\mathcal{G}_k}$ (intra-service region students not connected to the main component) are allowed to change cluster membership, while students in $\mathcal{S}_{\mathcal{G}_k}$ are fixed. 

The optimal solution to the reduced road network-aware constrained $k$-means clustering formulation, given by Equations \eqref{Eq:HybridkMeans_Objective}-\eqref{DoD:RNAwareKMeans2}, \eqref{students_in_same_cluster_certain}-\eqref{student_in_a_given_cluster_certain}, yields an updated assignment of students to buses, which induces a new set of service regions $\dot{H}_k$ for each bus $k \in \mathcal{K}$. Let $\dot{\mathcal{S}}_k$ denote the set of students assigned to bus $k$ in the optimal solution to the reduced road network-aware constrained $k$-means clustering formulation. Then $\dot{H}_k$ is defined as follows.
\begin{align}
\dot{H}_k = \left\{ q \in \mathbb{R}^2 : q = \sum_{s \in \dot{\mathcal{S}}_k} \alpha_s q_s, \ \sum_{s \in \dot{\mathcal{S}}_k} \alpha_s = 1, \ \alpha_s \geq 0 \:\:\: \forall s \in \dot{\mathcal{S}}_k \right\} \forall k \in \mathcal{K}.
\end{align}

To account for allowable walking distances, similar to what we did in Section \ref{subsubsec:BusStopClustering}, we expand each service region $\dot{H}_k$ by the maximum allowable walking distance $d^{\max}$, resulting in the updated service regions $\dot{H}_k'$, defined as follows.
\begin{align}
\dot{H}_k' = \left\{ q \in \mathbb{R}^2 : \exists \, q' \in \dot{H}_k \text{ such that } \|q - q'\|_2 \leq d^{\max} \right\}.
\end{align}

These expanded service regions define the feasible pickup areas for each bus, ensuring that every student is within a distance of $d^{\max}$ from at least one candidate stop within the same service region. In the subsequent routing stage, these regions are used to restrict the set of candidate stops and feasible assignments in each single-bus routing problem, which is solved using our SAS4 procedure to generate a high-quality route for each bus.
}

\subsection{Extended Solution Approach and Resulting Improvements}\label{enhanced_approach_plus_ablation}
\begin{figure}[h!]
    \centering
    \includegraphics[width=0.85\textwidth]{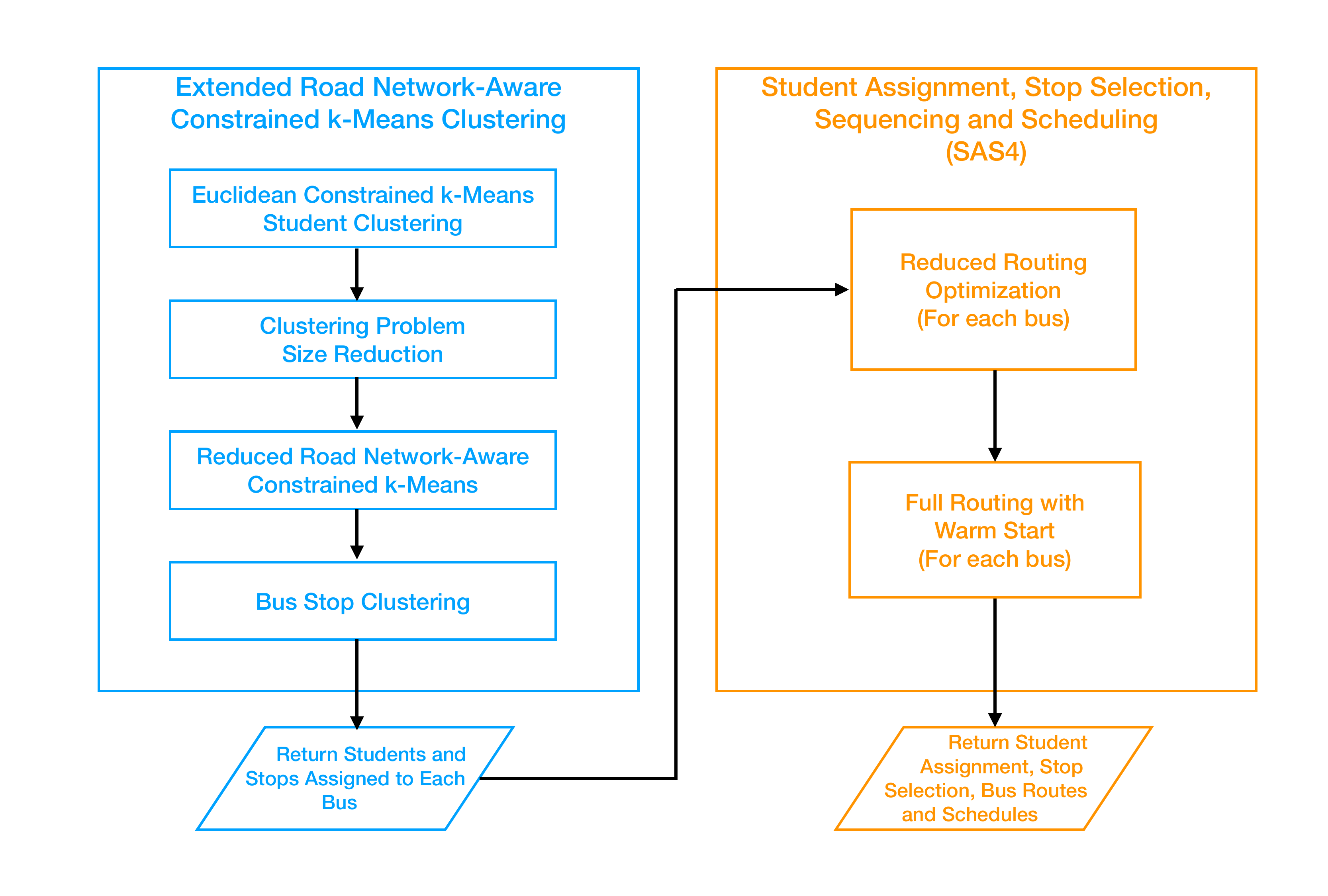}
    \caption{Extended solution approach: A road network-aware extension of our overall solution approach.}
    \label{fig:SBRPSolutionApproach_expanded}
\end{figure}

{Our extended solution approach, which incorporates road network-aware clustering, is shown in Figure \ref{fig:SBRPSolutionApproach_expanded}. We also present an additional ablation study reported in Table \ref{AblationTable_Full_SI} to show the value of each component of this extended approach. In this updated table, Ablation 2 corresponds to the ``Base" configuration reported in Section \ref{sec:BenefiitsAlgoDesign}. Finally, we provide visual illustrations showing the sequence of steps used to construct service regions under both the Euclidean Constrained $k$-means Clustering approach in Section \ref{subsection:constrained-$k$-means-main} and the Extended Road Network-Aware Constrained $k$-Means Clustering approach in this appendix.}

\begin{table}[h]
\centering
\caption{Extended ablation study results. ``Infeasible'' indicates configurations where the routing phase was infeasible for at least one cluster.``No sol.'' indicates configurations where no solution was found for any of the clusters. ``No sol. $^*$'' indicates configurations where no solution was found for at least one of the clusters.}
\label{AblationTable_Full_SI}
\begin{tabular}{@{}L{5.7cm}L{2.7cm}rrrr@{}}
\toprule
& \multirow{2}{*}{\begin{tabular}[c]{@{}r@{}}Feature(s)\\disabled\end{tabular}} & \multicolumn{2}{c}{District I} & \multicolumn{2}{c}{District II} \\
\cmidrule(lr){3-4} \cmidrule(lr){5-6}
Ablation & & Gap \% & Time (h) & Gap \% & Time (h) \\
\midrule
Extended base & - & 0 & 25.3 & 0 & 19.7 \\
1. No Euclidean constrained k-means \\ \:\: \& No clustering problem size reduction & \ref{subsection:constrained-$k$-means}, \ref{subsubsec:ClusteringReduction}  & Infeasible & 61.3  & Infeasible & 52.6  \\
2. No road network-awareness & \ref{subsubsec:ClusteringReduction}, \ref{Subsec:RoadNetworkAwareConstrainedKMeans} & 1.3 & 25.2 & 2.7 & 19.5\\
3. No clustering problem size reduction & \ref{subsubsec:ClusteringReduction} with warm-starting \ref{Subsec:RoadNetworkAwareConstrainedKMeans}  & No sol.$^*$ & 52.1 & No sol.$^*$ & 44.4 \\
4. No removal of arcs $\mathcal{A}_5$ & \textit{\underline{rdn.1 }} in \ref{Subsec:Reduced_Routing_Optimization} & 0 & 26.9 & 0 & 20.6 \\
5. No modification to the routing objective & \textit{\underline{rdn.2 }} in \ref{Subsec:Reduced_Routing_Optimization} & No sol.$^*$ & 28.4 & No sol.$^*$ & 22.1 \\
6. No student-to-stop pre-assignment & \textit{\underline{rdn.3 }} in \ref{Subsec:Reduced_Routing_Optimization} & No sol.$^*$ & 36.3 & No sol.$^*$ & 28.3 \\
7. No reduced routing but with pre-assignment & \ref{Subsec:Reduced_Routing_Optimization} with Eqs \eqref{stopminimizinfstudentassignment:OF}-\eqref{stopminimizinfstudentassignment:end} & 29.23 & 25.3 & 21.74 & 19.6 \\
8. No reduced routing &  \ref{Subsec:Reduced_Routing_Optimization} & No sol. & 25.2 & No sol. & 19.6 \\
\bottomrule
\end{tabular}
\end{table}

\clearpage
\begin{figure}[p]
\centering

\begin{minipage}{0.48\textwidth}
    \centering
    \includegraphics[width=\linewidth,height=0.25\textheight,keepaspectratio]{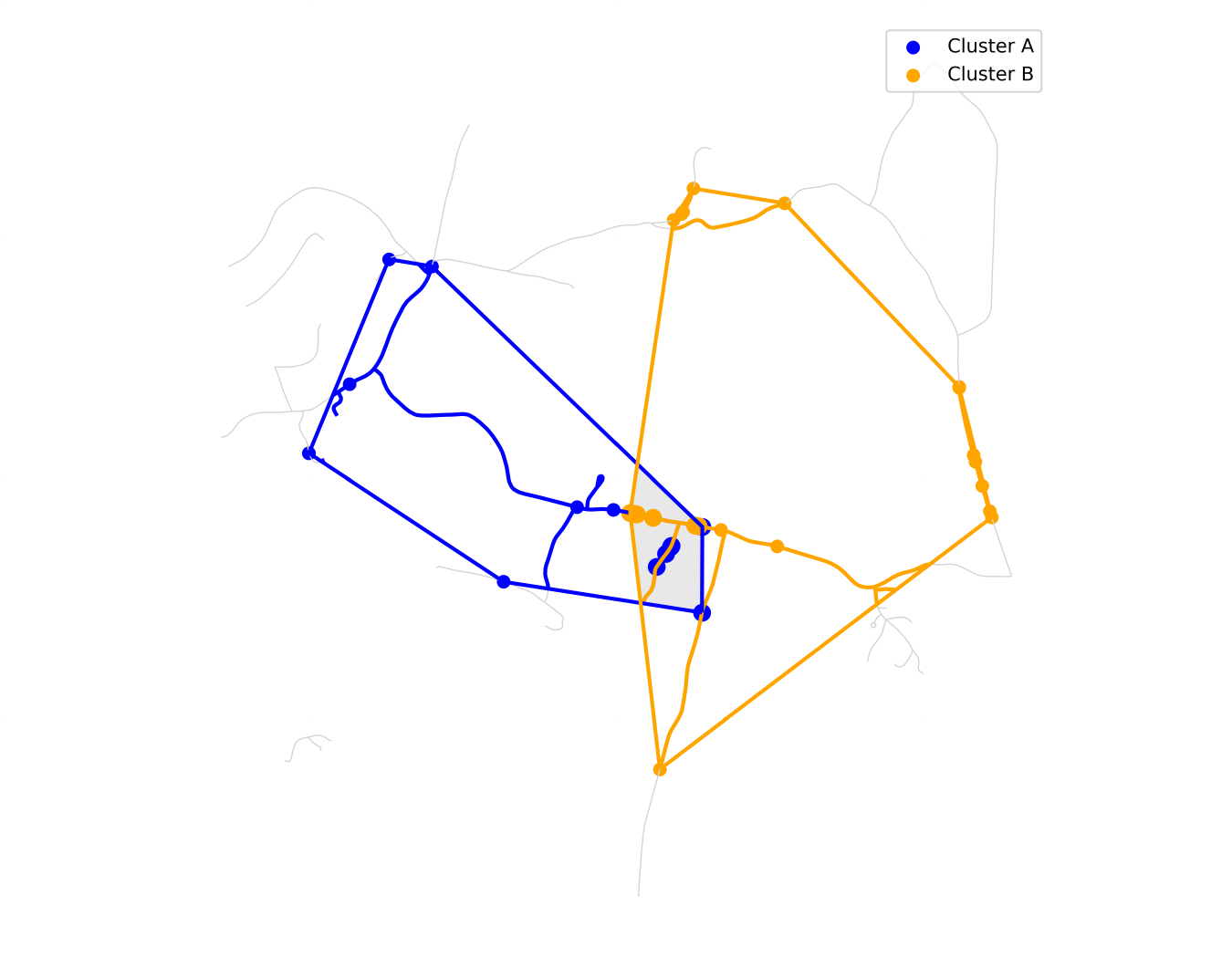}
    \\[0.5ex]
    (a) $H_k$: Euclidean constrained-$k$ means\\clusters of student locations
\end{minipage}\hfill
\begin{minipage}{0.48\textwidth}
    \centering
    \includegraphics[width=\linewidth,height=0.25\textheight,keepaspectratio]{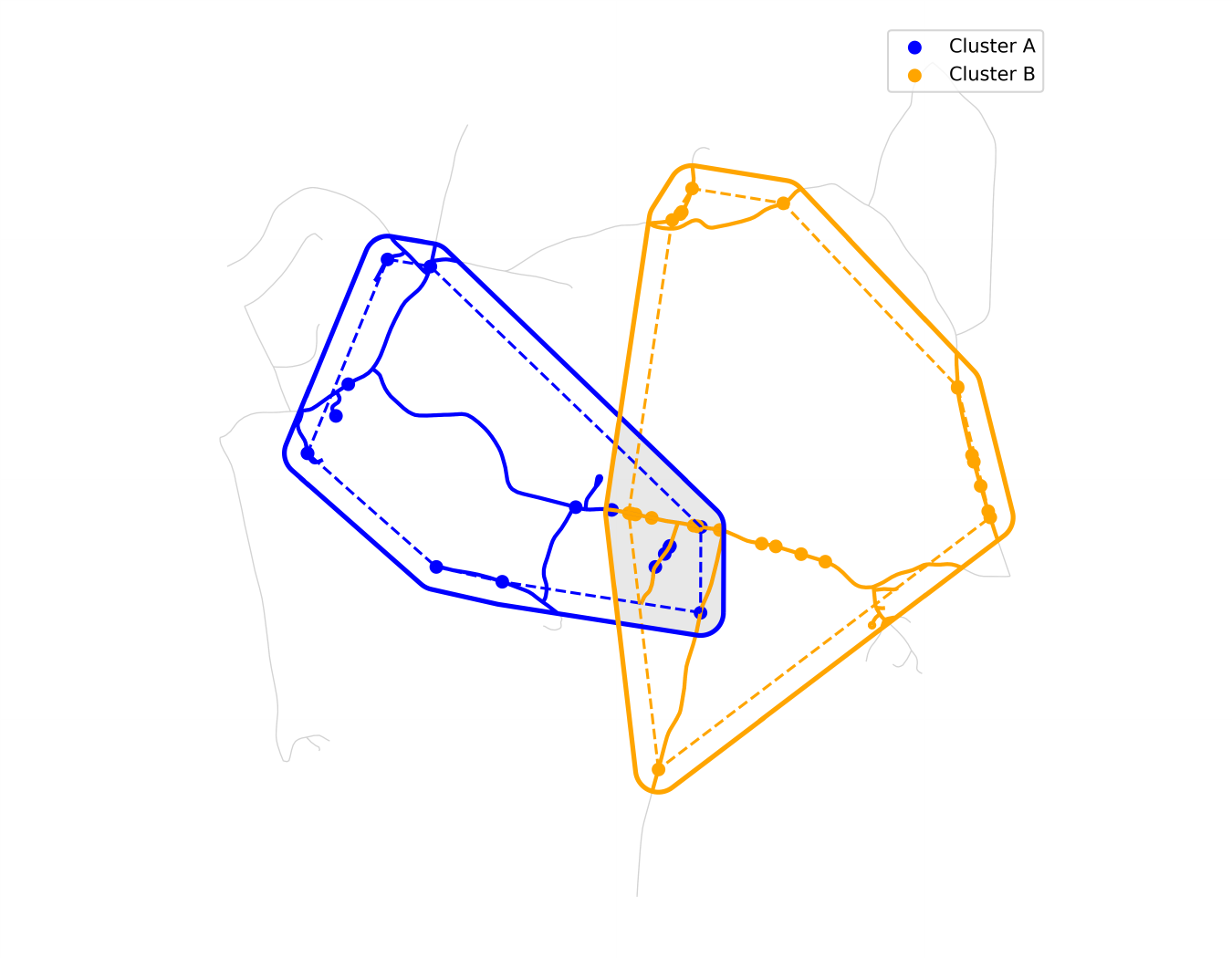}
    \\[0.5ex]
    (b) $H_k^{\prime}$: Service regions of buses\\in our solution approach
\end{minipage}

\vspace{0.02\textheight}

\begin{minipage}{0.48\textwidth}
    \centering
    \includegraphics[width=\linewidth,height=0.25\textheight,keepaspectratio]{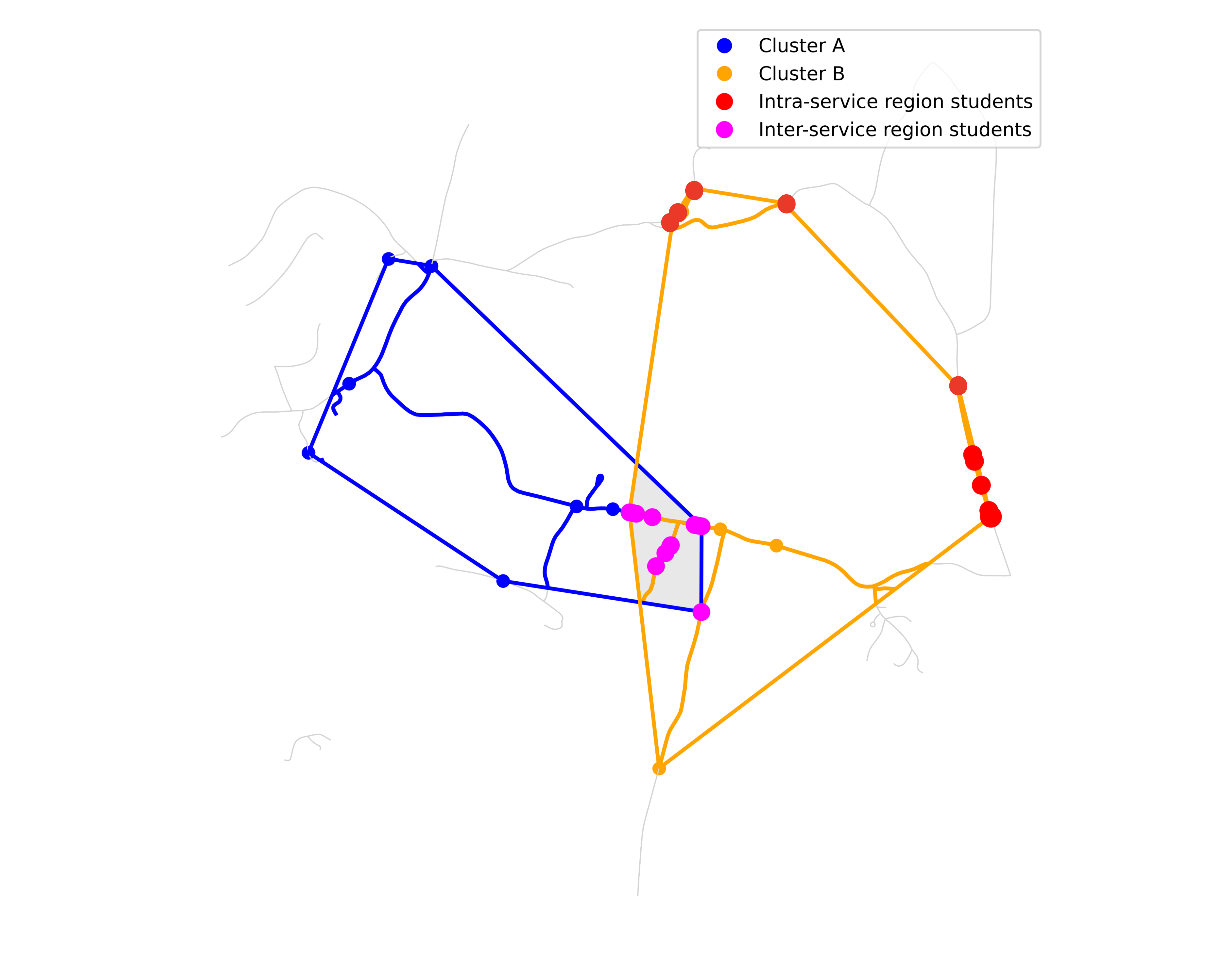}
    \\[0.5ex]
    (c) Identification of $\mathcal{S}_{H_k}$ (inter-service region students) and $\mathcal{S}_{\tilde{H}_k} \setminus \mathcal{S}_{\mathcal{G}_k}$ (intra-service region students not connected to the main component)
\end{minipage}\hfill
\begin{minipage}{0.48\textwidth}
    \centering
    \includegraphics[width=\linewidth,height=0.25\textheight,keepaspectratio]{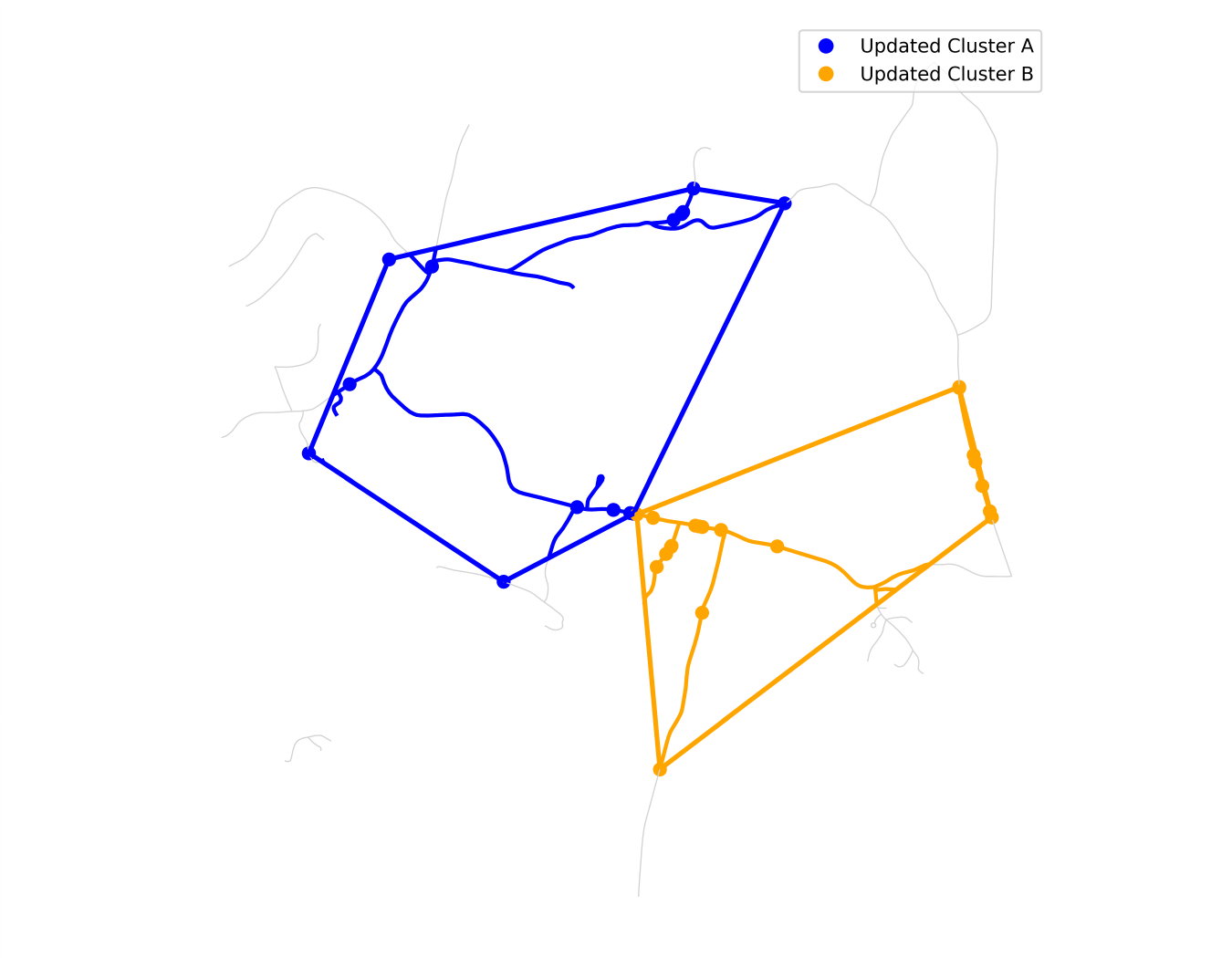}
    \\[0.5ex]
    (d) $\dot{H}_k$: Updated clusters after extended road network-aware constrained-$k$ means clustering 
\end{minipage}

\vspace{0.02\textheight}

\begin{minipage}{0.48\textwidth}
    \centering
    \includegraphics[width=\linewidth,height=0.25\textheight,keepaspectratio]{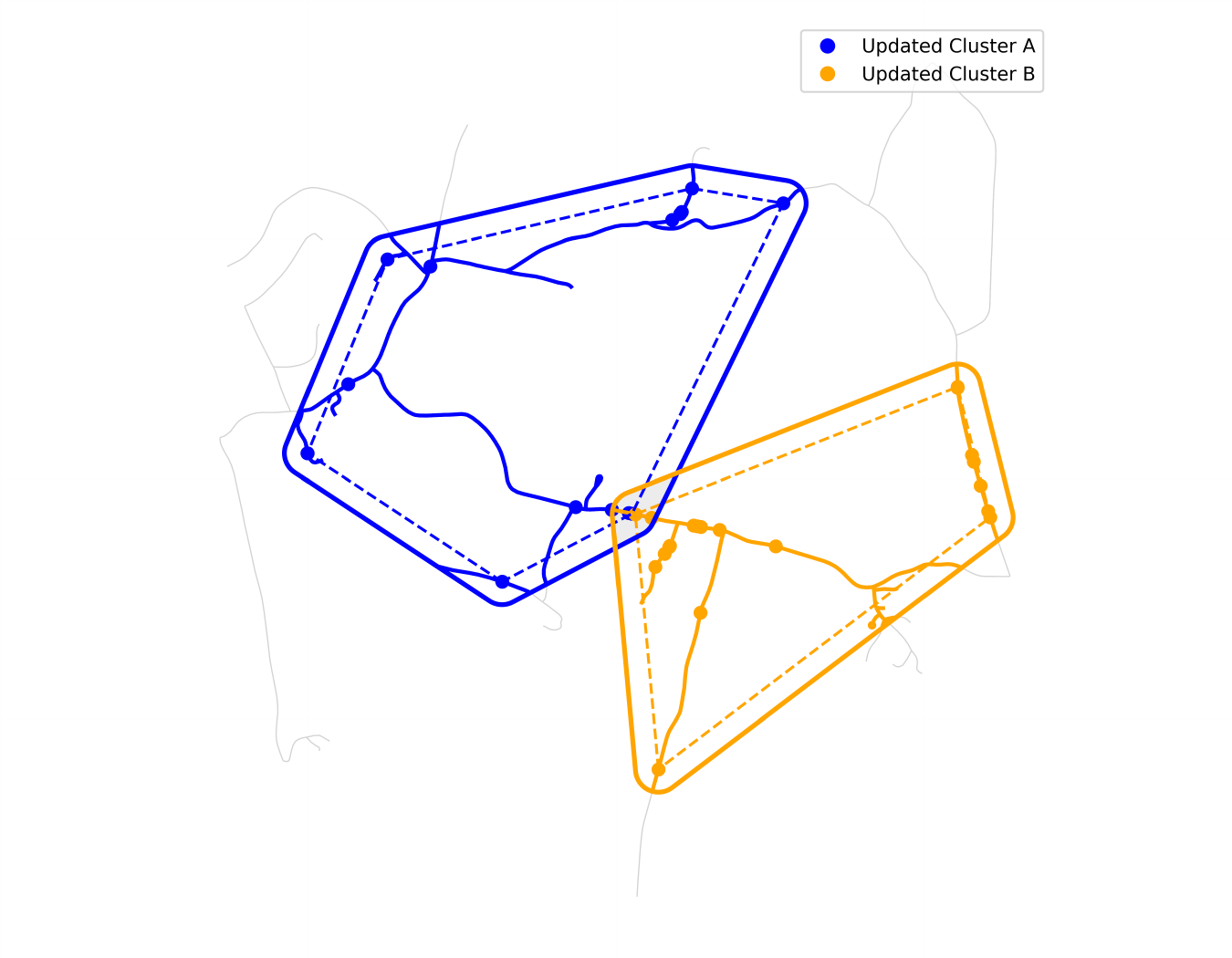}
    \\[0.5ex]
    (e) $\dot{H}_k^{\prime}$: Service regions of buses in our extended solution approach
\end{minipage}\hfill

\caption{Toy example illustrating the construction of service regions under the Euclidean Constrained \(k\)-Means Clustering approach (Subfigures (a) and (b)) and Extended Road Network-Aware Constrained $k$-Means Clustering approach (Subfigures (a), (c), (d), and (e)). %Subfigures (a) and (b) illustrate the steps in the Euclidean Constrained \(k\)-Means Clustering approach. Subfigures (a), (c), (d), and (e) illustrate the sequence of steps in our Extended Road Network-Aware constrained $k$-Means Clustering approach.
}
\label{fig:3x2_layout}
\end{figure}
\clearpage
\section{Case Study Details}\label {AppA}
\begin{figure}[ht]
    \centering
    \begin{minipage}{0.50\textwidth}
        \centering
        \includegraphics[width=\textwidth]{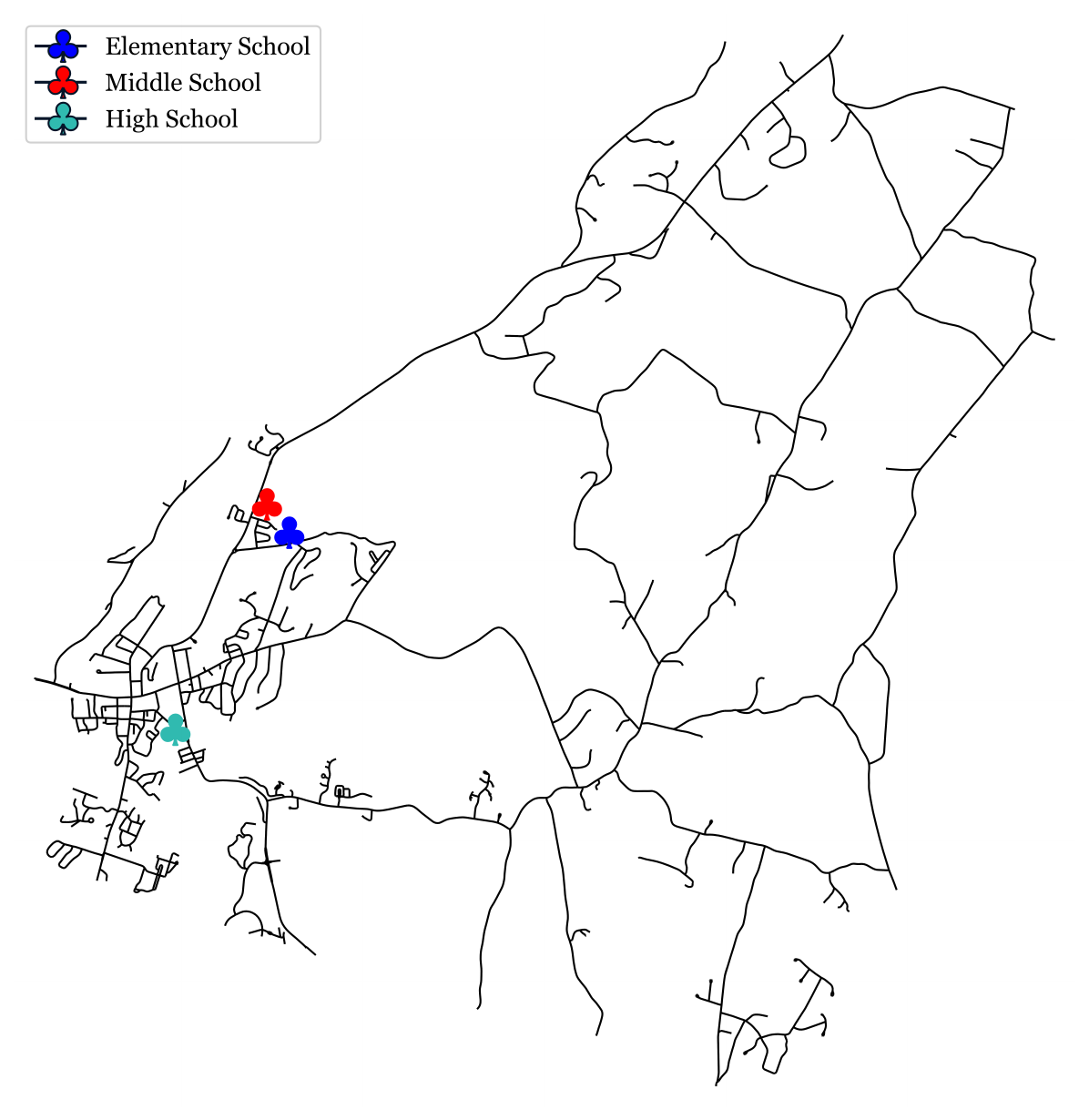}
    \end{minipage}\hfill
    \begin{minipage}{0.50\textwidth}
        \centering
        \includegraphics[width=\textwidth]{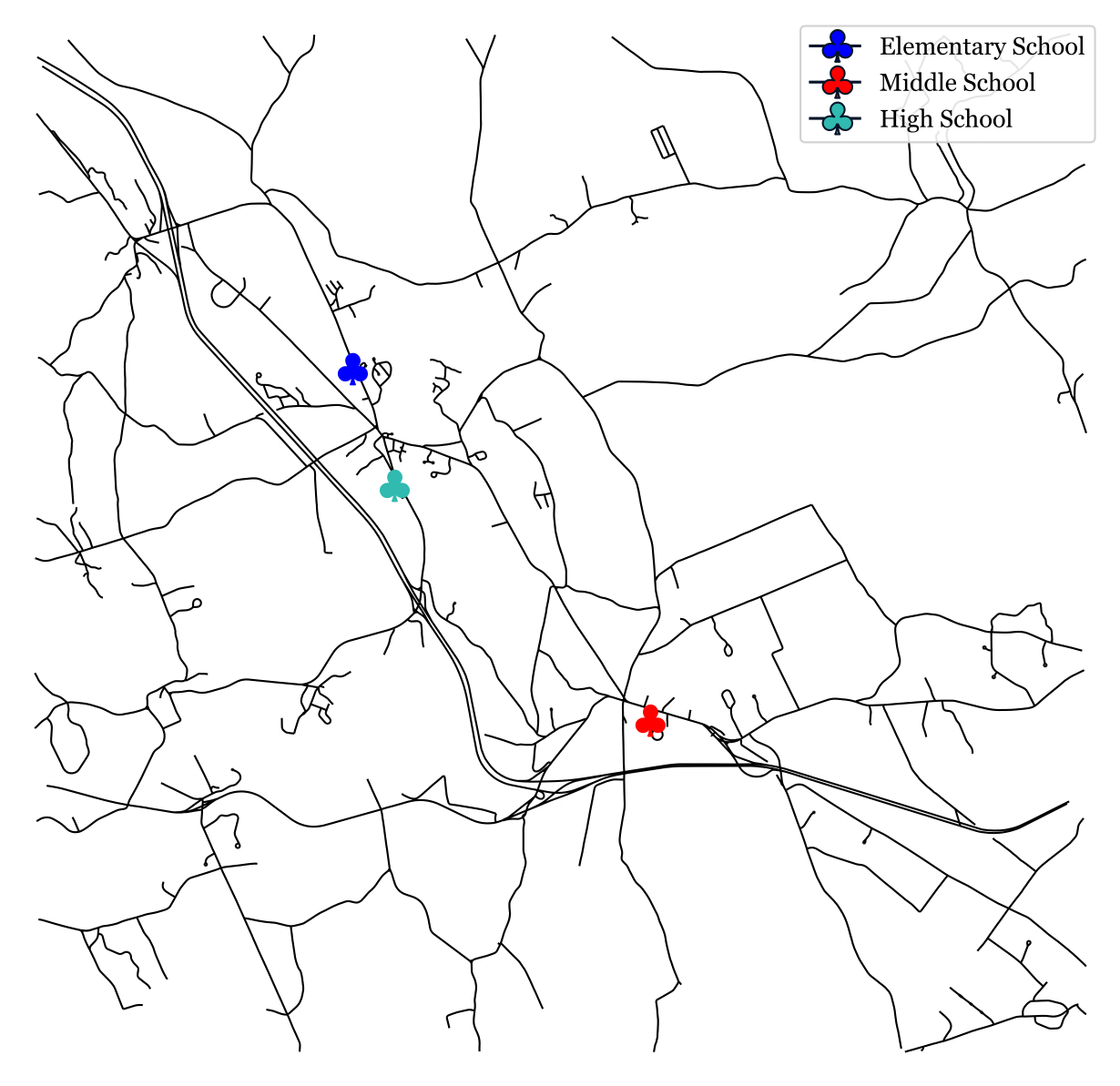}
    \end{minipage}
    \caption{Road network structures of our case study school districts}
    \label{fig:twocolumnfig}
\end{figure}

The road network topologies of the entire service areas where students are picked up in both school districts are shown in Figure \ref{fig:twocolumnfig}. Table \ref{Input_Para_table} lists the key input parameters of our case studies, as provided by our partner districts. 

We obtain information on existing bus routes and stop locations from publicly available data on school district websites. Information on the proportion of students who always, sometimes, or never use bus transportation is derived from the parent survey described in Section \ref{Subsec:SurveyandCounts}, assuming that these patterns apply to both school districts in our case studies. Information on the daily median ridership and the school-by-school split is obtained from the school districts. Information about the structure of the road network, travel speeds, and distances between student locations, candidate stop locations, and schools is obtained using OpenStreetMap \citep{Haklay:2008:00}, accessed through the OSMNX package in Python.  

\newpage
\begin{table}[h]
\centering
\caption{Input Parameters}
\label{Input_Para_table}
\begin{tabular}{lcc}
\toprule
\textbf{Parameter} & \textbf{District I} & \textbf{District II} \\
\midrule
Number of schools & 3 & 3 \\
Total number of students (excluding walkers/bikers) & 1,201 & 1,031 \\
%Percentage of students that always take the bus & 40\% & 40\% \\
%Percentage of students that sometimes take the bus & 20\% & 20\% \\
%Percentage of students that never take the bus & 40\% & 40\% \\
Median daily bus ridership & 562 & 434 \\
Elementary school bus commuters & 289 & 219 \\
Middle school bus commuters & 124 & 116 \\
High school bus commuters & 149 & 99 \\
Number of candidate bus stops & 179 & 208 \\
Number of buses & 9 & 7 \\
Bus capacity & 77 & 77 \\
Maximum time allowed for each bus route ({$T^B$}) (seconds) & 4,020 & 4,020 \\
{Maximum allowable student time on the bus ($T^S$) (seconds)} & {3,600} & {3,600} \\
Maximum allowed walking distance (miles) & 0.3 & 0.3 \\
Number of intersections & 409 & 382 \\
Total length of roads (miles) & 113.38 & 191.00 \\
Number of road segments & 1,206 & 1,117 \\
Average degree for nodes in road network & 2.04 & 2.03 \\
The value of tolerance $\epsilon$ for clustering (square-km) & $10^{-4}$ & $10^{-4}$ \\
\bottomrule
\end{tabular}
\end{table}
\newpage
\section{Structure of the Resulting Routes}\label{Section:Structure}

\begin{figure}[p]
\centering
\begin{minipage}{0.48\textwidth}
\centering
\includegraphics[width=0.9\linewidth,height=0.40\textheight,keepaspectratio]{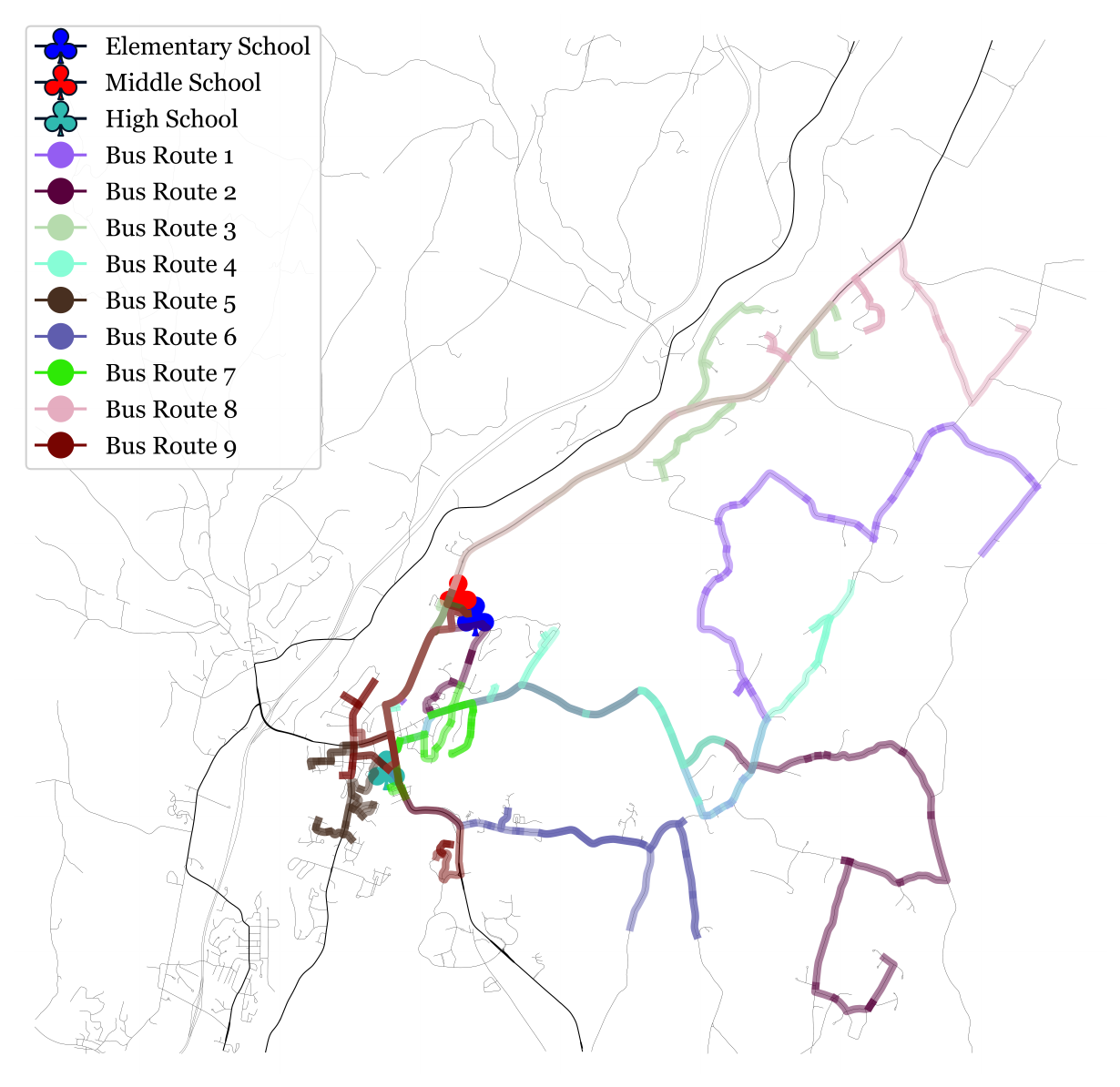}
\\
Status quo routes -- District I
\end{minipage}
\hfill
\begin{minipage}{0.48\textwidth}
\centering
\includegraphics[width=0.9\linewidth,height=0.40\textheight,keepaspectratio]{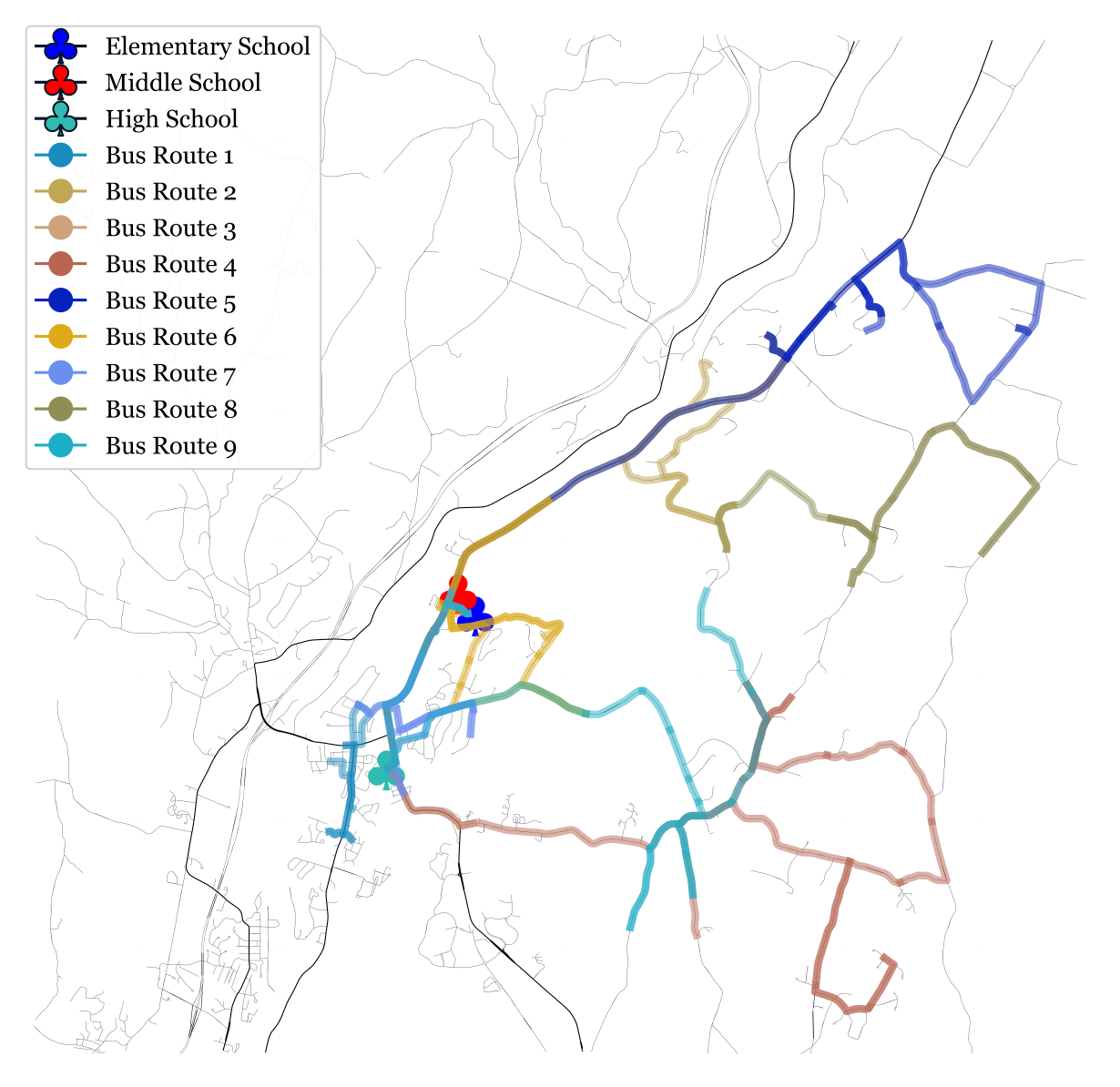}
\\
{Our solution -- District I}
\end{minipage}

\vspace{0.03\textheight}

\begin{minipage}{0.48\textwidth}
\centering
\includegraphics[width=0.9\linewidth,height=0.40\textheight,keepaspectratio]{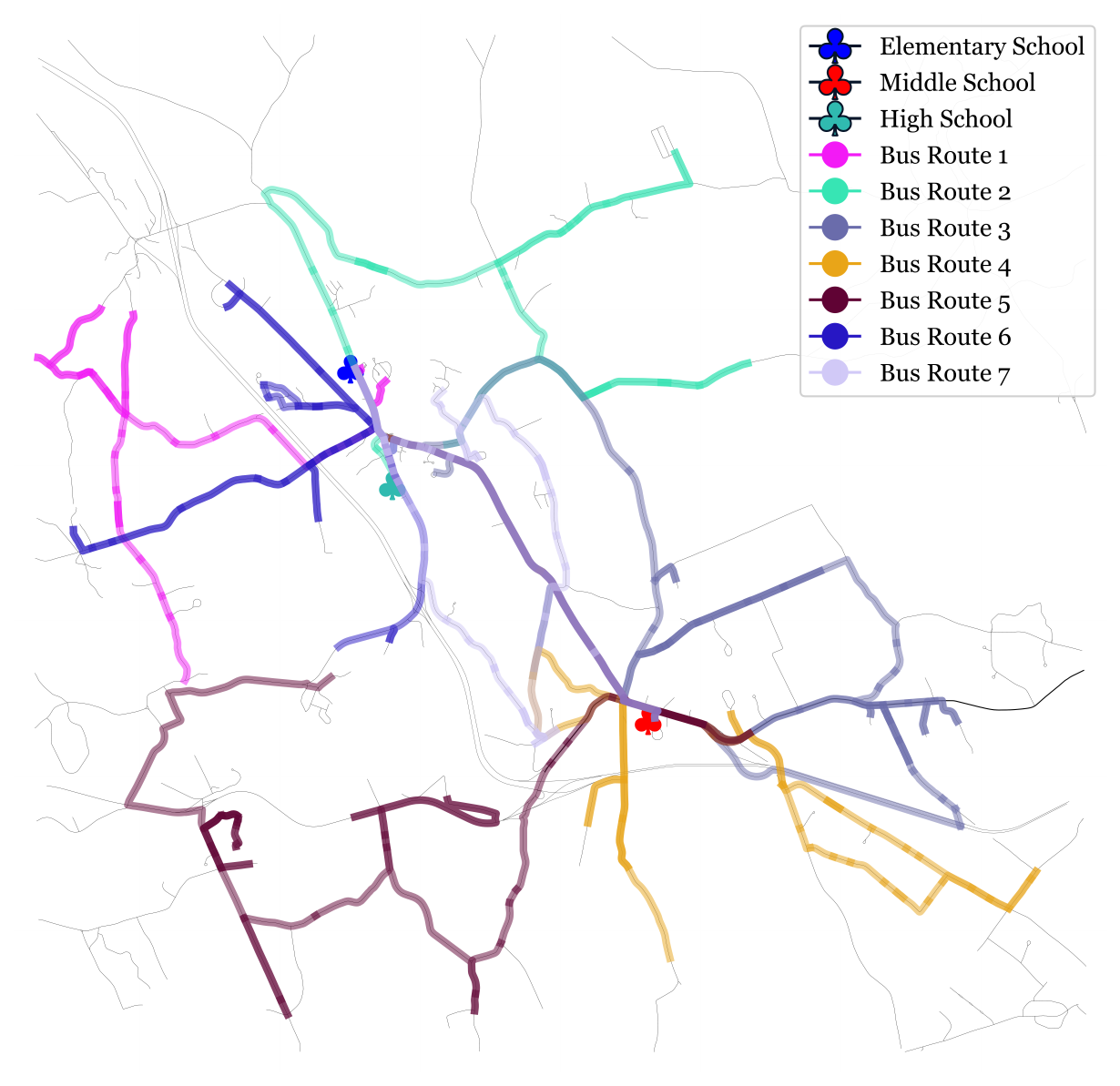}
\\
Status quo routes -- District II
\end{minipage}
\hfill
\begin{minipage}{0.48\textwidth}
\centering
\includegraphics[width=0.9\linewidth,height=0.40\textheight,keepaspectratio]{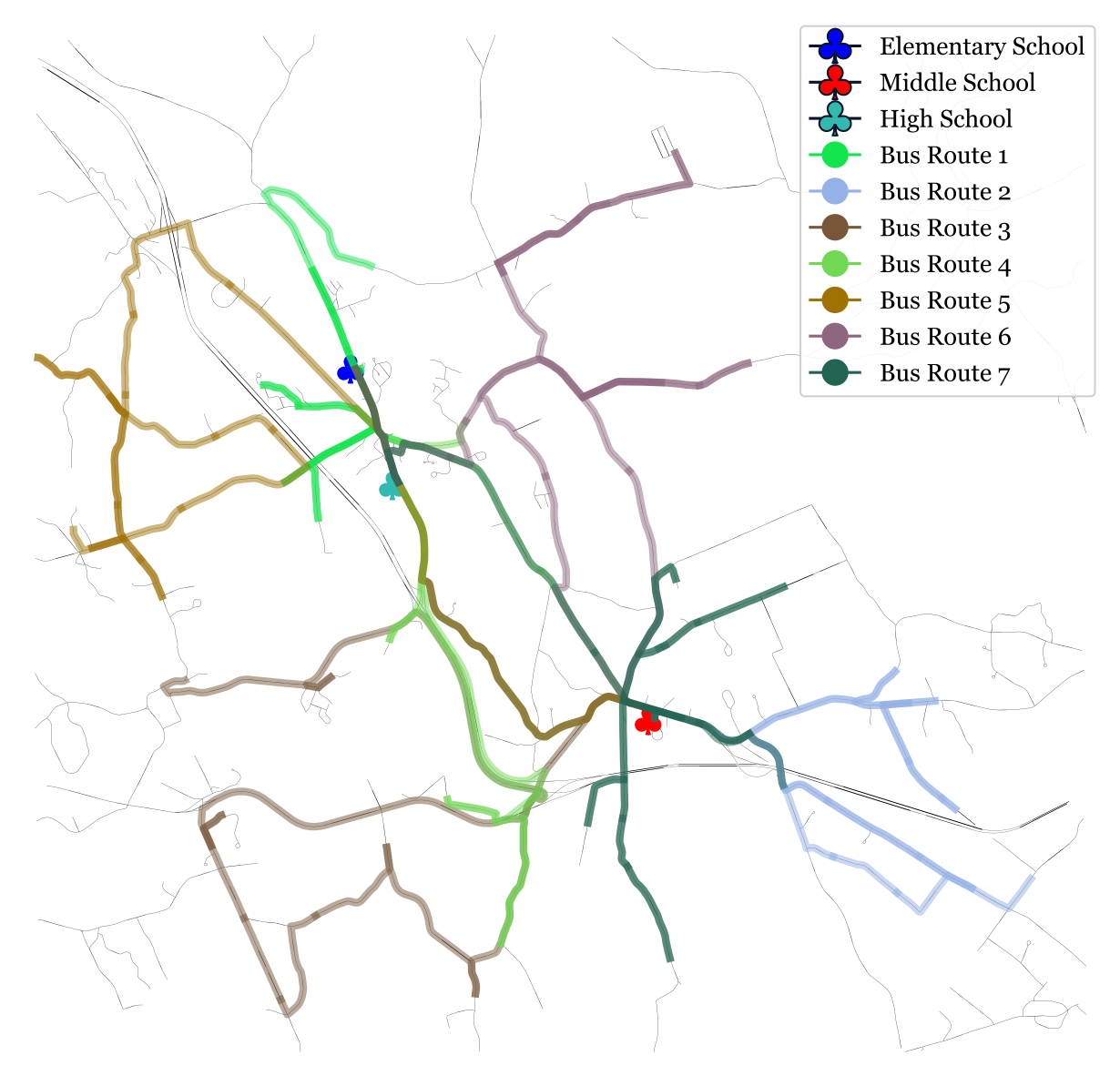}
\\
Our solution -- District II
\end{minipage}
\vspace{0.03\textheight}
\caption{Comparison of status quo routes and routes proposed in our solutions}
\label{fig:routecomparison}
\end{figure}

\begin{table}[H]
\centering
\scriptsize
\caption{Distributional Statistics of {Student Commute Times} by School and School District {(Min)}}
\label{Table:CommuteTimePercentilesMinutes}
\begin{tabular}{llrrrrrrrrr}
\hline
District & School & Commute Time & Mean & 25\% & 50\% & 75\% & 80\% & 85\% & 90\% & 95\% \\
\hline
District I & Elementary School & Status Quo & 32.7 & 23.9 & 33.3 & 40.3 & 41.6 & 43.4 & 46.7 & 50.4 \\
&  & Our Solution & 26.2 & 21.5 & 25.2 & 30.4 & 31.3 & 33.2 & 36.3 & 40.8 \\
&  & $\Delta$ & 6.5 & 2.4 & 8.1 & 9.9 & 10.3 & 10.2 & 10.4 & 9.6 \\
\cline{2-11}
& Middle School & Status Quo & 34.6 & 26.6 & 36.9 & 43.3 & 44.7 & 46.0 & 48.1 & 50.9 \\
&  & Our Solution & 25.8 & 20.1 & 25.0 & 30.9 & 32.2 & 34.6 & 37.1 & 38.1 \\
&  & $\Delta$ & 8.8 & 6.5 & 11.9 & 12.4 & 12.5 & 11.4 & 11.0 & 12.8 \\
\cline{2-11}
& High School & Status Quo & 34.4 & 26.3 & 33.9 & 42.3 & 43.6 & 45.1 & 46.5 & 51.8 \\
&  & Our Solution & 26.6 & 21.5 & 25.8 & 31.4 & 32.9 & 35.9 & 37.8 & 40.6 \\
&  & $\Delta$ & 7.8 & 4.8 & 8.1 & 10.9 & 10.7 & 9.2 & 8.7 & 11.2 \\
\hline
District II & Elementary School & Status Quo & 42.7 & 27.2 & 41.1 & 54.5 & 56.6 & 61.3 & 69.1 & 78.8 \\
&  & Our Solution & 32.2 & 23.3 & 29.8 & 38.6 & 40.2 & 45.3 & 55.7 & 62.9 \\
&  & $\Delta$ & 10.5 & 3.9 & 11.3 & 15.9 & 16.4 & 16.0 & 13.4 & 15.9 \\
\cline{2-11}
& Middle School & Status Quo & 39.2 & 22.6 & 37.0 & 52.1 & 54.3 & 62.3 & 72.2 & 76.3 \\
&  & Our Solution & 28.1 & 19.9 & 27.0 & 34.5 & 37.2 & 39.5 & 41.5 & 45.9 \\
&  & $\Delta$ & 11.1 & 2.7 & 10.0 & 17.6 & 17.1 & 22.8 & 30.7 & 30.4 \\
\cline{2-11}
& High School & Status Quo & 42.7 & 27.2 & 40.5 & 54.6 & 58.6 & 61.7 & 65.5 & 71.7 \\
&  & Our Solution & 28.8 & 22.0 & 26.8 & 35.5 & 38.0 & 40.7 & 44.1 & 52.4 \\
&  & $\Delta$ & 13.9 & 5.2 & 13.7 & 19.1 & 20.6 & 21.0 & 21.4 & 19.3 \\
\hline
\end{tabular}
\end{table}

\begin{figure}[h!]
    \centering
    \begin{minipage}{0.48\textwidth}
        \centering
        \textbf{District I}\\
        \includegraphics[width=\linewidth]{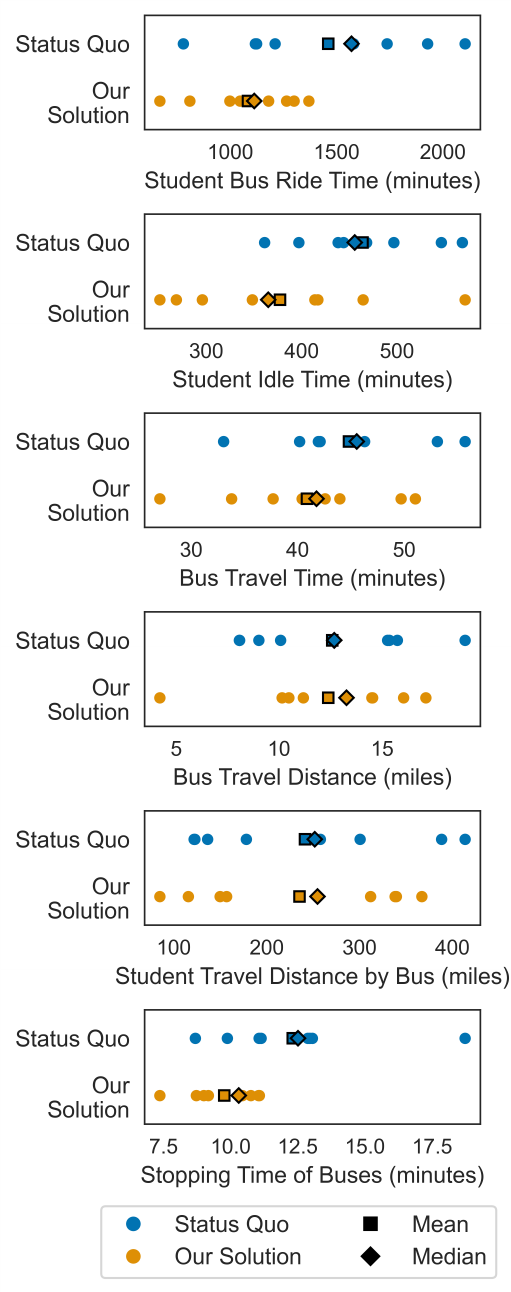} 
    \end{minipage}\hfill
    \begin{minipage}{0.48\textwidth}
        \centering
        \textbf{District II}\\
        \includegraphics[width=\linewidth]{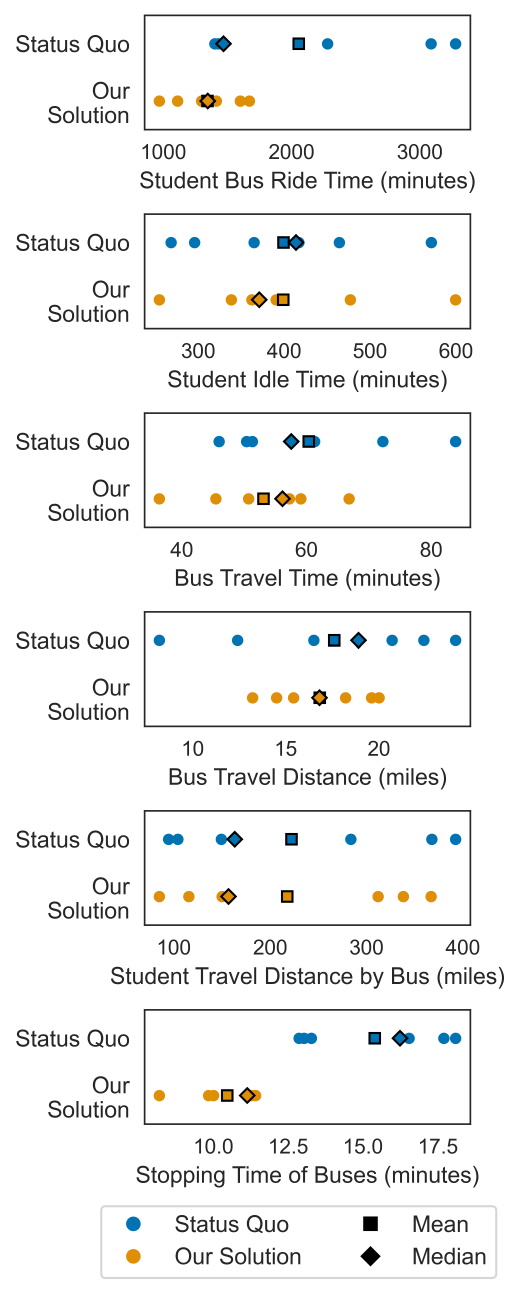} 
    \end{minipage}
    \caption{Structure of the resulting routes in the two districts. We show the student bus ride time (min), student idle time (min), bus travel time (min),  bus travel distance (miles), student travel distance by bus (miles), and stopping time of buses (min) for each bus under the status quo and our solution.}
    \label{fig:two_plots}
\end{figure}
\clearpage
% \begin{figure}[H]
%     \centering
% \includegraphics[width=0.75\textwidth]{sources/Images/dotplot_paper_updated_color_adjusted_annotated.pdf}
%     \caption{Structure of the resulting bus routes. The upper two figures plot the bus travel times, with one dot corresponding to each bus route. The lower six figures plot the maximum bus ride time across all students attending each particular school using each individual bus. Compared to the status-quo solution, our approach consistently reduces the average bus travel times and the longest bus ride times across the different schools in both school districts.}
%     \label{fig:1D-Dotplot}
% \end{figure}
\newpage

\section{Tradeoff Quantification Between {Average Commute Time} and Walking Time}\label{Walking_time_driving_time}
{
Table \ref{Table:Metrics_baseline} shows that our solution approach yields substantial reductions in the average bus ride times as well as the average commute times for students. These gains are accompanied by only a modest increase in walking time (0.8 min in District I and 0.96 min in District II). This increase is reasonable for the partner school districts, as discussed in Section \ref{sec3-Collaboration}, where current operations already impose a maximum walking distance of 0.3 miles—well below the state guideline of 1.2 miles. In this appendix, we treat walking time as a policy lever and examine how tightening or relaxing it impacts system performance. Specifically, we evaluate (i) whether the observed reductions in student commute time persist when walking times are tightly constrained to remain near status quo levels, and (ii) the extent to which additional improvements can be achieved when the maximum allowable walking distance constraint is relaxed.

\subsection{Constraining Walking Time to be Similar to Status Quo Solutions}\label{subses:constrain_WD}
To assess whether improvements persist under tighter walking-time control, we impose an explicit constraint on average walking time. Let $e^{SQ}_{isk}$ be the assignments according to the status quo routing solution. Then, the average walking time under the status quo, $\bar{F}^{\mathrm{SQ}}$, is defined as follows.
\begin{subequations}\renewcommand{\theequation}{\theparentequation.\arabic{equation}}
\begin{equation}
\bar{F}^{\mathrm{SQ}} = \frac{1}{|\mathcal{S}|} \sum_{i \in \mathcal{N}} \sum_{s \in \mathcal{S}_i} \sum_{k \in \mathcal{K}} F_{is} \cdot e^{SQ}_{isk}
\label{SQ_WT}
\end{equation}

Under our solution approach in Section \ref{sec5-SolutionApproach}, students are first assigned to buses (and in effect, service areas) in the clustering step. Then, for each bus, the routing problem is solved. To examine performance under stricter control of walking time, we augment the formulation in Appendix \ref{SingleBusMILO} (Equations \eqref{AppendixObjFun}--\eqref{eq:DoD:3_singleBus}) with an additional constraint that limits the average walking time of students assigned to a bus. Specifically, for each bus, we now force the average walking time of students in its service area to remain within a small tolerance \(\chi\) of the district-wide average walking time under the status quo. For our computational experiments, we set \(\chi = 20\) seconds. Note that we do not require this average to match the overall district-wide average status quo value exactly, because doing so would over-constrain the problem at the individual bus level. Moreover, reducing the maximum allowable walking distance below 0.3 miles would render some students ineligible for service by leaving them without any candidate stop within the allowable walking range. Accordingly, for each bus, we impose the following constraint.
\begin{equation}
\frac{1}{|\tilde{\mathcal{S}}|}\sum_{i \in \tilde{\mathcal{N}}}\sum_{s\in \tilde{S_{i}}} F_{is} \, e_{is} \le \bar{F}^{\mathrm{SQ}} + \chi  \label{SQ_WT_1}
\end{equation}

\end{subequations}

As shown in Figure \ref{fig:RideTimeWalkingTimeTradeOff}, even under this stricter walking-time constraint, our approach continues to achieve substantial reductions in average student commute time, yielding improvements in the total commute times of approximately 18\% (corresponding to an average of about 4–6 min per student) without materially changing average walking times.
\begin{figure}[H]
    \centering
\includegraphics[width=0.95\textwidth]{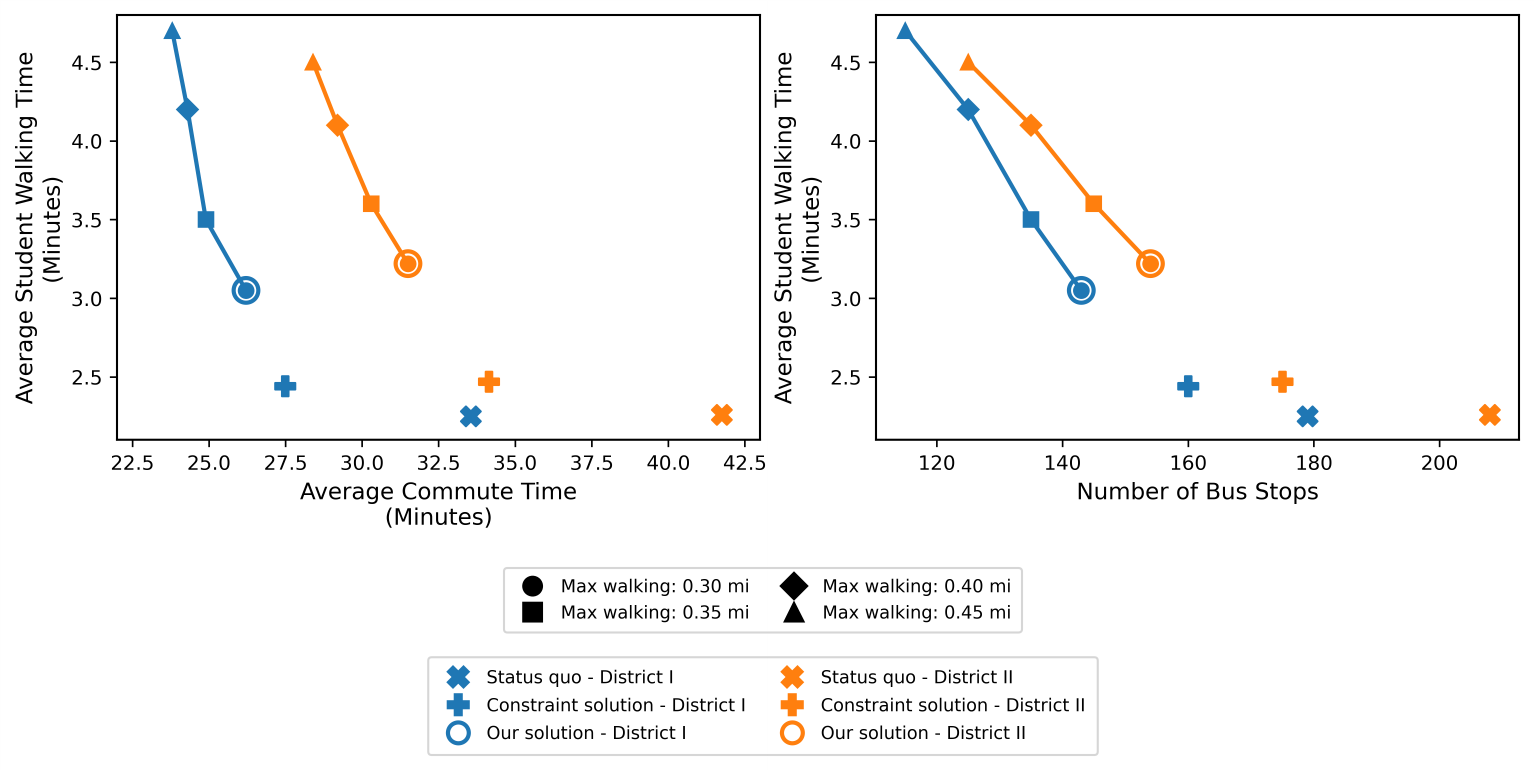}
\caption{ Impact of walking time (Y axis) on the average commute time and number of stops (X axis)}
\label{fig:RideTimeWalkingTimeTradeOff}
\end{figure}

\subsection{Effect of Increasing Maximum Allowable Walking Distance}
We next examine how relaxing the maximum allowable walking distance affects average student commute time. Increasing this limit from 0.30 to 0.45 miles provides additional flexibility in stop placement, reducing the number of stops (from 143 to 118 in District I and from 154 to 124 in District II) and enabling more direct routing.

As shown in Figure \ref{fig:RideTimeWalkingTimeTradeOff}, increasing the walking limit from 0.30 to 0.45 miles yields an additional reduction of approximately 3--4 min in average student commute time across both districts, while average walking time increases from about 3.0 to 4.5--4.7 min. The largest gains in terms of total commute time reductions occur between 0.30 and 0.35 miles, with diminishing improvements beyond 0.40 miles.

These findings, combined with the results in Section \ref{subses:constrain_WD}, highlight walking time as a controllable lever: substantial reductions in average commute time are achievable even under tight constraints on walking times, with further (though diminishing) gains available as these constraints are relaxed.}

\newpage
\section{Impact of Fleet Size on the {Total Commute Time}}\label{Impact_of_fleet_size}
\begin{figure}[h!]
    \centering
\includegraphics[scale=0.40]{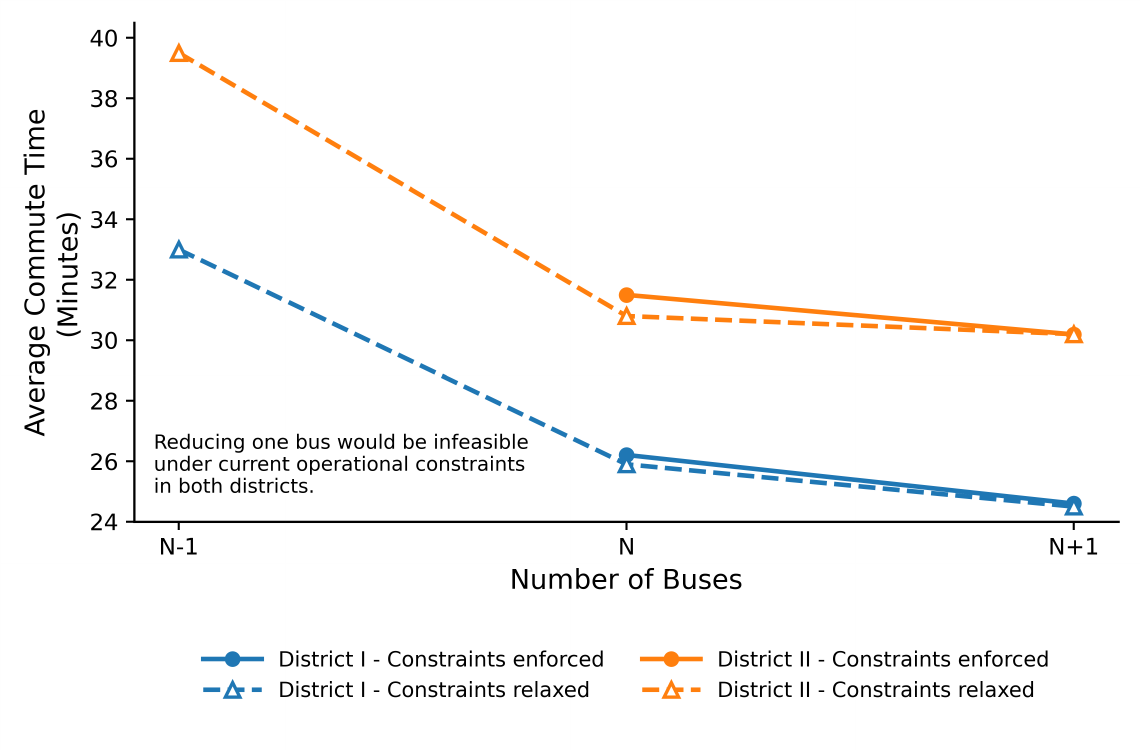}
\caption{Impact of fleet size on the average values of { students' total commute times}}
\label{fig:fleetsize}
\end{figure}
{As discussed in Section \ref{sec1-Introduction} and noted in Section \ref{sec3-Collaboration}, rural school districts typically operate with a small fleet of school buses. As a result, removing even a single bus can substantially increase student commute times and may render operations infeasible under the maximum student bus ride time constraint and the maximum bus travel time constraint. Conversely, adding a bus may reduce commute times, but districts must weigh these gains against the significant additional operating costs associated with fleet expansion. Motivated by this tradeoff, this appendix examines the impact of fleet size on average student commute time in both districts, with and without the maximum ride time and bus travel time constraints, as shown in Figure \ref{fig:fleetsize}. In Figure \ref{fig:fleetsize}, circular markers denote the results obtained with these two constraints enforced, while triangular markers denote the results obtained after relaxing them.

The results again show a clear asymmetry in how fleet size affects performance. In both districts, reducing the fleet by one bus leads to a sharp increase in average commute time when the maximum ride time and bus travel time constraints are relaxed. Under the constrained formulation, reducing the fleet by one bus is infeasible in both districts, indicating that the existing fleet sizes are already close to the minimum required to satisfy current operational requirements. This reinforces the practical difficulty of further fleet reductions in rural systems with limited slack.

In contrast, increasing the fleet by one bus yields only modest improvements. In District I, increasing the fleet from 9 to 10 buses reduces average commute time only slightly under both the constrained and relaxed cases. Similarly, in District II, increasing the fleet from 7 to 8 buses produces only a small reduction in average commute time, with almost no difference between the constrained and relaxed cases. This lack of meaningful difference between the constrained and relaxed cases is intuitive: once an additional bus is added, each bus serves fewer students, so the maximum ride time and bus travel time constraints are much less likely to be binding. As a result, relaxing these constraints provides little additional benefit in that setting.

Overall, the figure suggests that fleet reduction can impose disproportionate penalties and may even be infeasible, whereas fleet expansion offers only limited marginal gains. Taken together, these results indicate that the current fleet sizes in both districts are close to operationally efficient: reducing the fleet is either infeasible or leads to substantial increases in commute time, while adding a bus yields only modest improvements. Thus, for these districts, our results provide evidence against both fleet reduction and expansion under current operating conditions, suggesting that maintaining the existing fleet size is a well-balanced choice.}
\newpage

\section{MILO Formulation For a Single Bus} \label{SingleBusMILO}
In this section, we present the optimization formulation for the single school bus routing and scheduling problem. Section \ref{Subsec:Reduced_Routing_Optimization} employs a reduced version of this formulation (with the three reductions listed in Section \ref{Subsec:Reduced_Routing_Optimization}), while the full formulation is used in Section \ref{Subsec:WarmStart} in order to obtain the routing and scheduling solution for a single bus. Table \ref{Table:Notation_1bus} provides the notation used to define this MILO formulation for a single bus. This notation is a simplified version of the original notation used in Section \ref{sec4-Model} for the full formulation, but now adapted to the single bus setting.
\vspace{-0.2mm}
\renewcommand{\arraystretch}{0.99}
\begin{table}[H]
\caption{Notation for Single School Bus Routing and Scheduling Formulation}
\label{Table:Notation_1bus}
\scriptsize
\begin{tabular}{>{\raggedright}p{2.3cm} >{\raggedright}p{1cm} p{12.2cm}}
\toprule
\textbf{Component} & \textbf{Type} & \textbf{Description} \\
\midrule
%\{d\}$ & Set & Destination \\
$\tilde{\mathcal{N}}$ & Set & Nodes representing bus stops within the service region of the bus \\
$\mathcal{M}$ & Set & Nodes representing schools\\
$\tilde{\mathcal{S}}$ & Set & Students located within the service region of the bus\\
$\tilde{\mathcal{A}_1}$ & Set & Zero-cost directed arcs connecting origin $O$ to bus stops located within the service region of the bus. $\tilde{\mathcal{A}_1} = \{(O,i): \forall i \in \tilde{\mathcal{N}}\}$\\
$\tilde{\mathcal{A}_2}$ & Set & Directed arcs connecting bus stops to other bus stops located within the service region of the bus. $\tilde{\mathcal{A}_2} = \{(i,j), \forall i,j \in \tilde{\mathcal{N}}, i \neq j \}$\\
$\tilde{\mathcal{A}_3}$ & Set & Directed arcs connecting bus stops within the service region of the bus to schools. $\tilde{\mathcal{A}_3} = \{(i,j), 
    \forall i \in \tilde{\mathcal{N}}, \forall j \in \mathcal{M}\}$\\
$\tilde{\mathcal{A}_4}$ & Set & Directed arcs connecting schools to other schools. 
$\tilde{\mathcal{A}_4} = \{(i,j), 
    \forall i,j \in \mathcal{M}, i\neq j\}$\\
$\tilde{\mathcal{A}_5}$ & Set & Directed arcs connecting schools to bus stops within the service region of the bus. $\tilde{\mathcal{A}_5} = \{(i,j), \forall i \in \mathcal{M},  \forall j \in \tilde{\mathcal{N}}\}$\\
$\tilde{\mathcal{A}_6}$ & Set & Zero-cost directed arcs connecting schools to destination $D$. $\tilde{\mathcal{A}_6} = \{(i,D), \forall i \in \mathcal{M} \}$ \\
$\tilde{\mathcal{A}}$ & Set & Arcs. $\tilde{\mathcal{A}} = \tilde{\mathcal{A}_1} \cup \tilde{\mathcal{A}_2} \cup \tilde{\mathcal{A}_3} \cup \tilde{\mathcal{A}_4} \cup \tilde{\mathcal{A}_5} \cup \tilde{\mathcal{A}_6}$ \\
$\tilde{\mathcal{N}'}$ & Set & Nodes. $\tilde{\mathcal{N}}' = \{O\} \cup \tilde{\mathcal{N}} \cup \mathcal{M} \cup \{D\}$\\
$\tilde{\mathcal{S}_i}\subseteq \tilde{\mathcal{S}}$ & Set & Students located within the service region of the bus and within a short walking distance to stop $i\in\tilde{\mathcal{N}}$\\
$\tilde{\mathcal{S}_m}\subseteq \tilde{\mathcal{S}}$ & Set & Students located within service region of the bus that attend school $m\in\mathcal{M}$\\
\midrule
    $O, D \in \tilde{\mathcal{N}'}$ & Par. & Origin and destination depot nodes \\
    $A_1, A_2$ & Par. & Intercept of the linear expression for the time required to pick up and drop off students respectively\\
    $B_1, B_2$ & Par. & Slope of the linear expression for the time required to pick up or drop off students respectively\\
    $C$ & Par. & Capacity of the bus\\
    %$D$ & Par. & Maximum walking distance from home to bus stop\\
    % $T$ & Par. & Maximum travel time allowed for the bus route\\
    $G_{ms} \in\{0,1\}$ & Par. & 1 if student $s \in \tilde{S}$ attends school $m \in \mathcal{M}$, 0 otherwise\\
    $Q_i \in \{-1,0,1\}$ & Par. & -1 if $i=D$, 1 if $i=O$, 0 otherwise, $\forall i \in \tilde{\mathcal{N}'}$\\
    $\Delta_{ij} \in \mathbb{R}^+$ & Par. & Bus travel time for arc $(i,j) \in \tilde{\mathcal{A}}$\\
    $L_m \in \mathbb{Z}^+$ & Par. & Number of students attending school $m \in \mathcal{M}$\\
    {$T^B$} & {Par.} & {Maximum time allowed for the bus route}\\
    {$T^S$} & {Par.} & {Maximum allowable student time on the bus}\\
    {$F_{is}$} & {Par.} & {Walking time for student $s \in \tilde{\mathcal{S}_i}$ to stop $i \in \tilde{\mathcal{N}}$}\\
    % $M^{\prime} \in \mathbb{R}^+$ & Par. & Commute time between the farthest student that is likely to be picked while commuting in between schools and the school accessed by the student\\
    %$M^{\prime} \in \mathbb{R}^+$ & Par. & A large number\\
\midrule
    $x_{ij} \in\{0,1\}$ & Var. & 1 if the bus travels on arc $(i,j) \in \tilde{\mathcal{A}}$, 0 otherwise\\
    {$v_m \in \mathbb{Z}^+$} & Var. & Number of students dropped off by the bus at school $m \in \mathcal{M}$ \\
    % $u_{ik} \in \mathbb{R}^+$ & Var. & Number of students on bus $k \in \mathcal{K}$ right after the bus leaves after stopping at node $i \in \mathcal{N}'$\\
    $t_i \in \mathbb{R}^+$ & Var. & Time since departing from the first bus stop visited by the bus until when it reaches node $i \in \tilde{\mathcal{N}'}$\\
    $r_i \in\{0,1\} $ & Var. & 1 if the bus visits node $i \in \tilde{\mathcal{N}'}$, 0 otherwise\\
    ${w_{ij} \in \mathbb{Z}^+}$ & Var. & Number of students on the bus when traveling on arc $(i,j)\in \tilde{\mathcal{A}}$ if it travels on arc $(i,j)$, 0 otherwise\\
    $p_s \in \mathbb{R}^{+}$ & Var. &  Time when student $s \in \tilde{\mathcal{S}}$ is picked up by the bus\\
    $d_s \in \mathbb{R}^{+}$ & Var. &  Time when student $s \in \tilde{\mathcal{S}}$ is dropped off by the bus\\  
    $e_{is} \in\{0,1\} $ & Var. & 1 if student $s \in \tilde{\mathcal{S}_i}$ is picked up by the bus at bus stop $i \in \tilde{\mathcal{N}}$, 0 otherwise\\
    $\tau_{is} \in \mathbb{R}^+$  & Var. & Time at which student $s \in \tilde{\mathcal{S}_i}$ is picked up by the bus at bus stop $i \in \tilde{\mathcal{N}}$ if it picks up that student at that bus stop, 0 otherwise\\
    $\kappa_{ms} \in \mathbb{R}^+$  & Var. & Time at which student $s \in \tilde{\mathcal{S}_m}$ is dropped off by the bus at school $m \in \mathcal{M}$ if it drops off that student at that school, 0 otherwise\\
\bottomrule
\end{tabular}
\end{table}

\begin{subequations}
\renewcommand{\theequation}{\theparentequation.\arabic{equation}}
{\begin{equation}\label{AppendixObjFun}
    \text{Maximize} \sum_{s \in \tilde{\mathcal{S}}} p_s - \sum_{i \in \tilde{\mathcal{N}}} \sum_{s \in \tilde{\mathcal{S}_i}} F_{is} \cdot e_{is}
\end{equation}}
\begin{align}
% &\sum_{i \in \mathcal{N}} x_{oik} = r_{ok} \quad \forall k \in \mathcal{K} \label{ODConstraint:EachBusCanHaveOnlyOneStartPoint}\\
% &\sum_{m \in \mathcal{M}} x_{mdk} = r_{dk} \quad \forall k \in \mathcal{K} \label{ODConstraint:EachBusEndsAtASchool}\\
% &t_{i} \leq T\cdot (1-x_{Oi}) \quad \forall i \in \tilde{\mathcal{N}}\label{ODConstraint:StartTime_singleBus}\\
% &t_{D} \leq T \label{ODConstraint:TravelTimeLimit_SingleBus}\\
&r_{O} = 1,  r_{D} =1, t_D = {T^{B}} \quad \label{ODConstraint:DepotStartTime_SingleBus}\\
&\sum_{i \in \tilde{\mathcal{N}}}\sum_{s \in \tilde{\mathcal{S}_i} \cap \tilde{\mathcal{S}_m}} e_{is} = v_{m} \quad \forall m \in \mathcal{M}\label{PickUpConstraintNew:All_students_dropped_off_by_bus_are_picled_up_SingleBus}\\
& v_{m} = L_m \quad \forall m \in \mathcal{M}\label{PickUpConstraint:Students_Droppedoff_AtSchool_NumberOfStudentsThatGoToSchool_singleBus}\\
&\sum_{j\in \tilde{\mathcal{N}} \cup \mathcal{M} \cup \{D\}:(i,j)\in \tilde{\mathcal{A}}} x_{ij} - \sum_{j\in \tilde{\mathcal{N}} \cup \mathcal{M} \cup \{O\}:(j,i)\in \tilde{\mathcal{A}}} x_{ji} = Q_i \quad \forall i \in \tilde{\mathcal{N}}^\prime \label{FBConstraint:FlowBalanceConstraint_singleBus}\\
&r_{i} + r_{j} \geq 2 \cdot x_{ij} \quad  \forall (i,j) \in \tilde{\mathcal{A}}\label{BR:Constraint-Consecutive_stops_singleBus}\\
&\sum_{j \in \tilde{\mathcal{N}}^\prime : (i,j) \in \tilde{\mathcal{A}}} x_{ij} = r_{i} \quad \forall i \in \tilde{\mathcal{N}} \cup \mathcal{M}\label{BR:Constarint-Bus_has_to_come_from_somewhere_to_pick_up_or_drop-off_singleBus}\\
&r_{i} \leq \sum_{s \in \tilde{\mathcal{S}_{i}}} e_{is} \quad \forall i \in \tilde{\mathcal{N}}\label{BR:Constarint-Bus_stops_students_have_to_be_picked_singleBus}\\
&e_{is} \leq r_{i} \quad \forall i \in \tilde{\mathcal{N}}, \forall s \in \tilde{\mathcal{S}_i},\label{BR:Constarint-Bus_stops_students_have_to_be_picked_1_singleBus}\\
&\sum_{l\in \tilde{\mathcal{N}'}:(l,i) \in \tilde{\mathcal{A}_1} \cup \tilde{\mathcal{A}_2} \cup \tilde{\mathcal{A}_5}} w_{li}+ \sum_{s \in \tilde{\mathcal{S}_i}} e_{is} - w_{ij} \leq C\cdot (1-x_{ij}) \quad \forall (i,j) \in \tilde{\mathcal{A}_2} \cup \tilde{\mathcal{A}_3} \label{BR:Constraint-Bus_visits_stop_students_before_after_positive_singleBus}\\
&-\left(\sum_{l\in\tilde{\mathcal{N}'}:(l,i) \in \tilde{\mathcal{A}_1} \cup \tilde{\mathcal{A}_2} \cup \tilde{\mathcal{A}_5}} w_{li}+ \sum_{s \in \tilde{\mathcal{S}_i}} e_{is} - w_{ij}\right) \leq C\cdot (1-x_{ij}) \quad \forall (i,j) \in \tilde{\mathcal{A}_2} \cup \tilde{\mathcal{A}_3}\label{BR:Constraint-Bus_visits_stop_students_before_after_negative_singleBus}\\
&\sum_{l\in \tilde{\mathcal{N}'}:(l,i) \in \tilde{\mathcal{A}_3} \cup \tilde{\mathcal{A}_4}} w_{li} - v_{i} - w_{ij} \leq C\cdot (1-x_{ij}) \quad \forall (i,j) \in \tilde{\mathcal{A}_4} \cup \tilde{\mathcal{A}_5} \cup \tilde{\mathcal{A}_6}\label{BR:Constraint-Bus_visits_school_students_before_after_positive_singleBus}\\
&-\left(\sum_{l\in \tilde{\mathcal{N}'}:(l,i) \in \tilde{\mathcal{A}_3} \cup \tilde{\mathcal{A}_4}} w_{li} - v_{i} - w_{ij}\right) \leq C\cdot (1-x_{ij}) \quad \forall (i,j) \in \tilde{\mathcal{A}_4} \cup \tilde{\mathcal{A}_5} \cup \tilde{\mathcal{A}_6}\label{BR:Constraint-Bus_visits_school_students_before_after_negative_singleBus}\\
&G_{ms} + \sum_{i\in \tilde{\mathcal{N}}: s\in \tilde{\mathcal{S}_i}} e_{is} \leq r_{m}+1 \quad \forall m \in \mathcal{M}, s \in \tilde{\mathcal{S}}\label{BR:Constraint-Bus_visits_school_if_student_is_on_board_singleBus}\\
&\sum_{i \in \tilde{\mathcal{N}}: s \in \tilde{\mathcal{S}_i}} e_{is} = 1 \quad \forall s \in \tilde{\mathcal{S}}\label{BR:Student_unique_bus_unique_stop_singleBus}\\
&t_{i} + \Delta_{ij} + A_1 + B_1 \cdot \sum_{s \in \tilde{\mathcal{S}_i}}e_{is} - t_{j} \leq {T^B} \cdot (1-x_{ij})\quad  \forall (i,j) \in \tilde{\mathcal{A}_2} \cup \tilde{\mathcal{A}_3}\label{BR:Constraint-pick_up_time_singleBus}\\
&t_{i} + \Delta_{ij} + A_2 + B_2 \cdot v_{i} - t_{j} \leq {T^B} \cdot (1-x_{ij})\quad  \forall (i,j) \in \tilde{\mathcal{A}_4} \cup \tilde{\mathcal{A}_5}  \cup \tilde{\mathcal{A}_6}\label{BR:Constraint-drop_off_time_singleBus}\\
&{d_{s} - p_{s} \leq T^S \quad \forall s \in \tilde{S}}\\
& w_{ij} \leq C \cdot x_{ij} \quad \forall (i,j) \in \tilde{\mathcal{A}}\label{CapacityConstraint:NumberofStudentsInBusConsecutiveStops1_singleBus}\\
&\tau_{is} + {T^B} \cdot (1 - e_{is}) \geq p_{s} \quad \forall i \in \tilde{\mathcal{N}}, s \in \tilde{\mathcal{S}_i}\label{StudentTimeConstraintNew:PickUpTimeGreater_singleBus}\\
&\kappa_{ms} \leq d_{s} \quad \forall m \in \mathcal{M}, s \in \tilde{\mathcal{S}_m}\label{StudentTimeConstraint:Drop-offTimeLesser_singleBus}\\
% &T \cdot \sum_{i \in \tilde{\mathcal{N}}: s \in \tilde{\mathcal{S}_i}}e_{is} \geq {p_{s}} \quad \forall s \in \tilde{\mathcal{S}}\label{StudentTimeConstraintNew:PickUpTimeUpperLimit_singleBus}\\
% &T \cdot \sum_{i \in \tilde{\mathcal{N}}: s \in \tilde{\mathcal{S}_i}} e_{is} \geq {d_{s}} \quad \forall s \in \tilde{\mathcal{S}}\label{StudentTimeConstraintNew:Drop-offTimeUpperLimit_singleBus}\\
&\tau_{is} \leq {T^B} \cdot e_{is} \quad \forall i \in \tilde{\mathcal{N}}, s \in \tilde{\mathcal{S}_i}\label{StudentTimeConstraintNew:tau0=>t0_singleBus}\\
&\tau_{is} \leq t_{i}\quad \forall i \in \tilde{\mathcal{N}}, s \in \tilde{\mathcal{S}_i}\label{StudentTimeConstraint:TauCannotBeGreaterthanT_singleBus}\\
&\tau_{is} \geq t_{i} - {T^B} \cdot (1 - e_{is}) \quad \forall i \in \tilde{\mathcal{N}}, s \in \tilde{\mathcal{S}_i}\label{StudentTimeConstraintNew:Forcet=tau->p=1_singleBus}\\
&\kappa_{ms} \leq {T^B} \cdot \sum_{i \in \tilde{\mathcal{N}}: s \in \tilde{\mathcal{S}_i}} e_{is} \quad \forall m \in \mathcal{M}, s \in \tilde{\mathcal{S}_m}\label{StudentTimeConstraintNew:kappa0=>t0_singleBus}\\
&\kappa_{ms} \leq t_{m}\quad \forall m \in \mathcal{M}, s \in \tilde{\mathcal{S}_m}\label{StudentTimeConstraint:KappaCannotBeGreaterthanT_singleBus}\\
&\kappa_{ms} \geq t_{m} - {T^B} \cdot (1 - \sum_{i \in \tilde{\mathcal{N}}: s \in \tilde{\mathcal{S}_i}} e_{is}) \quad \forall m \in \mathcal{M}, s \in \tilde{\mathcal{S}_m}\label{StudentTimeConstraintNew:Forcet=Kappa->q=1_singleBus}\\
&{T^B}\cdot(G_{ms}+e_{is} – 2) \leq (t_{m} – t_{i}) \quad \forall m \in \mathcal{M}, i \in \tilde{\mathcal{N}}, s \in \tilde{\mathcal{S}_i}\label{StudentTimeConstraintNew:pickupafterdropoff_singleBus}\\
&e_{is} \in \{0,1\}, \tau_{is} \in \mathbb{R}^+ \quad \forall s \in \tilde{\mathcal{S}_i}, i \in \tilde{\mathcal{N}} \label{eq:DoD:1_singleBus}\\
& r_{i} \in \{0,1\}, {t_{i} \in \mathbb{R}^+} \quad \forall i \in \tilde{\mathcal{N}'}\label{eq:DoD:4_singleBus}\\
&x_{ij} \in \{0,1\}, {w_{ij} \in \mathbb{Z}^+} \quad \forall (i, j) \in \tilde{\mathcal{A}} \label{eq:DoD:2_singleBus}\\
&{\kappa_{ms}}, p_{s}, d_{s} \in \mathbb{R}^+, {v_{m} \in \mathbb{Z}^+ } \quad \forall m \in \mathcal{M}, s \in {\tilde{\mathcal{S}_m}}\label{eq:DoD:3_singleBus}
% &v_{m} \in \mathbb{Z}^+ \quad \forall m \in \mathcal{M}, s \in \tilde{\mathcal{S}_m}\label{eq:DoD:5_singleBus}
\end{align}
\end{subequations}